\documentclass[UTF-8,reqno]{amsart}
\usepackage{enumerate}
\usepackage{mhequ}
\usepackage{amssymb,url,color, booktabs,nccmath}
\usepackage[left=3.2cm,right=3.2cm,top=4cm,bottom=4cm]{geometry}
\usepackage{mathrsfs}
\usepackage{enumitem,dsfont}

\usepackage{color}
\usepackage[colorlinks=true]{hyperref}
\hypersetup{
    linkcolor=blue,          
    citecolor=red,        
    filecolor=blue,      
    urlcolor=cyan
}

\definecolor{darkergreen}{rgb}{0.0, 0.5, 0.0}

\numberwithin{equation}{section}
\def\theequation{\arabic{section}.\arabic{equation}}
\newcommand{\be}{\begin{eqnarray}}
\newcommand{\ee}{\end{eqnarray}}
\newcommand{\ce}{\begin{eqnarray*}}
\newcommand{\de}{\end{eqnarray*}}
\newtheorem{theorem}{Theorem}[section]
\newtheorem{lemma}[theorem]{Lemma}
\newtheorem{proposition}[theorem]{Proposition}
\newtheorem{Examples}[theorem]{Example}
\newtheorem{corollary}[theorem]{Corollary}

\newtheorem{definition}[theorem]{Definition}
\theoremstyle{definition}
\newtheorem{remark}[theorem]{Remark}

\DeclareMathOperator{\supp}{supp}

\def\[{{\Big[}}
\def\]{{\Big]}}
\def\<{{\langle}}
\def\>{{\rangle}}
\def\({{\Big(}}
\def\){{\Big)}}

\def\bx{{\mathbf{x}}}
\def\tr{\mathrm {tr}}

\def\dif{{\mathord{{\rm d}}}}

\def\={&\!\!=\!\!&}

\def\1{{\mathbf{1}}}

\def\geq{\geqslant}
\def\leq{\leqslant}

\def\div{\mathord{{\rm div}}}

\def\[{{\Big[}}
\def\]{{\Big]}}
\def\<{{\langle}}
\def\>{{\rangle}}
\def\({{\Big(}}
\def\){{\Big)}}

\def\bx{{\mathbf{x}}}
\def\tr{\mathrm {tr}}

\def\dif{{\mathord{{\rm d}}}}

\def\={&\!\!=\!\!&}
\def\bt{\begin{theorem}}
\def\et{\end{theorem}}
\def\bl{\begin{lemma}}
\def\el{\end{lemma}}
\def\br{\begin{remark}}
\def\er{\end{remark}}
\def\bx{\begin{Examples}}
\def\ex{\end{Examples}}
\def\bd{\begin{definition}}
\def\ed{\end{definition}}
\def\bp{\begin{proposition}}
\def\ep{\end{proposition}}
\def\bc{\begin{corollary}}
\def\ec{\end{corollary}}

\def\geq{\geqslant}
\def\leq{\leqslant}

\def\div{\mathord{{\rm div}}}

\def\Id{\textrm{Id}}

 \def\R{\mathbb R}
 \def\R{\mathbb R}    
\def\N{\mathbb N}  
   
\def\<{\langle} \def\>{\rangle}

\allowdisplaybreaks

\begin{document}

\title[Nonuniqueness for 2D NSE with Space-Time White Noise]{ Sharp Non-uniqueness of Solutions to 2D Navier-Stokes Equations with Space-Time White Noise}

\author{Huaxiang L\"u}
\address[H. L\"u]{Academy of Mathematics and Systems Science,
Chinese Academy of Sciences, Beijing 100190, China}
\email{lvhuaxiang22@mails.ucas.ac.cn }

\author{Xiangchan Zhu}
\address[X. Zhu]{ Academy of Mathematics and Systems Science,
Chinese Academy of Sciences, Beijing 100190, China}
\email{zhuxiangchan@126.com}

\thanks{
Research  supported   by National Key R\&D Program of China (No. 2022YFA1006300, 2020YFA0712700) and the NSFC (No.  12090014, 12288201) and
  the support by key Lab of Random Complex Structures and Data Science,
 Youth Innovation Promotion Association (2020003), Chinese Academy of Science. The financial support by the DFG through the CRC 1283 "Taming uncertainty and profiting
 from randomness and low regularity in analysis, stochastics and their applications" is greatly acknowledged.
}

\begin{abstract}
In this paper we are concerned with the 2D incompressible Navier-Stokes equations driven by space-time white noise. We establish existence of infinitely many global-in-time probabilistically strong
and analytically weak solutions $u$ for every divergence free initial condition $u_0\in L^p\cup C^{-1+\delta},\ p\in(1,2),\delta>0$. More precisely, there exist infinitely many solutions such that $u-z\in C([0,\infty);L^p)\cap L^2_{\rm{loc}}([0,\infty);H^\zeta)\cap L^1_{\rm{loc}}([0,\infty);W^{\frac13,1})$ for some $\zeta\in(0,1)$, where $z$ is the  solution to the linear equation. This result in
particular implies non-uniqueness in law. Our result is sharp in the sense that the solution satisfying $u-z\in C([0,\infty);L^2)\cap L^2_{\rm{loc}}([0,\infty);H^\zeta)$ for some $\zeta\in(0,1)$ is unique.
\end{abstract}


\date{\today}

\maketitle

\tableofcontents
\section{Introduction}

In this paper, we consider the following two dimensional Navier–Stokes system on $\mathbb{T}^2=\mathbb{R}^2/\mathbb{Z}^2$ driven by a space-time white noise
\begin{align}\label{1.1}
\dif u+\div(u\otimes u)\dif t+\nabla p\dif t&=\Delta u\dif t+\dif B_t,\notag\\
\div u&=0,\notag\\
    u(0)&=u_0,
\end{align}
where $p$ is the associated pressure. Here $B$ is a cylindrical Wiener process on some stochastic basis $(\Omega,\mathcal{F}, (\mathcal{F}_t)_{t\geq0},\mathbf{P})$. The time derivative of $B$ is the delta correlated space-time white noise. Such a noise appears in a  scaling limit of point vortex approximation and the vorticity form of the 2D Euler equations perturbed by a certain transport type noise (cf. \cite{FL20,FL21,LZ21}). More precisely, the scaling limit is described by the vorticity form of the 2D Navier-Stokes system driven by the curl of a space-time white noise, which is equivalent to the 2D Navier-Stokes equations driven by a space-time white noise in terms of velocity-pressure variables.

Specifically, it can be shown that under parabolic scaling the space-time white noise in spatial dimension $d$ is a random distribution of space-time regularity $-(d + 2)/2 - \kappa$ for any $\kappa>0$. According to Schauder’s estimates, we expect that solutions will have at most $-d/2+1-\kappa$ regularity. Therefore, in two dimensions, solutions are already not functions. As a result, the product in the convective term is analytically undefined, and probabilistic arguments are necessary to make sense of the equations.

Da Prato and Debussche \cite{DPD02} initially  solved this problem locally in time by decomposing \eqref{1.1} into a linear equation and nonlinear equation (see \eqref{3.1} and \eqref{3.2} below). Moreover, by utilizing the properties of the Gaussian invariant measure, they were able to obtain global-in-time existence for almost every initial condition with respect to the invariant measure. Through the strong Feller property in \cite{ZZ17}, global-in-time existence for every initial condition could be derived.

Recently, Hairer and Rosati \cite{HR23} employed dynamic high-low frequency decomposition and paraproduct to establish the existence of a unique global-in-time solution for initial conditions in $L^2\cup C^{-1+\delta}$, without relying on any knowledge of invariant measures. However, their methods depend on the solution to the linear equation being a log-correlated field and it is not clear that whether their method are applicable to even slightly more irregular noise.

In a more singular case, such as the 3D case, Hofmanov\'{a}, Zhu, and the second named author of this paper \cite{HZZ21b} established global existence and non-uniqueness of strong solutions in a paracontrolled sense, with the help of the convex integration method. The first goal of this paper is to extend the result in \cite{HZZ21b} to the two-dimensional case. However, the convex integration method in the 2D case is more complex than in the 3D case since there is less freedom in dimensions, and a space-time intermittent variant of the framework is required (see \cite{CL20,CL22}).

Additionally, the nonlinearity in the Navier-Stokes system appears similar to the one in the Langevin dynamics for the Yang-Mills measure, which is the stochastic quantization of the Yang-Mills field. Local-in-time solutions for this case were constructed in \cite{She21, CCHS22a,CCHS22b}. However, the existence of global solutions is still an open question. The idea is to use dynamics and PDE techniques to study properties of the field. Formally, these equations have the law of the associated field as an invariant measure. In the case of stochastic quantization of the Euclidean field theory, it was possible to use dynamics to construct and study properties of the corresponding measure (see \cite{MW17a,MW17b,GH21} and \cite{SZZ21}). As stochastic quantization of the 2D Yang-Mills field is much simpler than the 3D case, and global existence for the Langevin dynamics for the 2D Yang-Mills measure is unattainable using classical PDE techniques, we hope that our extension to the 2D case can provide some insight into this problem.

It is also natural to investigate the sharpness of uniqueness/non-uniqueness in this setting. In the deterministic case, the classical Ladyzhenskaya-Prodi-Serrin criterion ensures that there exists at most one solution in the class $C_tL^2$, without any additional regularity assumptions \cite{FLRT00, LM01}. However, if the solution $u$ to \eqref{1.1} minus the solution to the linear equation has no further regularity, the nonlinear term cannot be well-defined. Although \cite{HR23} provides a unique solution space for \eqref{1.1}, it is interesting to determine whether such a space is sharp, i.e., whether we could find more solutions in a slightly larger space. Therefore, the second aim of this paper is to provide an answer to this problem. Our results are twofold:
\begin{itemize}
  \item We prove that there exists at most one  solution $u$ satisfying  $u-z\in C([0,\infty);L^{2})\cap L^2_{\textrm{loc}}([0,\infty);H^\zeta)$ for some $\zeta\in(0,1)$, which can be viewed as Ladyzhenskaya-Prodi-Serrin criteria in this setting. Here $z$ is defined in (\ref{3.1*}).
\item The above result is sharp in the sense that there exists infinitely many solutions in a slightly larger space. More precisely, for every divergence free initial condition in $L^p,p\in (1,2)$ we show existence of infinite many global in time probabilistically strong and analytically weak solutions $u$ such that\begin{align*}
    u-z\in C([0,\infty);L^{p})\cap L^2_{\rm{loc}}([0,\infty);H^\zeta)
\end{align*} $\mathbf{P}$-a.s. for some $\zeta\in(0,1)$. We also emphasize that our method  also applies to more rough noise.
    \end{itemize}

\subsection{Main result}
Using Da Prato-Debussche's trick, we divide the equation (\ref{1.1}) into two parts: the following linear equation:
\begin{align}\label{3.1}
    \dif z-\Delta z \dif t+\nabla p_z&=\dif B,\notag\\
    \div z&=0,\notag\\
    z(0)&=0,
\end{align}
and the following random PDE:
\begin{align}\label{3.2}
    \dif v-\Delta v\dif t+\div((v+z)\otimes(v+z))\dif t+\nabla p_v\dif t&=0,\notag\\
    \div v&=0,\\
   v(0)&=u_0.\notag
\end{align}
The solution to (\ref{3.1}) can be written as
\begin{align}z(t):=\int_{0}^te^{(t-s)\Delta}\dif B_s,\label{3.1*}\end{align}where $ e^{t\Delta}$ is the heat semigroup. Note that here we do not choose $z$ being stationary, since the zero initial condition will be convenient to deal with the $L^p$ initial condition. Similar to \cite{DPD02}, Gaussian computation guarantees we still have $z\in C_TC^{-\kappa}$ for $\kappa>0$ and give a meaning to $\div(z\otimes z)$ by wick product. More details can be seen in Section \ref{so}.

Our first main result concerns the existence of many global solutions for any initial condition in $L^p\cup C^{-1+\delta},$ with $p\in(1,2),\delta>0$.

\bt\label{thm:globalu}
Let $u_0 \in L^p\cup C^{-1+\delta}$ with some  $p\in(1,2),\delta>0$ $\mathbf{P}$-a.s. be a divergence free initial condition independent of the Wiener process $B$. Then there exist infinitely many probabilistically strong and analytically weak solutions $u$ to the equation (\ref{1.1}) on $[0,\infty)$.

Moreover, if $u_0 \in L^p$ with some  $p\in(1,2)$, then for any $0<\kappa<(1-\frac{1}{p})\wedge(\frac2p-1)\leq\frac{1}{3}$ it holds that
\begin{align*}
z&\in C_{\rm{loc}}([0,\infty);C^{-\kappa}),\\
    u-z&\in C([0,\infty);L^{p})\cap L^2_{\rm{loc}}([0,\infty);H^\zeta)\cap L^1_{\rm{loc}}([0,\infty);W^{\frac13,1})
\end{align*} $\mathbf{P}$-a.s. for some $\zeta\in(0,1)$ independent of $\kappa$, where $z$ is defined in (\ref{3.1*}).
\et

\br
We could also cover the case introduced in  \cite{HR23}
\begin{align*}
  \dif u+\div(u\otimes u)\dif t+\nabla p\dif t&=\Delta u\dif t+\mathcal{Z}\dif t+\dif B,\\
\div u&=0,\ \ u(0)=u_0,
\end{align*}
where $\mathcal{Z}$ is in the parabolically scaled H\"older-Besov space $\mathcal{C}_{\rm{parab}}^{-2+3\kappa}$. In this case, we could set $z=\int_0^te^{(t-s)\Delta}\dif B_s+\int_0^te^{(t-s)\Delta}\mathcal{Z}\dif s$ and then do the same decomposition. Here by the regularity of $\mathcal{Z}$ and the Schauder estimate $\int_0^te^{(t-s)\Delta}\mathcal{Z}\dif s\in C_TC^{2\kappa}$. Then we define
 $$z^{:2:}(t):=:(\int_0^te^{(t-s)\Delta}\dif B_s)^2:+2\int_0^te^{(t-s)\Delta}\dif B_s\int_0^te^{(t-s)\Delta}\mathcal{Z}\dif s+(\int_0^te^{(t-s)\Delta}\mathcal{Z}\dif s)^2.$$
 Then $z\in C_TC^{-\kappa}$ and $z^{:2:}\in C_TC^{-2\kappa}$, where we used renormalization of $:(\int_0^te^{(t-s)\Delta}\dif B_s)^2:$ in $C_TC^{-2\kappa}$.
 \er

Furthermore, our construction directly implies the following result.
\bc\label{thm:noninlaw}
Non-uniqueness in law holds for the Navier–Stokes system (\ref{1.1}) for every given initial law supported on divergence free vector fields in $L^p\cup C^{-1+\delta}$ for $p\in(1,2),\delta>0$.
\ec

In the present paper, our idea is to utilize convex integration method in order
to construct global-in-time solutions. This iterative technique allows us to construct solutions step by step, with each iteration explicitly addressing a specific scale. This method heavily relies on the structure of the nonlinearity, which propagates oscillations and reduces the Reynolds stress error term as we progress through the iteration, gradually bringing us closer to a solution. In this study, we employ the accelerating jets technique proposed in \cite{CL20} for the two-dimensional case. Compared to the spatially intermittent framework used in \cite{HZZ21b} that was based on \cite{BV19}, the accelerating jets with space-time intermittent variant allow us more flexibility in terms of the non-uniqueness range scaling. The accelerating jets used time periodicity and in the stochastic setting we have to modify some estimate (see
Remark \ref{rem:gxi} and Theorem \ref{ihiot}) to adjust the stopping time. For more information on convex integration, please refer to section \ref{sec:ci}.

To enable convex integration in the singular setting of equation \eqref{3.2}, we adopt a further decomposition of the Navier-Stokes system \eqref{3.2}, following the approach in \cite{HZZ21b}. Specifically, we split $v$ into $e^{t\Delta}u_0+v^1+v^2$, where $v^1$ represents the irregular part and $v^2$ represents the regular part. The equation for $v^1$ is linear, while the equation for $v^2$ contains the quadratic nonlinearity. The equation for $v^1$ can be solved using a fixed point argument, and we apply convex integration on the level of $v^2$. Since the equation for $v^2$ is coupled with the equation for $v^1$, we propose a joint iterative procedure that approximates both equations. Compared to \cite{HZZ21b}, we separate the initial value part as our initial value is only in the $L^p$ space. This also leads to a different extension for the initial value part with $e^{|t|\Delta}u_0$ to negative time for the mollification around $t=0$, which is required in the convex integration estimate. Additionally, in \cite{HZZ21b}, to handle the singularity of $v^1$ around zero, the authors put some regular terms into the equation of $v^1$. In this paper, as we use the accelerating jets introduced in \cite{CL20} and only control the $L^1_TL^1$ norm of the Reynolds stress, we could keep these regular terms in the equation of $v^2$. To maintain the same initial value during the iteration, the oscillation can only be added for positive times. In our setting, we need to take care the term that arises from the noise part in a small time interval near zero (see Remark \ref{rem:3.7}).

Finally we obtain the following sharp uniqueness result in the space $ C([0,\infty);L^2)\cap L^2_{\textrm{loc}}([0,\infty);H^\zeta)$ which can be viewed as Ladyzhenskaya-Prodi-Serrin criteria in this setting.

\bt\label{thm:unique}
There exists at most one probabilistically strong and analytically weak solution $u$ to the equation (\ref{1.1}) on $[0,\infty)$ such that $u-z\in C([0,\infty);L^2)\cap L^2_{\rm{loc}}([0,\infty);H^\zeta)$ $\mathbf{P}$-a.s. for some $\zeta\in(0,1)$, where $z$ is defined in (\ref{3.1*}).
\et

\br If $u_0\in L^\infty$ there exists exactly one solution $u$ such that $u-z\in C([0,\infty);L^2)\cap L^2_{\rm{loc}}([0,\infty);H^\zeta)$ by \cite{HR23}.
\er

\subsection{Convex integration}\label{sec:ci}
Convex integration was first introduced to fluid dynamics by De Lellis and Sz\'ekelyhidi Jr. \cite{DLS09, DLS10, DLS13}. Since then, this method has led to several breakthroughs in the field, including the proof of Onsager's conjecture for the incompressible Euler equations \cite{Ise18,BDLSV19}. Also the question of
well/ill-posedness of the three dimensional Navier-Stokes equations has experienced an immense
breakthrough: Buckmaster and Vicol \cite{BV19b} established non-uniqueness of weak solutions with
finite kinetic energy, Buckmaster, Colombo and Vicol \cite{BCV18} were able to connect two arbitrary
strong solutions via a weak solution.
Sharp non-uniqueness results for the Navier-Stokes equations in dimension $d \geq 2$ were obtained by
Cheskidov and Luo \cite{CL20, CL22}. We refer to the reviews \cite{BV19, BV21} for further details and
references.

The convex integration method has been successfully applied in the stochastic setting as well. In \cite{HZZ19, HLP22, Pap23}, non-uniqueness in law of weak solutions to the Navier-Stokes equation driven by additive, linear multiplicative, transport noise, and nonlinear noise of cylindrical type was established. In \cite{HZZ20}, similar results were obtained, along with the existence and non-uniqueness of strong Markov solutions for the Euler equations. The long-standing problem of constructing probabilistically strong global solutions to the Navier-Stokes equation perturbed by trace-class noise, along with non-uniqueness statements for such solutions, was solved in \cite{ HZZ,CDZ22}. By combining convex integration techniques with tools from paracontrolled calculus \cite{GIP15} for singular SPDEs, \cite{HZZ21b} was able to study the Navier-Stokes equation perturbed by space-time white noise, yielding global existence and non-uniqueness of weak solutions. Remarkably, the interplay of techniques allows pushing the solution theory even further into the regimes of so-called supercritical equations \cite{HZZ22c} that are inaccessible by standard theories for singular SPDEs such as paracontrolled calculus or regularity structures \cite{Hai14}. Finally, in the most recent works \cite{CDZ22, HZZ22}, a stochastic convex integration theory is developed, leading to non-uniqueness of stationary ergodic solutions to stochastic perturbations of the Navier-Stokes and Euler equations.

Convex integration has also been applied to other related equations with weaker diffusion, such as power-law fluids with a small parameter $p$ (see \cite{BMS}, and \cite{LZ22} for the stochastic setting), and hypodissipative type Navier-Stokes equations (see \cite{CDLDR18, RS21, Yam21a, Yam21b, Yam22a, Yam22b, Yam22c}).

In summary, the convex integration method has been successfully applied in fluid dynamics to study the incompressible Euler and Navier-Stokes equations, as well as their stochastic counterparts.  However, constructing Leray solutions using convex integration remains an open problem, and a recent non-uniqueness result for Leray solution with forced Navier-Stokes system was obtained by a different method \cite{ABC22}.

\noindent{\bf Organization of the paper.}
In Section~\ref{s:not} we collect the basic notations and preliminaries
used throughout the paper. In  Section \ref{cogs}  we recall the construction of stochastic
objects and present a formal decomposition
of the system into the system for $v^1$ and $v^2$ as discussed above. Then we explain the set-up of the iterative convex integration procedure and provide proofs of our main results, namely, Theorem \ref{thm:globalu} and Corollary \ref{thm:noninlaw}. We give estimates of $v^1_q$ in Section \ref{est:vq1}. Section \ref{proofof3.1} is devoted to the core of the convex integration construction, namely, the iteration Proposition \ref{prop:3.1}. The proof of Theorem \ref{thm:unique} is shown in Section \ref{potuni}. Finally in Appendix we collect several auxiliary results.

\section{Preliminaries}
\label{s:not}

  Throughout the paper, we employ the notation $a\lesssim b$ if there exists a constant $c>0$ such that $a\leq cb$.
 \subsection{Function spaces}
  Given a Banach space $E$ with a norm $\|\cdot\|_E$ and $T>0$, we write $C_TE=C([0,T];E)$ for the space of continuous functions from $[0,T]$ to $E$, equipped with the supremum norm $\|f\|_{C_TE}=\sup_{t\in[0,T]}\|f(t)\|_{E}$.
  For $\alpha\in(0,1)$ we  define $C^\alpha_TE$ as the space of $\alpha$-H\"{o}lder continuous functions from $[0,T]$ to $E$, endowed with the norm $\|f\|_{C^\alpha_TE}=\sup_{s,t\in[0,T],s\neq t}\frac{\|f(s)-f(t)\|_E}{|t-s|^\alpha}+\sup_{t\in[0,T]}\|f(t)\|_{E}.$ Here we use $C_T^\alpha$ to denote the case when $E=\mathbb{R}$.
  For $p\in [1,\infty]$ we write $L^p_TE=L^p([0,T];E)$ for the space of $L^p$-integrable functions from $[0,T]$ to $E$, equipped with the usual $L^p$-norm.
    We use $L^p$ to denote the set of  standard $L^p$-integrable functions from $\mathbb{T}^2$ to $\mathbb{R}^2$. For $s>0$, $p>1$ we set $W^{s,p}:=\{f\in L^p; \|(I-\Delta)^{\frac{s}{2}}f\|_{L^p}<\infty\}$ with the norm  $\|f\|_{W^{s,p}}=\|(I-\Delta)^{\frac{s}{2}}f\|_{L^p}$.
For $T>0$ and a domain $D\subset\R^{+}$ we denote by  $C^{N}_{T,x}$ and $C^{N}_{D,x}$, respectively, the space of $C^{N}$-functions on $[0,T]\times\mathbb{T}^{2}$ and on $D\times\mathbb{T}^{2}$, respectively,  $N\in\N_0:=\N\cup\{0\}$. The spaces are equipped with the norms
$$
\|f\|_{C^N_{T,x}}=\sum_{\substack{0\leq n+|\alpha|\leq N\\ n\in\N_{0},\alpha\in\N^{2}_{0} }}\|\partial_t^n D^\alpha f\|_{L^\infty_T L^\infty},\qquad \|f\|_{C^N_{D,x}}=\sum_{\substack{0\leq n+|\alpha|\leq N\\ n\in\N_{0},\alpha\in\N^{2}_{0} }}\sup_{t\in D}\|\partial_t^n D^\alpha f\|_{ L^\infty}.
$$
We define the projection onto null-mean functions
is $\mathbb{P}_{\neq0} f := f -\int_{\mathbb{T}^d}\!\!\!\!\!\!\!\!\!\!\!\!\!\; {}-{}\ \ f\dif x$. For a tensor $T$ , we denote its traceless part by $\mathring{T} := T- \frac{1}{d} \tr (T )\Id$. By
$\mathcal{S}^{d\times d}$ we denote the space of symmetric matrix and by $\mathcal{S}^{d\times d}_0$ the space of symmetric trace-free matrix.

We use $(\Delta_i)_{i\geq 1}$ to denote the Littlewood-Paley blocks corresponding to a dyadic partition of unity. Besov spaces on the torus with general indices $\alpha \in \mathbb{R}, p, q \in[1, \infty]$ are defined as the completion of $C^\infty(\mathbb{T}^2)$ with respect to the norm
$$\|u\|_{B_{p,q}^\alpha}:=\Big{(}\sum_{j\geq-1}2^{j\alpha q}\|\Delta_ju\|_{L^p}^q\Big{)}^{1/q}.$$
The H\"{o}lder–Besov space $C^\alpha$ is given by $C^\alpha = B^\alpha_{\infty,\infty}$, and we also set $H^\alpha = B_{2,2}^\alpha,\ \alpha \in \mathbb{R}$.

The following embedding results will be frequently used.
\bl\label{lem:2.2}
$(1).($\cite[Lemma A.2]{GIP15}$)$  Let $1\leq p_1\leq p_2\leq \infty$ and $1\leq q_1\leq q_2\leq \infty$, and let $\alpha\in\mathbb{R}$. Then $B_{p_1,q_1}^\alpha\subset B_{p_2,q_2}^{\alpha-2(1/p_1-1/p_2)}.$\\
$(2).($\cite[Theorem 4.6.1]{Tri78}$)$Let $s\in\mathbb{R},1<p<\infty,\epsilon>0$. Then $W^{s,2}=B_{2,2}^s=H^s$, and $B_{p,1}^s\subset W^{s,p}\subset B_{p,\infty}^s\subset B_{p,1}^{s-\epsilon}.$
\el

\subsection{Paraproducts}
Paraproducts were introduced by Bony in
\cite{Bon81} and they permit to decompose a product of two distributions into three parts which behave differently in terms of regularity. More precisely, using the Littlewood-Paley blocks, the product $fg$ of two Schwartz distributions $f,g \in \mathcal{S}' (\mathbb{T}^2)$ can be formally decomposed as
$$fg=f\prec g+f\succ g+f\circ g,$$
with
$$f\prec g=g\succ f=\sum_{j\geq-1}\sum_{i<j-1}\Delta_if\Delta_jg,\ \ f\circ g=\sum_{|i-j|\leq1}\Delta_if\Delta_jg.$$
Here, the paraproducts $\prec$ and $\succ$ are always well-defined and critical is the resonant product denoted by $\circ$. In general, it is only well-defined provided the sum of the regularities of $f$ and $g$ in terms of Besov spaces is strictly positive. Moreover, we have the following paraproduct estimates.
\bl\label{lem:2.3}$($\cite[Lemma 2.1]{GIP15},\cite[Proposition A.7]{MW17a}$)$
Let $\beta\in\mathbb{R},p,p_1,p_2,q\in[1,\infty]$ such that $\frac{1}{p}=\frac{1}{p_1}+\frac{1}{p_2}$. Then it holds
$$\|f\prec g\|_{B_{p,q}^\beta}\lesssim\|f\|_{L^{p_1}}\|g\|_{B_{p_2,q}^\beta},$$
and if $\alpha<0$ then
$$\|f\prec g\|_{B_{p,q}^{\alpha+\beta}}\lesssim\|f\|_{B_{p_1,q}^\alpha}\|g\|_{B_{p_2,q}^\beta}.$$
If $\alpha+\beta>0$ then it holds
$$\|f\circ g\|_{B_{p,q}^{\alpha+\beta}}\lesssim\|f\|_{B_{p_1,q}^\alpha}\|g\|_{B_{p_2,q}^\beta}.$$
\el
We denote $\succeq=\circ+\succ,\ \preceq=\circ+\prec$.

Analogously to the the real-valued case, we may define paraproducts for vector-valued distributions. More precisely, for two vector-valued distributions $f, g \in \mathcal{S}'(\mathbb{T}^d; \mathbb{R}^m)$, we use the following
tensor paraproduct notation
\begin{align*}
f \otimes g = (f_ig_j)^m_{i,j=1} = f \prec\!\!\!\!\!\!\!\bigcirc g + f \circledcirc\ g + f \succ\!\!\!\!\!\!\!\bigcirc \ g
= (f_i \prec g_j)_{i,j=1}^m + (f_i \circ g_j )^m_{i,j=1} + (f_i \succ g_j )^m_{i,j=1},
\end{align*}
and note that Lemma \ref{lem:2.3} carries over mutatis mutandis. We also denote $\succcurlyeq\!\!\!\!\!\!\!\bigcirc=\circledcirc+\succ\!\!\!\!\!\!\!\bigcirc,\ \preccurlyeq\!\!\!\!\!\!\!\bigcirc=\circledcirc+\prec\!\!\!\!\!\!\!\bigcirc$.

We also recall the following lemma for the Helmholtz projection $\mathbb{P}_{\rm{H}}$
\bl\label{lem:est:leray}$($\cite[Lemma 2.5]{HZZ21b}$)$
Assume that $\alpha\in\mathbb{R}$ and $p \in [1, \infty]$. Then for every $k, l = 1, 2$
\begin{align*}
    \|\mathbb{P}_{\rm{H}}^{kl}f\|_{B_{p.\infty}^\alpha}\lesssim\|f\|_{B_{p.\infty}^\alpha}.
\end{align*}
\el

Finally, we introduce localizers in terms of Littlewood-Paley expansions. Let $J \in\mathbb{N}_0$. For
$f\in \mathcal{S}'(\mathbb{T}^2)$ we define
$$\Delta_{>J}f=\sum_{j>J}\Delta_jf,\ \ \Delta_{\leq J}f=\sum_{j\leq J}\Delta_jf.$$
Then it holds in particular for $\alpha\leq\beta\leq \gamma$
\begin{align}\label{deltaf}
\|\Delta_{>J}f\|_{C^\alpha}\lesssim2^{-J(\beta-\alpha)}\|f\|_{C^\beta},\ \ \|\Delta_{\leq J}f\|_{C^\gamma}\lesssim2^{J(\gamma-\beta)}\|f\|_{C^\beta}.
\end{align}

\subsection{Anti-divergence operators}

 We recall the following anti-divergence operator
 $\mathcal{R}$ from \cite[Appendix D]{CL20}, which acts on vector fields $v$ as
$$(\mathcal{R}v)_{ij}=\mathcal{R}_{ijk}v_k$$
where
$$\mathcal{R}_{i j k}=- \Delta^{-1} \partial_{k} \delta_{i j}+\Delta^{-1} \partial_{i} \delta_{j k}+\Delta^{-1} \partial_{j} \delta_{i k}.$$
Then $\mathcal{R}v(x)$ is a symmetric trace-free matrix for each $x \in \mathbb{T}^2$, and $\mathcal{R}$ is a right
inverse of the $\div$ operator, i.e. $\div(\mathcal{R}v) = v-\int_{\mathbb{T}^2}\!\!\!\!\!\!\!\!\!\!\!\!\!\!\; {}-{}\ v\dif x$
By a direct computation we have for any divergence-free $v\in C^\infty(\mathbb{T}^2,\mathbb{R}^2)$
\begin{align}\label{rdeltav}
\mathcal{R}\Delta v=\nabla v+\nabla^{T} v.
\end{align}

Let $C_0^\infty(\mathbb{T}^2,\mathbb{R}^{2})$ be the set of smooth functions with zero mean. By \cite[Theorem D.2]{CL20} we know for any $1\leq p \leq\infty$, $\sigma\in\mathbb{N}$, $f\in C_0^\infty(\mathbb{T}^2,\mathbb{R}^{2})$
\begin{align}
\|\mathcal{R}f(\sigma\cdot)\|_{ L^p}\lesssim\sigma^{-1} \| f\|_{ L^p}.\label{bb1}
\end{align}

Let $C_0^\infty(\mathbb{T}^2,\mathbb{R}^{2\times 2})$ be the set of smooth matrix valued functions with zero mean. We also introduce the bilinear version $\mathcal{B}: C^\infty(\mathbb{T}^2, \mathbb{R}^2) \times C^\infty_0(\mathbb{T}^2, \mathbb{R}^{2\times 2})\to C^\infty(\mathbb{T}^2, \mathcal{S}_0^{2\times 2})$ by $$(\mathcal{B}(v,A))_{ij}=v_l\mathcal{R}_{ijk}A_{lk}-\mathcal{R}(\partial_iv_l\mathcal{R}_{ijk}A_{lk}).$$

Then by \cite[Theorem D.3]{CL20} we have $\div(\mathcal{B}(v,A))=vA-{\int_{\mathbb{T}^2}\!\!\!\!\!\!\!\!\!\!\!\!\; {}-{}\ vA\dif x}$ and for any $1 \leq p \leq\infty$
\begin{align}
\|\mathcal{B}(v,A)\|_{L^p}\lesssim\|v\|_{C^1}\|\mathcal{R}A\|_{L^p}.    \label{bb}
\end{align}

\subsection{Probabilistic elements}
Regarding the driving noise, we assume that $B$ is a  vector-valued $L^2$-cylindrical Wiener process on some stochastic basis $(\Omega, \mathcal{F},\mathbf{P})$. The time derivative of $B$ is the space-time white noise.

 \section{Convex integration set-up and results}\label{cogs}
This section is devoted to the proof of Theorem \ref{thm:globalu} and Corollary \ref{thm:noninlaw}. The goal is to establish existence of non-unique global-in-time probabilistically strong solutions to
(\ref{3.2}) for every given divergence free initial condition in $L^p\cup C^{-1+\delta}$ for $p\in(1,2),\delta>0$.

\subsection{Stochastic objects}\label{so}
Let us recall that due to \cite{DPD02}, the equation (\ref{3.2}) is locally well-posed for initial conditions in $C^{-1+\delta}$. The solution $u$ belongs
to $C([0, t_0); C^{-1+\delta}
)$ where $t_0$ is a strictly positive stopping time so that
$\|u_\epsilon-u\|_{C_{t_0}C^{-1+\delta}}\to0$
in probability. Here, $u_{\epsilon}$ denotes the solution to the regularized Navier–Stokes system
\begin{align*}
\partial u_{\epsilon}+\div(u_{\epsilon}\otimes u_{\epsilon})+\nabla p_{\epsilon}=\Delta u_{\epsilon}+\zeta_{\epsilon},\ \
\div u_{\epsilon}=0,\notag
\end{align*}
where $\zeta_{\epsilon}$ is a mollification of the space-time white noise $\zeta=\frac{\dif B}{\dif t}$.
To summarize the main ideas, let $z_{\epsilon}$ be the solution to
$$\partial_tz_{\epsilon} + \nabla p^{z_{\epsilon}} = \Delta z_{\epsilon} + \zeta_{\epsilon},\ \ \div z_{\epsilon} = 0,\ \ z_\epsilon(0)=0.$$

Then $z_\epsilon\to z$ in $L^p(\Omega,C_TC^{-\kappa})$ for every $p \in [1, \infty),\kappa\in(0,\infty)$. Moreover, the renormalized product $z^{:2:}$ can be defined as a Wick product in the sense that for any $t\in[0,T]$ there exist diverging constants $C_{\epsilon}(t) \in \mathbb{R}^{2\times2}, C_{\epsilon}^{ij}(t)\to\infty$, so that
$$z^{:2:}_\epsilon(t)= z_{\epsilon} \otimes z_{\epsilon}(t)-C_{\epsilon}(t)$$
has a well-defined limit in $L^p(\Omega,C_TC^{-\kappa})$ for $p \in [1, \infty)$. In fact, $C_{\epsilon}(t) =\mathbf{E}[z_{\epsilon} \otimes z_{\epsilon}(t)]$, which implies $\div(z^{:2:}_\epsilon(t))= \div(z_{\epsilon} \otimes z_{\epsilon}(t))$. We recall the following result from \cite[Lemma 3.2, Theorem 5.1]{DPD02}
\bp\label{prop:Z}
For every $\kappa>0,0<\delta<\frac1{4}$, there exist random distribution $z,z^{:2:}$ such that
$z_\epsilon\to z$ in $C_TC^{-\kappa} \cap C_T^{\delta/2}C^{-\kappa-\delta}$ and $z^{:2:}_\epsilon\to z^{:2:}$ in $C_TC^{-2\kappa} \cap C_T^{\delta/2}C^{-2\kappa-\delta}\ \ \mathbf{P}$-a.s. as $\epsilon\to0$.
\ep

\br
In fact   for equation (\ref{3.2}) if regularity of $z$ in $C_TC^{-\kappa}$ and $z^{:2:}$ in $C_TC^{-2\kappa}$, $\kappa\in(0,\frac{1}{3})$ our argument also works, which can be seen in (\ref{3.7}) below. Here $\kappa<\frac{1}{3}$ is necessary otherwise the term $v\otimes z+z\otimes v$ in (\ref{3.2}) is not well-defined and we need to perform a further decomposition.
\er

\subsection{Formal decomposition}
Following the way in \cite{HZZ21b,HZZ} we begin by constructing local probabilistically strong solutions with Cauchy problem for the initial condition in $L^p$ for $p\in(1,2)$. Then
we use the final value at the stopping time of any
such convex integration solution as a new initial condition for the convex integration procedure. This way we are able to extend the convex integration solutions as probabilistically strong solutions defined on the whole time interval $[0, \infty)$.

To this end, we intend to prescribe an arbitrary random initial condition $u_0\in L^p$ independent of the given Wiener process $B$. Let $(\mathcal{F}_t)_{t\geq 0}$ be the augmented joint canonical filtration on $(\Omega, \mathcal{F})$ generated by $B$ and $u_0$. Then $B$ is a $(\mathcal{F}_t)_{t\geq 0}$-Wiener process and $u_0$ is $\mathcal{F}_0$-measurable.

To deal with (\ref{3.2}) we split $v=z^{in}+v^1+v^2$ where $z^{in}=e^{t\Delta}u_0$ and $e^{t\Delta}$ is the heat semigroup. The equation for $v^1$ contains all the irregular terms of the product $(v + z) \otimes (v +z)$, whereas the regular ones are put in $v^2$. Additionally, the equation for $v^1$ shall be linear so that it
can be solved by a fixed point argument. The decomposition can be done as follows. The product $z\otimes z$ needs to
be constructed by renormalization as a Wick product denoted by $z^{:2:}$ and it is of spatial regularity $C^{-2\kappa}$.
So we have the first irregular term $z^{:2:}$. Then we write with the help of paraproducts and
Littlewood-Paley projectors
\begin{align*}
 (v^1+v^2+z^{in})\otimes z&=(v^1+v^2+z^{in})\prec\!\!\!\!\!\!\!\bigcirc\ \Delta_{>R}z \\
 &+(v^1+v^2+z^{in})\succcurlyeq\!\!\!\!\!\!\!\bigcirc\ \Delta_{>R}z+(v^1+v^2+z^{in})\otimes\Delta_{\leq R}z,
\end{align*}
where the first line on the right hand side is irregular and the second and the third are regular. We treat the symmetric term $z\otimes(v^1+v^2+z^{in})$ the same way. Here we included a
suitable cut-off $R$ to be chosen appropriately in Section \ref{est:vq1}. This eventually simplifies the fixed point argument
used to establish, for a given convex integration iteration $v_q^2$, the existence and uniqueness of $v_q^1$.
Finally, we have the regular term $(v^1 + v^2+z^{in})\otimes(v^1 + v^2+z^{in}
)$. This leads to
\begin{align}\label{3.3}
(\partial_t-\Delta)v^1+\nabla p^1+\div(z^{:2:}+(v^1+v^2+z^{in})\prec\!\!\!\!\!\!\!\bigcirc\Delta_{>R}z+\Delta_{>R}z\succ\!\!\!\!\!\!\!\bigcirc(v^1+v^2&+z^{in}))=0,\notag\\
\div v^1&=0,\\
v^1(0)&=0.\notag
\end{align}
\begin{align}
(\partial_t-\Delta)v^2+\nabla p^2&+\div((v^1+v^2+z^{in})\otimes\Delta_{\leq R}z+\Delta_{\leq R}z\otimes(v^1+v^2+z^{in}))\notag\\&+\div((v^1+v^2+z^{in})\otimes(v^1+v^2+z^{in}))\notag\\
&+\div((v^1+v^2+z^{in})\succcurlyeq\!\!\!\!\!\!\!\bigcirc\Delta_{>R}z+\Delta_{>R}z\preccurlyeq\!\!\!\!\!\!\!\bigcirc\ (v^1+v^2+z^{in}))=0,\notag\\
&\ \ \ \ \ \ \ \ \ \ \ \ \ \ \ \ \ \ \ \ \ \ \ \ \ \ \ \ \ \ \ \ \ \ \ \ \ \ \ \ \ \ \ \ \ \ \ \ \ \ \ \ \ \ \ \ \ \ \ \ \div v^2=0,\label{3.4}\\
&\ \ \ \ \ \ \ \ \ \ \ \ \ \ \ \ \ \ \ \ \ \ \ \ \ \ \ \ \ \ \ \ \ \ \ \ \ \ \ \ \ \ \ \ \ \ \ \ \ \ \ \ \ \ \ \ \ \ \ \
v^2(0)=0.\notag
\end{align}
It is easy to see that $v=z^{in}+v^1+v^2$ solves (\ref{3.2}). As in \cite{HZZ21b}, these equations need to be considered together with in the convex integration scheme and we put forward a joint iterative procedure.

\br
Compared to \cite{HZZ21b} we do not include the terms $\Delta_{>R}z\preccurlyeq\! \!\!\!\!\!\!\bigcirc\ (v^1+z^{in})$ and $(v^1+z^{in})$\\
$\succcurlyeq\!\!\!\!\!\!\!\bigcirc\ \Delta_{>R}z$ in (\ref{3.3}). These terms require regularity of $v^1+z^{in}$ and lead to a negative power of $t$ in the estimate of $\mathring{R}_q$. In \cite{HZZ21b}, the authors bound the $C_tL^1$ norm of $\mathring{R}_q$, whereas here we bound it in a weaker $L_t^1L_x^1$-norm.
\er

\subsection{Convex integration set-up}
The convex integration iteration is indexed by a parameter $q\in\mathbb{N}_0$. It will be seen that the
Reynolds stress $\mathring{R}_q$ is only required for the approximations $v_q^2$ of $v^2$, whereas the approximations $v_q^1$
of $v^1$ are obtained by a fixed point argument.

We consider an increasing sequence $\{\lambda_q\}_{q\in\mathbb{N}_0}\subset\mathbb{N}$ which diverges to $\infty$, and a sequence $\{\delta_q\}_{q\in\mathbb{N}_0\cup\{-1\}}$ $\subset (0,1)$
 which decreases to 0. We choose $a\in\mathbb{N},b\in\mathbb{N},\beta\in(0,1)$ and let
$$\lambda_q=a^{(b^q)},\ q\in \mathbb{N}_0,\ \ \delta_q=\frac{1}{2}\lambda_1^{2\beta}\lambda_q^{-2\beta},\ q\in\mathbb{N},\ \ \delta_0=\delta_{-1}=1,$$
where $\beta$ will be chosen sufficiently small and $a$ as well as $b$ will be chosen sufficiently large. In addition, we used
\begin{align}
   2\delta_{q+1}\leq\delta_q \label{para1}
\end{align}
for $q\in\mathbb{N}_0$  and
\begin{align}
  \sum_{q\geq1} \delta_q^{1/2}\leq  \frac{1}{\sqrt2}\sum_{q\geq1}a^{b\beta-qb\beta}\leq\frac{1}{\sqrt2} \frac{1}{1-a^{-b\beta}}\leq1,  \label{para2}
\end{align}
which boils down to
\begin{align}\label{ab2}
a^{2\beta(b-1)}\geq2,\ \ a^{b\beta}\geq4>2+\sqrt2,
\end{align}
which we assumed from now on. Then we choose $\alpha\in(0,1)$ small enough and define $f(q)$ satisfying $2^{f(q)}=\lambda_q^{\alpha/8}$. We require that $f(q)\in\mathbb{N}$, which can be satisfied by choosing an appropriate value of $a$. More details on the choice of these four parameters $a,b,\alpha,\beta$ will be given in Section \ref{cop} below.

At each step $q$, a triple $(v_q^1, v_q^2,\mathring{R}_q)$ is constructed solving the following system:
\begin{align}
(\partial_t-\Delta)v_q^1+\nabla p_q^1+&\div(z^{:2:}+(v_q^1+v_q^2+z^{in})\prec\!\!\!\!\!\!\!\bigcirc\Delta_{\leq f(q)}\Delta_{> R}z)\notag\\
+&\div(\Delta_{\leq f(q)}\Delta_{>R}z\succ\!\!\!\!\!\!\!\bigcirc(v_q^1+v_q^2+z^{in}))=0,\notag\\
&\ \ \ \ \ \ \ \ \ \ \ \ \ \ \ \ \ \ \ \ \ \ \ \ \ \ \ \ \ \ \ \ \ \ \ \ \ \div v_q^1=0,\notag\\
&\ \ \ \ \ \ \ \ \ \ \ \ \ \ \ \ \ \ \ \ \ \ \ \ \ \ \ \ \ \ \ \ \ \ \ \ \ v^1_q(0)=0.\label{3.5}
\end{align}
\begin{align}
(\partial_t-\Delta)v_q^2+\nabla p_q^2&+\div((v_q^1+v_q^2+z^{in})\otimes(v_q^1+v_q^2+z^{in}))\notag\\
&+\div((v_q^1+v_q^2+z^{in})\otimes\Delta_{\leq R}z+\Delta_{\leq R}z\otimes(v_q^1+v_q^2+z^{in}))\notag\\
&+\div((v_q^1+v_q^2+z^{in})\succcurlyeq\!\!\!\!\!\!\!\bigcirc
\Delta_{>R}z+\Delta_{>R}z\preccurlyeq\!\!\!\!\!\!\!\bigcirc(v_q^1+v_q^2+z^{in}))=\div \mathring{R}_q,\notag\\
&\ \ \ \ \ \ \ \ \ \ \ \ \ \ \ \ \ \ \ \ \ \ \ \ \ \ \ \ \ \ \ \ \ \ \ \ \ \ \ \ \ \ \ \ \ \ \ \ \ \ \ \ \ \ \ \ \ \ \div v_q^2=0,\notag\\
&\ \ \ \ \ \ \ \ \ \ \ \ \ \ \ \ \ \ \ \ \ \ \ \ \ \ \ \ \ \ \ \ \ \ \ \ \ \ \ \ \ \ \ \ \ \ \ \ \ \ \ \ \ \ \ \ \ \ v^2_q(0)=0.\label{3.6}
\end{align}
Here we add the localizers $\Delta_{\leq f(q)}$ in the equation of $v_{q+1}^1$, which are
used to control the blow up of a certain norm of $v_{q}^1$ as $q\to\infty$. Note that the Reynolds stress $\mathring{R}_q$ is only included in the equation for $v_q^2$. Indeed, $v_q^2$ contains the
quadratic nonlinearity which is used in the convex integration to reduce the stress.

At each iteration step $q+1$, we first use the previous stress $\mathring{R}_q$  to define the principle part
of the corrector $w_{q+1}^{(p)}$, the incompressibility corrector $w_{q+1}^{(c)}$ and the time corrector $w_{q+1}^{(t)}$ in terms
of the accelerating jet flows introduced in Appendix \ref{s:appB}. This gives the next iteration $v_{q+1}^2$ and consequently we
obtain $v_{q+1}^1$ by a fixed point argument.

As the next step, we define a stopping time which controls suitable norms of all the required
stochastic objects. Namely, by Proposition \ref{prop:Z},
 for $L \geq 2$ we define $\kappa_0:=(1-\frac{1}{p})\wedge(\frac2p-1)\leq\frac13$ and let for some $0<\kappa<\kappa_0$
\begin{align}
T_L:&=T_L^1\wedge T_L^2\wedge T_L^3,\label{3.7}\\
T_L^1:&=\inf\{t\geq0,\|z(t)\|_{C^{-\kappa}}\geq L^{1/4}\}\wedge L,\notag\\
T_L^2:&=\inf\{t\geq0,\|z\|_{C_t^{\kappa_0/2}C^{-\kappa-\kappa_0}}\geq L^{1/4}\}\wedge\inf\{t\geq0,\|z\|_{C_t^{\frac{1}{4}(1-2\kappa-\kappa_0)} C^{-\frac12+\frac{\kappa_0}2}}\geq L^{1/4}\}\wedge L,\notag\\
T_L^3:&=\inf\{t\geq0,\|z^{:2:}(t)\|_{C^{-2\kappa}}\geq L^{1/2}\}\wedge L.\notag
\end{align}
Then $T_L$ is $\mathbf{P}$-a.s. strictly positive stopping time and it holds that $T_L\to \infty$ as $L\to\infty$ $\mathbf{P}$-a.s.

We intend to solve (\ref{3.2}) for any given divergence free initial condition $u_0 \in L^p$ measurable
with respect to $\mathcal{F}_0$. However, in the first step, we take the following additional assumption: Let $N\geq2$ be given and assume in addition that $\mathbf{P}$-a.s.
\begin{align}
\|u_0\|_{L^p}\leq N.\label{omegan}
\end{align}
We keep this additional assumption on the initial condition throughout the convex integration
step in Proposition \ref{prop:3.1}. In Theorem \ref{thm:globalu*} it is relaxed to $u_0 \in L^p
\ \mathbf{P}$-a.s and finally, Corollary \ref{coro:globalu}
proves the result if $u_0 \in L^p\cup C^{-1+\delta},\ \delta>0 \ \mathbf{P}$-a.s. We suppose that there is a deterministic constant
$M_L(N)\geq\max\{L^{9},N^2L^2\},$
In the following we write $M_L$ instead of $M_L(N)$ for simplicity. We denote $A=([\frac{3p}{p-1}]+1)M_L,\sigma_q=\delta_{q},q\in\mathbb{N}_0\cup\{-1\},\gamma_q=\delta_q
,q\in\mathbb{N}_0\backslash\{3\},\gamma_3=K$. $K>1$ is
a large constant which is used to distinguish different solutions.

\br
Comparing to \cite{HZZ21b}, here we choose $\sigma_q=\delta_q$ for $q\in\mathbb{N}_0\cup\{-1\}$ and $\gamma_q=\delta_q$ for $q\in\mathbb{N}_0\backslash\{3\}$ to deduce $v\in L^2([0,T_L];H^\zeta)$ for some $\zeta>0$ (see (\ref{dvq2hzeta}) below for details).
\er

To handle the mollification around $t = 0$ needed in the convex integration below, we require that
(\ref{3.6}) is satisfied also for some negative times, namely, it holds on an interval
$[t_q, T_L]$ for certain $t_q < 0$. More precisely, we let $t_q := -1 + \sum_{ 1\leq r\leq q} \delta_r^{1/2}$ and by (\ref{para2}) we obtain $-1\leq t_q\leq0$. Here we defined $\sum_{1\leq r\leq 0}:=0.$ For $q\in\mathbb{N}_0,t\in[t_q,0)$ we assume $z^{in}(t)=e^{|t|\Delta}u_0, z(t)=v_q^1(t)=v_q^2(t)=0,\mathring{R}_q(t)=z^{in}(t)\mathring{\otimes} z^{in}(t)$.
As $v_q^2$ equals zero near zero, $\partial_t v_q^2(0)=0$, which implies by our extension that the equation (\ref{3.6}) holds
also for $t\in[t_q,0)$.
\br\label{BMSno}
Comparing to \cite{HZZ21b}, we do not extend $z^{in}$ and $\mathring{R}_q$ by their values at time zero as $u_0$ is only in $L^p$ for $p\in(1,2)$ and we need to ensure  $z^{in}{\otimes} z^{in}$ in $L^1_tL_x^1$ within the estimate for $\mathring{R}_q$.
\er
The key result is the following iterative proposition, which we prove below in Section \ref{proofof3.1} below.

\bp\label{prop:3.1}
Let $p\in(1,2)$,  $L,N\geq 2$ and assume (\ref{omegan}). There exists a choice of parameters $a,b,\alpha,\beta$ such
that the following holds true: Let $(v_q^1, v_q^2,\mathring{R}_q)$ for some $q\in\mathbb{N}_0$ be an $(\mathcal{F}_t)_{t\geq0}$-adapted solution to (\ref{3.5}) and (\ref{3.6}) satisfying
\begin{align}\label{3.8}
\left\{
\begin{aligned}
\|v_q^2\|_{L^2_{[0,T_L]}L^2}&\leq M_0(M_L^{3/4}\sum_{m=1}^{q}\delta_{m}^{1/2}+\sqrt{2}M_L^{1/4}\sum_{m=1}^{q}\gamma_{m}^{1/2})+\sqrt{2}M_0(M_L+A)^{1/2}\sum_{m=1}^{q-1}(m\sigma_{m-1})^{1/2}\\
&\leq M_0M_L^{3/4}+\sqrt{2} M_0M_L^{1/4}(K^{1/2}+1)+17M_0(M_L+A)^{1/2},\\
v_q^2(t)&=0,\ \ t\in[t_q,\sigma_q\wedge T_L],
\end{aligned}
\right.
\end{align}
for a universal constant $M_0$, where we define $\sum_{1\leq r\leq -1}:=0,\sum_{1\leq r\leq 0}:=0$, and
\begin{align}
\|v_q^2\|_{C_{t,x}^1}&\leq\lambda_q^{4}M_L^{1/2},\ \ t\in[0,T_L],\label{3.9}\\
\|v_q^2\|_{C_{T_L}L^p}&\leq M_L^{1/2}\sum_{m=1}^q\delta_m^{1/2}\leq M_L^{1/2},\label{3.10}
\end{align}
\begin{align}
\left\{
\begin{aligned}
&\|\mathring{R}_{q}\|_{L^1_{(\sigma_{q-1}\wedge T_L,T_L]}L^1}\leq M_L\delta_{q+1},\\
&\|\mathring{R}_{q}\|_{L^1_{[0, T_L]}L^1}\leq M_L\delta_{q+1}+ 2(q+1)A\sigma_{q}^{2-\frac{2}{p}},\\
&  \sup_{t_{q}\leq a\leq ({\sigma_{q}}\wedge T_L)-h} \|\mathring{R}_{q}\|_{L^1_{[a,a+h]}L^1}\leq 2(q+1)A(\frac{h}{2})^{2-\frac{2}{p}}\ for\ any\ h\in(0,(\sigma_{q}\wedge T_L)-t_q].
\end{aligned}
\right.\label{3.11}
\end{align}
Then there exists an $(\mathcal{F}_t)_{t\geq0}$-adapted process $(v_{q+1}^1,v_{q+1}^2, \mathring{R}_{q+1})$ which solves  (\ref{3.5}) and (\ref{3.6}) {on the level $q+1$} and satisfies
\begin{align}
\|v_{q+1}^2-v_q^2\|_{C_{T_L}L^{p}}&\leq M_L^{1/2}\lambda_{q+1}^{-\alpha}\leq M_L^{1/2}\delta_{q+1}^{1/2},\label{3.14}\\
\|v_{q+1}^2-v_q^2\|_{C_{T_L}W^{1/2,6/5}}&\leq M_L^{1/2}\lambda_{q+1}^{-\alpha}\leq M_L^{1/2}\delta_{q+1}^{1/2},\label{3.16}
\end{align}

\begin{align}\label{3.15}
\left\{
\begin{aligned}
\|v_{q+1}^2-v_q^2\|_{L^2_{((2\sigma_{q-1})\wedge T_L,T_L]}L^2}\leq&M_0(M_L^{1/2}\delta_{q+1}^{1/2}+\gamma_{q+1}^{1/2})(M_L^{1/2}-2\sigma_{q-1})^{1/2},\\
\|v_{q+1}^2-v_q^2\|_{L^2_{(\sigma_{q+1}\wedge T_L,(2\sigma_{q-1})\wedge T_L]}L^2}\leq&M_0((M_L+qA)^{1/2}+\gamma_{q+1}^{1/2})(2\sigma_{q-1})^{1/2},\\
v_{q+1}^2(t)=&0,\ for\ t\in[t_{q+1},\sigma_{q+1}\wedge T_L].
\end{aligned}
\right.
\end{align}

Consequently, $(v_{q+1}^1,v_{q+1}^2, \mathring{R}_{q+1})$ obeys (\ref{3.8})-(\ref{3.11}) at the level $q+1$. Furthermore,
\begin{align}
|\|v_{q+1}^2\|_{L^2_{(2\wedge T_L,T_L]}L^2}^2-\|v_q^2\|_{L^2_{(2\wedge T_L,T_L]}L^2}^2-2\gamma_{q+1}(T_L-2\wedge T_L)|\leq5M_L\delta_{q+1}.\label{3.17}
\end{align}
\ep

\br\label{rem:3.7}
As we used accelerating jets introduced in \cite{CL20}, the bound for $v^2_q$ is $L^2$ in time and for $\mathring{R}_q$ is $L^1$ in time compared to \cite{HZZ21b}. Moreover, we have extra $(q+1)A$ term in  (\ref{3.11}), which   comes from the noise part and propagates to all the other estimates. To still obtain $\mathring{R}_q\to0$ in $L^1([0,T_L];L^1)$, we keep track of the length of interval in the estimate of $\mathring{R}_q$, which present in the last line of  (\ref{3.11}).
\er

Note that Proposition \ref{prop:3.1} does not include any bounds on $v_q^1,v_{q+1}^1$. Indeed,
the definition of the new velocity $v_{q+1}^2$ does not require $v_{q+1}^1$. Once $v_{q+1}^2$ is obtained, all the necessary bounds for $v_{q+1}^1$ can be deduced from Section \ref{est:vq1} below. Specifically, in Section \ref{est:vq1}, we demonstrate the following.
\bp\label{prop:3.2}
Under the assumptions of Proposition \ref{prop:3.1}, it holds for $0<\kappa<\kappa_0$
\begin{align}
    \|v_{q+1}^1-v_q^1\|_{C_{T_L}L^2}\lesssim  \|v_{q+1}^1-v_q^1\|_{C_{T_L}B_{p,\infty}^{1-\kappa-\kappa_0}}&\lesssim M_L^{1/2}(\lambda_{q+1}^{-\alpha}+L^{5/4}\lambda_{q}^{-\frac{\alpha}{8}\kappa_0})\leq M_L^{1/2}\delta_{q+1}^{1/2},\label{3.12}\\
\|v_{q}^1\|_{C_{T_L}L^2}\lesssim\|v_q^1\|_{C_{T_L}B^{1-\kappa-\kappa_0}_{p,\infty}}&\lesssim L^{3/2}+L^{-3/4}(M_L^{1/2}+N)\leq M_L^{1/2},\label{3.13}\\
\|v_q^1\|_{C_{T_L}C^{\frac12}}&\lesssim
 L^2M_L^{1/2}\lambda_q^{\alpha}.\label{4.4}
\end{align}
\ep

We intend to start the iteration from $v_0^2 = 0$ on $[0, T_L]$. Then (\ref{3.8})-(\ref{3.10}) hold. In that case,  we obtain $v_0^1$ by (\ref{3.5}) and $\mathring{R}_0$ is the trace-free part of the matrix
\begin{align}
(v_0^1+z^{in})&{\otimes}\Delta_{\leq R}z+\Delta_{\leq R}z{\otimes}(v_0^1+z^{in})+(v_0^1+z^{in}){\otimes}(v_0^1+z^{in})\notag\\
&+(v_0^1+z^{in})\ {\succcurlyeq\!\!\!\!\!\!\!\bigcirc}\Delta_{>R}z+\Delta_{>R}z\ {\preccurlyeq\!\!\!\!\!\!\!\bigcirc}(v_0^1+z^{in}).\notag
\end{align}
By (\ref{deltaf}), (\ref{3.7}) and the choice of $R$ in Section \ref{est:vq1} we obtain
\begin{align}\label{est:deltaz}
\|\Delta_{\leq R}z\|_{C_tL^\infty}&+\|z\|_{C_tC^{-\kappa}}\lesssim\|\Delta_{\leq R}z\|_{C_tC^{\kappa_0/8}}+\|z\|_{C_tC^{-\kappa}}\notag\\
&\lesssim   (2^{(\frac{1}{8}\kappa_0+\kappa) R}+1)\|z\|_{C_tC^{-\kappa}}\lesssim(2^{\frac{9}{8}\kappa_0 R}+1)\|z\|_{C_tC^{-\kappa}}\lesssim L^{5/2}.
\end{align}
By \cite[Lemma 9]{DV} we obtain for $t\in[0,T_L],p\in(1,2)$
\begin{align}
  \|z^{in}(t)\|_{W^{2/p-1,p}}+\|z^{in}(t)\|_{L^2}\lesssim (1+t^{\frac12-\frac1p})\|u_0\|_{L^p}.\notag
\end{align}
With these bounds in hand, by Lemmas \ref{lem:2.2}, \ref{lem:2.3} with $2/p-1-\kappa>0,1-2\kappa-\kappa_0>0$ and Proposition \ref{prop:3.2} we have for $t\in[0,T_L]$
\begin{align}
\|\mathring{R}_{0}(t)\|_{L^1}&\leq 2(\|v_0^1\|_{C_tB^{1-\kappa-\kappa_0}_{p,\infty}}+\|z^{in}(t)\|_{W^{2/p-1,p}})(\|\Delta_{\leq R}z\|_{C_tL^\infty}+\|z\|_{C_tC^{-\kappa}})+\|v_0^1(t)+z^{in}(t)\|_{L^2}^2\notag\\
&\lesssim (\|v_0^1\|_{C_tB^{1-\kappa-\kappa_0}_{p,\infty}}+\|z^{in}(t)\|_{W^{2/p-1,p}})^2+(\|\Delta_{\leq R}z\|_{C_tL^\infty}+\|z\|_{C_tC^{-\kappa}})^2+\|v_0^1(t)+z^{in}(t)\|_{L^2}^2\notag\\
&\lesssim (L^{3/2}+L^{-3/4}(M_L^{1/2}+N))^2+(1+|t|^{1/2-1/p})^2N^2+L^{5}\notag\\
&\lesssim L^{5}+L^{-3/2}M_L+N^2(1+|t|^{1-2/p}),\label{est:r0t+}
\end{align}
which implies that
$$\|\mathring{R}_{0}\|_{L^1_{[0, T_L]}L^1}\lesssim L^{6}+L^{-1/2}M_L+N^2(L+L^{2-2/p})\leq
\frac{1}{2}M_L.$$
Here and in the following we used $L$ large enough to absorb the implicit constant.

For $0< t\leq \sigma_0\wedge T_L\leq 1$, we have  $t^{1-2/p}\geq1$ as $1-2/p<0$. Then we have $$\|\mathring{R}_0(t)\|_{L^1}\leq (2-2/p)At^{1-2/p}.$$
For $-1\leq t_0\leq t<0$ we have  $\mathring{R}_{0}(t)=z^{in}\mathring{\otimes} z^{in}$ and $|t|^{1-2/p}\geq1$ as $1-2/p<0$. Therefore we obtain
\begin{align}
\|\mathring{R}_{0}(t)\|_{L^1}
&\lesssim \|z^{in}(t)\|_{L^2}^2\lesssim N^2(1+|t|^{1-2/p})\leq M_L|t|^{1-2/p}.\label{est:r0t-}
\end{align}
Combining the two inequalities above and by straightforward calculations we obtain for any $0<h\leq \sigma_0\wedge T_L-t_0$
$$\sup_{t_{0}\leq a<a+h \leq \sigma_0\wedge T_L}\|\mathring{R}_{0}\|_{L^1_{[a,a+h]}L^1}\leq 2A(\frac{h}{2})^{2-2/p}.$$
Hence (\ref{3.11}) is satisfied on the level $q = 0$ as $\delta_1 = 1/2$.

\subsection{Proof of main result}
We first deduce the following result by using Propositions \ref{prop:3.1} and \ref{prop:3.2}.

\bt\label{thm:1}
There exists a $\mathbf{P}$-a.s. strictly positive stopping time $T_L$, arbitrarily large by choosing $L$ large, such that for any $\mathcal{F}_0$-measurable divergence-free initial condition $u_0\in L^p,p\in(1,2)\  \mathbf{P}$-a.s. the following holds true: There exists an $(\mathcal{F}_t)_{t\geq0}$-adapted process $(v^1
, v^2)$ such that
\begin{align*}
    v^1&\in C([0, T_L];H^\zeta)\cap C([0,T_L];B_{p,\infty}^{1-\kappa-\kappa_0}),\\
    v^2&\in C([0, T_L];L^{p})\cap L^2([0, T_L];H^\zeta)\cap C([0,T_L];W^{\frac12,\frac65})
\end{align*}
 $\mathbf{P}$-a.s. for $0<\kappa<\kappa_0$, $\zeta\in(0,1)
 $ independent of $\kappa$, and is an analytically weak solution to (\ref{3.3}) and (\ref{3.4}) with $z^{in}(t)=e^{t\Delta}u_0,
v^1(0) = v^2(0) = 0$. There are infinitely many such solutions $v^2$ and also infinitely many solutions $v=v^1+v^2+z^{in}$ to (\ref{3.2}) on $[0,T_L]$ satisfying
\begin{align}
   v\in C([0, T_L];L^{p})\cap L^2([0, T_L];H^\zeta) \cap L^1([0,T_L];W^{\frac13,1})\notag
\end{align} $\mathbf{P}$-a.s. for some $\zeta\in(0,1)$.
\et

\begin{proof}
First we consider $u_0\in L^p$ satisfying $\|u_0\|_{L^p}\leq N$. Letting $v_0^2\equiv0$, we repeatedly apply Proposition \ref{prop:3.1} and obtain $(\mathcal{F}_t)_{t\geq0}$-adapted processes $(v_q^1, v_q^2,\mathring{R}_q),q\in\mathbb{N}$. By (\ref{3.12}), (\ref{4.4}) and interpolation we obtain
\begin{align*}
    \sum_{q\geq0}\|v_{q+1}^1-v_q^1\|_{C_{T_L}H^\zeta}&\lesssim \sum_{q\geq0}\|v_{q+1}^1-v_q^1\|_{C_{T_L}L^2}^{1-3\zeta}\|v_{q+1}^1-v_q^1\|_{C_{T_L}H^{1/3}}^{3\zeta}\\
    &\lesssim \sum_{q\geq0}\|v_{q+1}^1-v_q^1\|_{C_{T_L}L^2}^{1-3\zeta}\|v_{q+1}^1-v_q^1\|_{C_{T_L}C^{1/2}}^{3\zeta}\\
    &\lesssim  \sum_{q\geq0}(M_L\delta_{q+1})^{\frac{1-3\zeta}2}(L^2M_L^{1/2}\lambda_{q+1})^{3\zeta}\lesssim M_L\sum_{q\geq0}\delta_{q+1}^{\frac{1-3\zeta- 3\zeta/\beta}{2}}\lesssim M_L,
\end{align*}
where we chose $0<\zeta<\frac{\beta}{3+3\beta}$ and used $\delta_q\leq 2^{-q}$ to deduce $\sum_{q\geq0}\delta_{q+1}^{\frac{1-3\zeta- 3\zeta/\beta}{2}}<\infty$. As a result a limit $v^1=\lim_{q\to\infty}v_q^1$ exists and lies in $
C([0, T_L];H^\zeta)$.  By (\ref{para2}) and (\ref{3.12}) we obtain $v_q^1\to v^1$ in $C([0,T_L];B_{p,\infty}^{1-\kappa-\kappa_0})$ as $q\to\infty$. Using (\ref{para2}), (\ref{3.14}) and (\ref{3.16})
we obtain $v^2=\lim_{q\to\infty}v_q^2$ exists and lies in $C([0, T_L]; L^{p} )\cap C([0,T_L];W^{1/2,6/5})$. By (\ref{3.15}) we obtain
\begin{align*}
\|v_{q+1}^2-v_q^2\|_{L^2_{[0,T_L]}L^2}&\leq 2M_0(M_L^{1/2}\delta_{q+1}^{1/2}+\gamma_{q+1}^{1/2})M_L^{1/4}+M_0(M_L+qA)^{1/2}(2\sigma_{q-1})^{1/2}\lesssim M_0M_L\delta_{q-1}^{1/3},
\end{align*}
where in the last inequality we used (\ref{para1}) to deduce $\sqrt{q}\delta_{q-1}^{1/6}\lesssim \sqrt{q}2^{-q/6}\lesssim1$.
And we note that the implicit constant depends on $K$ as $\gamma_3=K$.
Then together with (\ref{3.9}), interpolation  and H\"older's inequality we obtain that
\begin{align}
\int_0^{T_L} \|v_{q+1}^2-v_q^2\|_{H^\zeta}^2\dif t&\leq \int_0^{T_L}\|v_{q+1}^2-v_q^2\|_{L^2}^{2(1-\zeta)}\|v_{q+1}^2-v_q^2\|_{H^1}^{2\zeta}\dif t\notag\\
 &\leq \left(\int_0^{T_L} \|v_{q+1}^2-v_q^2\|_{L^2}^2\dif t\right)^{1-\zeta}\left( \int_0^{T_L} \|v_{q+1}^2-v_q^2\|_{H^1}^2\dif t\right)^\zeta\notag\\
&\lesssim( M_0^2M_L^2\delta_{q-1}^{2/3})^{1-\zeta}(M_L\lambda_{q+1}^8)^\zeta\lesssim M_0^2M_L^{2}\delta_{q-1}^{\frac23(1-\zeta)-\frac{4\zeta b^2}{\beta}},  \label{dvq2hzeta}
\end{align}
where we chose  $0<\zeta<\frac{\beta}{\beta+6 b^2}$ such that $\frac23(1-\zeta)-\frac{4\zeta b^2}{\beta}>0$. Thus we obtain $v_q^2\to v^2$ in $L^2([0, T_L];H^\zeta)$.

For $z^{in}$ part, by \cite[Lemma 9]{DV} we have
\begin{align*}
\|z^{in}(t)\|_{L^p}\lesssim \|u_0\|_{L^p},\ \|z^{in}(t)\|_{H^{\zeta}}\lesssim (1+t^{\frac12-\frac1p-\frac{\zeta}{2}})\|u_0\|_{L^p},\ \|z^{in}(t)\|_{W^{\frac13,1}}\lesssim (1+t^{-\frac16})\|u_0\|_{L^p}.
\end{align*}
Then we chose $0<\zeta<2-\frac2p$ and obtained
$z^{in}\in C([0, T_L];L^{p})\cap L^2([0, T_L];H^\zeta) \cap L^1([0,T_L];W^{\frac{1}{3},1})$ by calculation. Thus we obtain $v\in C([0, T_L];L^{p})\cap L^2([0, T_L];H^\zeta) \cap L^1([0,T_L];W^{\frac{1}{3},1})$\ $\mathbf{P}$-a.s.

Furthermore by (\ref{3.11}) we obtain $\mathring{R}_{q}\to0$ in $L^1([0, T_L];L^1)$ as $q\to\infty$. Thus $(v^1,v^2)$ satisfies equations (\ref{3.3}) and (\ref{3.4}) before $T_L$ in the analytic weak sense. Since $v_q^1(0)=v_q^2(0) = 0$
we deduce $v^1(0)=v^2(0) = 0$. Hence $v$ defined above solves (\ref{3.2}).

Next, we prove non-uniqueness of the constructed solutions. In view of (\ref{para2}) and (\ref{3.17}), we have
\begin{align}
|\ \|v^2&\|_{L^2_{(2\wedge T_L,T_L]}L^2}^2-2K(T_L-2\wedge T_L)|\notag\\
&\leq
\sum_{q=0}^\infty\left| \|v_{q+1}^2\|_{L^2_{(2\wedge T_L,T_L]}L^2}^2-\|v_q^2\|_{L^2_{(2\wedge T_L,T_L]}L^2}^2-2\gamma_{q+1}(T_L-2\wedge T_L)\right|+2\sum_{q\neq2}\gamma_{q+1}M_L^{1/2}\notag\\
&\leq5M_L\sum_{q=0}^\infty \delta_{q+1}+2\sum_{q\neq2}\gamma_{q+1}M_L^{1/2}:=c,\notag
\end{align}
and by (\ref{para2}), (\ref{3.12}) and (\ref{3.13}) we obtain
\begin{align}
 \|v^1\|_{L^2_{(2\wedge T_L,T_L]}L^2}&\leq  \|v_0^1\|_{L^2_{(2\wedge T_L,T_L]}L^2}+\sum_{q=0}^\infty \|v_{q+1}^1-v_q^1\|_{L^2_{(2\wedge T_L,T_L]}L^2}\notag\\
 &\leq M_L^{1/2}(1+\sum_{q=0}^\infty\delta_{q+1}^{1/2})T_L^{1/2}\leq M_L,\notag
\end{align}
which imply non-uniqueness by choosing different $K$. More precisely, for a given $L \geq  2$ sufficiently large it holds $\mathbf{P}(4< T_L) > 0$. The parameters $L, N$ determine $M_L(N)$ and consequently
by choosing different $K = K(L, N)$ and $K' = K'(L, N)$ so that $\sqrt{4K'
- c} -M_L>M_L +\sqrt{2KM_L^{1/2} + c}$,  we deduce that the
corresponding solutions $v_K = v_K^1 + z^{in} + v^2_K$ and $v_{K'} = v_{K'}^1 + z^{in} + v^2_{K'}$  have different $L^2L^2$-norms on the set $\{4
< T_L\}$. In fact, it is easy to see
$$ \|v^1_K + v^2_K\|_{L^2_{(2
,T_L]}L^2} \leq M_L +\sqrt{2KM_L^{1/2} + c}<\sqrt{4K'
- c} -M_L\leq \|v^1_{K'} + v^2_{K'}\|_{L^2_{(2
,T_L]}L^2} .$$
By the choice of  $K,K'$ we have different solutions.

The choice of  $\zeta$  depends on $p,\beta,b$, which are independent of $N, L,K$ and $\kappa$, see Section \ref{cop}. Thus for a general initial condition $u_0 \in L^p\ \mathbf{P}$-a.s., define $\Omega_N:= \{N-1 \leq \|u_0\|_{L^p} < N\} \in \mathcal{F}_0$. Then
the first part of this proof gives the existence of infinitely many adapted solutions $(v^{1,N},v^{2,N})$ on each $\Omega_N$.
Letting $v^{1}:= \sum_{N\in\mathbb{N}} v^{1,N} 1_{\Omega_N},v^{2}:= \sum_{N\in\mathbb{N}} v^{2,N} 1_{\Omega_N}$
concludes the proof.
\end{proof}

\bt\label{thm:globalu*}
Let $u_0 \in L^p, p\in(1,2)\ \mathbf{P}$-a.s. be an $\mathcal{F}_0$-measurable  divergence free initial condition.  Then there exist infinitely many probabilistically strong and analytically weak solutions $v$ to (\ref{3.2}) on $[0,\infty)$. Moreover, it holds $$v\in C([0,\infty);L^{p})\cap L^2_{\rm{loc}}([0,\infty);H^\zeta)\cap L^1_{\rm{loc}}([0,\infty);W^{\frac13,1})$$ $\mathbf{P}$-a.s. for some $\zeta\in(0,1)$.
\et

\begin{proof}
For $u_0\in L^p, p\in(1,2)$, by Theorem \ref{thm:1} we constructed a probabilistically strong solution $v$ before
the stopping time $T_L$ starting from the given initial condition $u_0\in L^p\ \mathbf{P}$-a.s. Since $T_L > 0\ \mathbf{P}$-a.s.,we know
\begin{align}
\|v^2(T_L)\|_{L^p}+
\|v^1(T_L)\|_{L^p}+\|z^{in}(T_L)\|_{L^p}<\infty.\notag
\end{align}
 Hence we can use the value $(v^1+v^2+z^{in})(T_L)$ as a new initial condition in Theorem \ref{thm:1}. More precisely, now we aim at solving the original system on the time interval $[T_L, T_L +\hat{T}_{L}]$ by applying Theorem \ref{thm:1} with some stopping time $\hat{T}_{L} > 0\ \mathbf{P}$ -a.s.

To this end, we define $\hat{z}^{in}(0) = (v^1+v^2+z^{in})(T_L)$, and we consider $\hat{z}(t)=z(t+T_L),\hat{z}^{:2:}(t)=z^{:2:} (t+T_L)$. Then we define the stopping time
\begin{align}
\hat{T}_L:&=\hat{T}_L^1\wedge \hat{T}_L^2\wedge \hat{T}_L^3,\notag\\
\hat{T}_L^1:&=\inf\{t\geq0,\|\hat{z}(t)\|_{C^{-\kappa}}\geq (L+1)^{1/4}\}\wedge (L+1),\notag\\
\hat{T}_L^2:&=\inf\{t\geq0,\|\hat{z}\|_{C_t^{\kappa_0/2} C^{-\kappa-\kappa_0}}\geq (L+1)^{1/4}\}\notag\\
&\ \ \ \ \ \  \ \ \ \ \wedge\inf\{t\geq0,\|\hat{z}\|_{C_t^{\frac{1-2\kappa-\kappa_0}{4}}C^{-\frac12+\frac{\kappa_0}2}}\geq (L+1)^{1/4}\}\wedge (L+1),\notag\\
\hat{T}_L^3:&=\inf\{t\geq0,\|\hat{z}^{:2:}(t)\|_{C^{-2\kappa}}\geq (L+1)^{1/2}\}\wedge (L+1).\notag
\end{align}
Then we find $T_{L+1}-T_L\leq\hat{T}_{L}$.
Similarly we construct a solution $(\hat{v}^1,\hat{v}^2,\hat{z}^{in})$  on $[0,\hat{T}_{L}]$, which is still adapted to $(\mathcal{F}_{t+T_L})_{t\geq0}$. We define $v_1(t)=(v^1+v^2+z^{in})(t)1_{\{t\leq T_L\}}+(\hat{v}^1+\hat{v}^2+\hat{z}^{in})(t-T_L)1_{\{t>T_L\}}$. By a similar argument as \cite[Theorem 1.1]{HZZ} we obtain $v_1$ satisfies the equation (\ref{3.2})  before $T_{L+1}$ and is adapted to the natural filtration $(\mathcal{F}_t)_{t\geq0}$.

Now, we can iterate the above steps to obtain $\bar{v} = v1_{\{t\leq T_L\}}+\sum_{k=1}^\infty v_k1_{\{T_{L+k-1}<t\leq T_{L+k}\}}$ and $\bar{v}\in  C([0,\infty);L^{p})\cap L^2_{\rm{loc}}([0,\infty);H^\zeta)\cap L^1_{\rm{loc}}([0,\infty);W^{\frac13,1})$ for some $\zeta\in(0,1)$, is a probabilistically strong solution.

Furthermore, as in the proof of Theorem \ref{thm:1}
we obtain infinitely many such solutions by choosing different $K$.
\end{proof}

\bc\label{coro:globalu}
Let $u_0 \in L^p\cup C^{-1+\delta},\ p\in(1,2),\ \delta>0\ \mathbf{P}$-a.s. be an $\mathcal{F}_0$-measurable divergence free initial condition.  Then there exist infinitely many probabilistically strong and analytically weak solutions $v$ to the equation (\ref{3.2}) on $[0,\infty)$.
\ec
\begin{proof}
If $u_0 \in C^{-1+\delta}$ with $\delta > 0$, by the same argument as \cite{DPD02} there exists a stopping time $0<\mathbf{t}\leq T_L$ and a local solution $v$ to (\ref{3.2}). Now $v(\mathbf{t}) \in C^{1-\kappa}$
and we can start from $v(\mathbf{t})$ and obtain infinitely many global solutions by using Theorem \ref{thm:globalu*}.
\end{proof}
Accordingly Theorem \ref{thm:globalu} is proved.

\begin{proof}[Proof of Corollary \ref{thm:noninlaw}]
Let $L > 1$ be such that $\mathbf{P}(4< T_L) > 1/2$. With the notation from the proof of Theorem \ref{thm:1} and particularly in view of (\ref{3.17}), we choose again $K, K'$ so that $$\sqrt{4K- c} >M_L,\ \ \sqrt{4K'- c} -M_L>M_L +\sqrt{2KM_L^{1/2} + c}.$$
The corresponding solutions then satisfy
$$\mathbf{P}(\sqrt{4K- c} -M_L \leq \|v^1_K + v^2_K\|_{L^2_{(2,T_L]}L^2} \leq M_L +\sqrt{2KM_L^{1/2} + c})>\frac{1}{2},$$
$$\mathbf{P}(\sqrt{4K'- c} -M_L \leq \|v^1_{K'} + v^2_{K'}\|_{L^2_{(2,T_L]}L^2} \leq M_L +\sqrt{2K'M_L^{1/2} + c})>\frac{1}{2}.$$
Since the two intervals $[\sqrt{4K- c} -M_L, M_L +\sqrt{2KM_L^{1/2} + c}]$ and $[\sqrt{4K'- c} -M_L ,M_L +\sqrt{2K'M_L^{1/2} + c}]$ are disjoint, the laws of $v_K$ and $v_{K'}$
are different. This carries over to the solutions $v_K$ and $v_{K'}$ as well as to the final solutions obtained at the end of the proof of Theorem \ref{thm:globalu*}.
\end{proof}

\section{Estimate of $v_{q}^1$}\label{est:vq1}

In this section, we work under the assumptions of Proposition \ref{prop:3.1}. The main aim is to prove the
bounds (\ref{3.12})-(\ref{4.4}). Moreover, we recall that the equation for $v_q^1$ is linear. Hence, for a given $v_q^2$ we obtain the existence and uniqueness of solution $v_q^1$ to (\ref{3.5})
by a fixed point argument together with the uniform estimate derived in the sequel. Also, if $v_q^2$
    is $(\mathcal{F}_t)_{t\geq0}$-adapted, so is $v_q^1$. This in particular gives the existence of $v_{q+1}^1$ in Proposition \ref{prop:3.1},
once the new velocity $v_{q+1}^2$ was constructed in Section \ref{proofof3.1}.

In the following, similar to \cite{HZZ21b} we make use of the localizers $\Delta_{>R}$ present in the equation for $v_q^1$ in (\ref{3.3}).
Namely, by an appropriate choice of $R$ we can always apply (\ref{deltaf}) to get a small constant in front
of terms which contain $v_q^1$. We are therefore able to absorb them into the left hand sides of the
estimates without a Gronwall argument. In the following, all the estimates are pathwise and valid before the stopping time $T_L$.

By Lemma \ref{C.3}, Lemma \ref{lem:2.3}, (\ref{deltaf}), (\ref{3.7}) and  interpolation we
have for $1<p<2,0<\kappa<\kappa_0$
\begin{align}
\|v_q^1\|_{C_t^{\frac{1}{4}(1-2\kappa-\kappa_0)}B^{1/2-\kappa_0/2}_{p,\infty}}&\lesssim
\|v_q^1\|_{C_tB_{p,\infty}^{1-\kappa-\kappa_0}}+\|v_q^1\|_{C_t^{\frac{1}{2}(1-\kappa-\kappa_0)}L^p}\notag\\
&\lesssim L\|z^{:2:}\|_{C_tC^{-\kappa-\kappa_0}}+L\|v_q^1+v_q^2+z^{in}\|_{C_tL^p}\|\Delta_{>R}z\|_{C_tC^{-\kappa-\kappa_0}}\notag\\
&\lesssim L\|z^{:2:}\|_{C_tC^{-2\kappa}}+\|v_q^1+v_q^2+z^{in}\|_{C_tL^p}L^{5/4}2^{-\kappa_0 R}\notag\\
&\lesssim L^{3/2}+(\|v_q^1\|_{C_tL^2}+\|v_q^2\|_{C_tL^p}+\|z^{in}\|_{C_tL^p})L^{5/4}2^{-\kappa_0 R}.\notag
\end{align}
Then we choose $R$ such that $2^{\kappa_0 R} = 4CL^2$ with $C$ being the implicit constant. Together with (\ref{omegan}) and (\ref{3.10}) we obtain
\begin{align}
\|v_q^1\|_{C_tL^2}&\lesssim\|v_q^1\|_{C_tB_{p,\infty}^{1-\kappa-\kappa_0}}\lesssim L^{3/2}+L^{-3/4}(\|v_q^2\|_{C_tL^{p}}+\|u_0\|_{C_tL^p})\notag\\
&\lesssim L^{3/2}+L^{-3/4}(M_L^{1/2}+N)\leq M_L^{1/2},\notag
\end{align}
where in the first inequality we used the embedding Lemma \ref{lem:2.2} to deduce $B_{p,\infty}^{1-\kappa-\kappa_0}\subset L^2$ with  $1-\kappa-\kappa_0>2(\frac{1}{p}-\frac{1}{2})$, i.e. $\kappa< \kappa_0\leq 1-\frac{1}{p}$. In the last inequality we used $M_L\geq L^4$ and chose $L$ large enough to absorb the universal constants. Then (\ref{3.13}) follows.
Moreover we have
\begin{align}
\|v_q^1\|_{C_t^{\frac{1}{4}(1-2\kappa-\kappa_0)}B^{1/2-\kappa_0/2}_{p,\infty}}\lesssim M_L^{1/2}.\label{4.2}
\end{align}

Similarly we have
\begin{align}
&\|v_{q+1}^1-v_q^1\|_{C_tB_{p,\infty}^{1-\kappa-\kappa_0}}\notag\\
&\lesssim L\|v_{q+1}^1+v_{q+1}^2-v_q^1-v_q^2\|_{C_tL^p}\|\Delta_{>R}z\|_{C_tC^{-\kappa-\kappa_0}}\notag\\
&\ \ \ \ \ \ \ \ \ \ +L\|v_q^1+v_q^2\|_{C_tL^p}\|(\Delta_{\leq f(q+1)}-\Delta_{\leq f(q)})\Delta_{>R}z\|_{C_tC^{-\kappa-\kappa_0}}\notag\\
&\lesssim L^{5/4}2^{-\kappa_0 R}\|v_{q+1}^1+v_{q+1}^2-v_q^1-v_q^2\|_{C_tL^p}+\|v_q^1+v_q^2\|_{C_tL^p}L^{5/4}\lambda_{q}^{-\frac{\alpha}{8}\kappa_0},\notag
\end{align}
which by (\ref{3.10}), (\ref{3.14}) and (\ref{3.13}) implies that
\begin{align}
\|v_{q+1}^1-v_q^1\|_{C_tL^{2}}&\lesssim\|v_{q+1}^1-v_q^1\|_{C_tB_{p,\infty}^{1-\kappa-\kappa_0}}\notag\\
&\lesssim \|v_{q+1}^2-v_q^2\|_{C_tL^{p}}+\|v_{q}^1+v_q^2\|_{C_tL^{p}} L^{5/4}\lambda_{q}^{-\frac{\alpha}{8}\kappa_0}\notag\\
&\lesssim M_L^{1/2}\lambda_{q+1}^{-\alpha}+M_L^{1/2}L^{5/4}\lambda_{q}^{-\frac{\alpha}{8}\kappa_0}\leq M_L^{1/2}\delta_{q+1}^{1/2} .\notag
\end{align}
where we used the conditions on the parameters in Section \ref{cop} below to deduce  $\alpha>\beta,\frac{\alpha}{8}\kappa_0>\beta b$ and $L^{5/4}\lambda_{q}^{-\frac{\alpha\kappa_0}{8}+b\beta}\ll1$ by choosing $a$ large enough. Then (\ref{3.12}) follows.

Finally we estimate the following bound used below. These norms may priori blow up during the iteration. By Lemma \ref{C.3}, (\ref{3.7}), (\ref{omegan}), (\ref{3.10}) and (\ref{3.13})
\begin{align}
\|v_q^1\|_{C_t^{\frac14}L^\infty}+\|v_q^1\|_{C_tC^{\frac12}}&\lesssim L\|z^{:2:}\|_{C_tC^{-2\kappa}}+L\|v_q^1+v_q^2+z^{in}\|_{C_tC^{-2/p}}\|\Delta_{\leq f(q)}\Delta_{>R}z\|_{C_tC^{\frac2p-\frac12}}\notag\\
&\lesssim L^{3/2}+L^{5/4}\lambda_q^{\frac{\alpha}{8}(\frac2p-\frac12+\kappa)}\|v_q^1+v_q^2+z^{in}\|_{C_tL^{p}}\notag\\
&\lesssim L^2M_L^{1/2}\lambda_q^{\frac{\alpha}{8}(\frac2p-\frac12+\kappa_0)}\lesssim L^2M_L^{1/2}\lambda_q^{\alpha},\notag
\end{align}
where we used the embedding $L^p\subset C^{-2/p}$ by Lemma \ref{lem:2.2} in the second step and used $\frac2p-\frac{1}{2}+\kappa_0<2$ in the last step.

\section{The main iteration-proof of Proposition \ref{prop:3.1}}\label{proofof3.1}
The proof proceeds in several main steps which are the same in many convex integration schemes. First of all, we start the construction by fixing the parameters in Section~\ref{cop} and proceed with a mollification step in Section~\ref{moll}. Section~\ref{covq+12} introduces the new iteration $v_{q+1}^2$. This is the main part of the construction which used the accelerating jets introduced in Appendix \ref{s:appB}.  Section~\ref{guina4} contains the inductive estimates of $v_{q+1}^2$, whereas in Section~\ref{subsub5} we show how the energy is controlled. Finally, in Section~\ref{dotrs}, we define the new stress $\mathring{R}_{q+1}$ and  establish the inductive moment estimate on $\mathring{R}_{q+1}$ in Section~\ref{votiefr}.

\subsection{Choice of parameters}\label{cop}
In the sequel, additional parameters will be indispensable and their value has to be carefully chosen in order to respect all the compatibility conditions appearing in the estimates below. First we denote $\kappa_1:=\frac{\kappa_0}{4}\wedge(\frac{1}{6}-\frac{\kappa_0}{2})$. For a sufficiently small $\alpha\in(0,1)$ to be chosen below, we define $l:=\lambda_{q+1}^{-\frac{3\alpha}{2}}\lambda_q^{-2}$ and have $l^{-1}\leq\lambda_{q+1}^{2\alpha}$ and $l\lambda_{q}^4\leq\lambda_{q+1}^{-\alpha} $  provided $\alpha b>4$.

In addition, we also introduce a small constant $\gamma\in(0,1)$ such that  $\gamma^*:=\frac{1}{\gamma}\in\mathbb{N}$ and satisfies the following bounds
$$\gamma<\frac{1}{56},\ \  (\frac{1}{2}-\frac{1}{p})(2-8\gamma)<-10\gamma.$$
In the sequel, we use the following bounds
$$ \alpha<\frac{1}{144}\gamma,\ \  \alpha>\frac{4}{3\kappa_1}\beta b,\ \ \alpha>\frac{16}{\kappa_0}\beta b^2,\ \  \alpha b>\frac{16}{3\kappa_1}.$$
In addition, we postulate
\begin{align}
M_L + K + qA\leq l^{-1},\ \
(\lambda_{q+1}^{-\alpha+2\beta b}+L^{5/4}\lambda_{q}^{-\frac{\alpha}{8}\kappa_0+2\beta b^2})(M_L+qA+K)^{\frac12}+M_L^{1/2}\lambda_{q+1}^{-\frac{3}2\kappa_1\alpha+2\beta b}\ll1, \label{a2}
\end{align}
which can be satisfied by choosing $a$ large enough.

The above can be obtained by choosing $\alpha=\frac{1}{145\gamma^*}$, and choosing $b\in \gamma^*\mathbb{N}$ large enough such that
$b>\frac{16}{3\kappa_1\alpha}$, and finally choosing $\beta$ small such that $\alpha>\frac{4}{3\kappa_1}\beta b, \alpha>\frac{16}{\kappa_0}\beta b^2$.
In the end we choose $a\in2^{8\times 145\gamma^*\mathbb{N}}$ large enough such that (\ref{ab2}) and (\ref{a2}) hold. The choice of $a$ above ensures $f(q)\in\mathbb{N}$. In the sequel, we increase $a$ to absorb various implicit and universal constants. From here we could choose the value of $b$ and $\beta$ depend only on $p$, and are independent of $L$, $N,K$ and $\kappa$.

\subsection{Mollification}\label{moll}
We intend to replace $v^2_q$ by a mollified velocity field $v_l$, and we define
\begin{align}
v_l=(v_q^2*_x\phi_l)*_t\varphi_l,\ \ \mathring{R}_l=(\mathring{R}_q*_x\phi_l)*_t\varphi_l\notag,
\end{align}
where $\phi_l:=\frac{1}{l^2}\phi(\frac{\cdot}{l})$ is a family of standard mollifiers on $\mathbb{R}^2$, and $\varphi_l:=\frac{1}{l}\varphi(\frac{\cdot}{l})$ is a family of standard mollifiers with support on $(0,1)$. Here the one side mollifier is used to preserve adaptedness. Then since $l\leq \delta_{q+1}^{1/2}=t_{q+1}-t_q$, by calculating and (\ref{3.6}) we obtain on $[t_{q+1},T_L]$
\begin{align}
(\partial_t-\Delta)v_l+\nabla p_l+\div N_{com}&=\div(\mathring{R}_l+{R}_{com}),\\
\div v_l&=0.\notag
\end{align}\label{4.5}
where
\begin{align*}
N_{com}:=(v_q^1&+v_q^2+z^{in})\otimes\Delta_{\leq R}z+\Delta_{\leq R}z\otimes(v_q^1+v_q^2+z^{in})+(v_q^1+v_q^2+z^{in})\otimes(v_q^1+v_q^2+z^{in})\\
&+(v_q^1+v_q^2+z^{in})\succcurlyeq\!\!\!\!\!\!\!\bigcirc\Delta_{>R}z+\Delta_{>R}z\preccurlyeq\!\!\!\!\!\!\!\bigcirc(v_q^1+v_q^2+z^{in}),
\end{align*}
and $R_{com}$ is the trace-free part of
${N}_{com}-{N}_{com}*_x\phi_l*_t\varphi_l,$
$$p_l:=(p_q^2*_x\phi_l)*_t\varphi_l-\frac{1}{2}\tr(N_{com}-N_{com}*_x\phi_l*_t\varphi_l).$$

\subsection{Construction of $v_{q+1}^2$}\label{covq+12}
Let us now proceed with the construction of the perturbation $w_{q+1}$ which then defines the next iteration by $v_{q+1}^2 := v_l +w_{q+1}$. To this end, we employ the accelerating jet flows introduced in \cite[Section 3]{CL20}, which we recall in Appendix \ref{s:appB}. In particular, the building blocks $W_{(\xi)}$ for $\xi\in\Lambda$ are defined in (\ref{Wxi}) and the set $\Lambda$ is introduced  in Lemma \ref{lem:B.1}. The necessary estimates are collected in (\ref{2.4}). For the accelerating jet flows we choose the following parameters
$$ \sigma=\lambda_{q+1}^{2\gamma},\ \ {\eta}=\lambda_{q+1}^{16\gamma},\ \ \nu=\lambda_{q+1}^{1-8\gamma},\ \ \mu=\lambda_{q+1},\ \ {\theta}=\lambda_{q+1}^{1+5\gamma}.$$
It is required that $b\in \gamma^*\mathbb{N}$ to ensure the above parameters are integers.

As the next step, we shall define certain amplitude functions used in the definition of the perturbations $w_{q+1}$ similarly as \cite[Section 5.1.2]{HZZ},
$$\rho:=2\sqrt{l^2+|\mathring{R}_l|^2}+\gamma_{q+1},$$
 it follows that $\rho$ is $(\mathcal{F}_ t)_{t\geq 0}$-adapted. Then by the estimate of $\rho$ in Appendix \ref{s:appA} we have for $N\geq1$
\begin{align}
\|\rho\|_{C_{[(2\sigma_{q-1})\wedge T_L,T_L],x}^N}\lesssim &\ l^{2-7N}\delta_{q+1}M_L+\gamma_{q+1},\notag\\
\|\rho\|_{C_{[0,T_L],x}^N}\lesssim &\ l^{2-7N}(M_L+qA)+\gamma_{q+1}.\label{4.6}
\end{align}

We define the amplitude functions
\begin{align}
a_{(\xi)}(\omega,x,t):=\rho^{1/2}(\omega,x,t)\gamma_{\xi}(\Id-\frac{\mathring{R}_l(\omega,x,t)}{\rho(\omega,x,t)}), \label{4.7}
\end{align}
where $\gamma_\xi$ is introduced in Lemma \ref{lem:B.1}. By $|\mathrm{Id}-\frac{\mathring{R}_l}{\rho}-\Id|\leq\frac{1}{2}$, it holds by Lemma \ref{lem:B.1}
\begin{align}
\rho\Id-\mathring{R}_l=\sum_{\xi\in\Lambda}\rho\gamma_{\xi}^2(\Id-\frac{\mathring{R}_l}{\rho})\xi\otimes\xi=\sum_{\xi\in\Lambda}a_{(\xi)}^2\xi\otimes\xi.\label{phoid-r}
\end{align}
Then by the estimate of $a_{(\xi)}$ in Appendix \ref{s:appA} we have for $N\geq0$

\begin{align}
\|a_{(\xi)}\|_{C_{[(2\sigma_{q-1})\wedge T_L,T_L],x}^N}&\lesssim l^{-8-7N}(\delta_{q+1}^{1/2}M_L^{1/2}+\gamma_{q+1}^{1/2}),\notag\\
\|a_{(\xi)}\|_{C_{[0,T_L],x}^N}&\lesssim l^{-8-7N}((M_L+qA)^{1/2}+\gamma_{q+1}^{1/2}).\label{4.8}
\end{align}

With these preparations in hand, similarly as \cite[Section 3.4]{CL20} we define the principal part
\begin{align}\label{4.9}
w_{q+1}^{(p)}:=\sum_{\xi\in\Lambda}a_{(\xi)}g_{(\xi)}W_{(\xi)},
\end{align}
where $W_{(\xi)},g_{(\xi)}$ are introduced in Appendix \ref{s:appB}. Then by (\ref{wpingjun}), (\ref{wbujiao}) and (\ref{phoid-r}) we obtain
\begin{align}\label{4.10}
w_{q+1}^{(p)}\otimes w_{q+1}^{(p)}+\mathring{R}_l&=\sum_{\xi\in\Lambda}a_{(\xi)}^2g_{(\xi)}^2
(W_{(\xi)}\otimes W_{(\xi)})-\sum_{\xi\in\Lambda}a_{(\xi)}^2\xi\otimes\xi+\rho \Id\notag\\
&=\sum_{\xi\in\Lambda}a_{(\xi)}^2g_{(\xi)}^2\mathbb{P}_{\neq0}(W_{(\xi)}\otimes W_{(\xi)})+\sum_{\xi\in\Lambda}a_{(\xi)}^2(g_{(\xi)}^2-1)\xi\otimes\xi+\rho \Id,
\end{align}
where we recall the notation $\mathbb{P}_{\neq0}f:=f-\int_{\mathbb{T}^2}\!\!\!\!\!\!\!\!\!\!\!\!\; {}-{}\  f\dif x$.

We define the incompressibility corrector by
\begin{align}
w_{q+1}^{(c)}:=\sum_{\xi\in\Lambda}\left(a_{(\xi)}g_{(\xi)}W_{(\xi)}^{(c)}+\sigma^{-1}\nabla^{\perp}a_{(\xi)}g_{(\xi)}\Psi_{(\xi)}\right),\label{4.11}
\end{align}
where $W_{(\xi)}^{(c)}$ and $\Psi_{(\xi)}$ are introduced in Appendix \ref{s:appB}. By (\ref{2.6}) we have
\begin{align}\label{wp+wc}
  w_{q+1}^{(p)}+w_{q+1}^{(c)}=\sigma^{-1}\sum_{\xi\in\Lambda}\nabla^{\perp}[a_{(\xi)}g_{(\xi)}\Psi_{(\xi)}],
\end{align}
therefore $\div(w_{q+1}^{(p)}+w_{q+1}^{(c)})=0$.

The temporal corrector $w_{q+1}^{(t)}$ is defined as
\begin{align}
    w_{q+1}^{(t)} = w_{q+1}^{(o)} + w_{q+1}^{(a)},\label{t=0+a}
\end{align}
where $w_{q+1}^{(o)}$ is the temporal oscillation corrector
\begin{align}\label{4.12}
w_{q+1}^{(o)}:=-\sigma^{-1}\mathbb{P}_{\rm{H}}\mathbb{P}_{\neq0}\sum_{\xi\in\Lambda}h_{(\xi)}\div(a_{(\xi)}^2\xi\otimes\xi),
\end{align}
 and $w_{q+1}^{(a)}$ is the acceleration corrector
\begin{align}\label{4.13}
w_{q+1}^{(a)}:=-\theta^{-1}\sigma\mathbb{P}_{\rm{H}}\mathbb{P}_{\neq0}\sum_{\xi\in\Lambda}a_{(\xi)}^2g_{(\xi)}|W_{(\xi)}|^2\xi.
\end{align}
Here $\mathbb{P}_{\rm{H}}$ is the Helmholtz projection, and $h_{(\xi)}, g_{(\xi)}$ are given in Appendix \ref{s:appB}.

Since $\rho$ and $\mathring{R}_l$ are $(\mathcal{F}_t)_{t\geq0}$-adapted, we know $a_{(\xi)}$ is $(\mathcal{F}_t)_{t\geq0}$-adapted. Moreover, as $W_{(\xi)},W_{(\xi)}^{(c)}$, $\Psi_{(\xi)}, g_{(\xi)}$ are deterministic, $w_{q+1}^{(p)},w_{q+1}^{(c)},w_{q+1}^{(t)}$ are also $(\mathcal{F}_t)_{t\geq0}$-adapted.

Let us introduce a smooth cut-off function
\begin{align}\label{def:chi}
\chi(t)=\left\{
\begin{aligned}
&0,\ \ \ \ \ \ \ t\leq\sigma_{q+1},\\
&\in(0,1),t\in(\sigma_{q+1},2\sigma_{q+1}),\\
&1,\ \ \ \ \ \ \ t\geq2\sigma_{q+1}.
\end{aligned}
\right.
\end{align}
Note that $\|\chi'\|_{C_t^0}\leq \sigma_{q+1}^{-1}$ which has to be taken into account in the estimates of $C_{t,x}^1$ and $\mathring{R}_{q+1}$ below.

We define the perturbations $\tilde{w}_{q+1}^{(p)},\tilde{w}_{q+1}^{(c)},\tilde{w}_{q+1}^{(t)}$ as follows:
$$\tilde{w}_{q+1}^{(p)}=w_{q+1}^{(p)}\chi,\ \ \tilde{w}_{q+1}^{(c)}=w_{q+1}^{(c)}\chi,\ \ \tilde{w}_{q+1}^{(t)}=w_{q+1}^{(t)}\chi^2.$$
Then $\tilde{w}_{q+1}^{(p)},\tilde{w}_{q+1}^{(c)},\tilde{w}_{q+1}^{(t)}$ are $(\mathcal{F}_t)_{t\geq0}$-adapted.

Finally, the total perturbation $w_{q+1}$ is defined by
$$w_{q+1}:=\tilde{w}_{q+1}^{(p)}+\tilde{w}_{q+1}^{(c)}+\tilde{w}_{q+1}^{(t)},$$
which is $(\mathcal{F}_t)_{t\geq0}$-adapted and divergence-free with mean zero.

The new velocity $v_{q+1}^2$ is defined as
\begin{align}
v_{q+1}^2=v_l+w_{q+1}.\notag
\end{align}
which is $(\mathcal{F}_t)_{t\geq0}$-adapted and divergence-free with mean zero.

\subsection{Verification of the inductive estimates for $v_{q+1}^2$}\label{guina4}
First we estimate $L^m$-norm for $m\in(1,\infty)$ and then verify the inductive estimates (\ref{3.14})-(\ref{3.15}) and (\ref{3.8})-(\ref{3.10}) at the level $q+1$.

We first consider bound $w^{(p)}_{q+1}$ in $L^2$.
We know that $W_{(\xi)}$ is $(\mathbb{T}/\sigma)^2$-periodic, so by (\ref{4.7}), (\ref{4.8}), (\ref{4.9}), (\ref{2.4}) and applying Theorem \ref{ihiot} with
$p=2$ we obtain for $t\in((2\sigma _{q-1})\wedge T_L,T_L]$
\begin{align}
\|\tilde{w}_{q+1}^{(p)}(t)\|_{L^2}&\lesssim\sum_{\xi\in\Lambda}|g_{(\xi)}(t)|(\|a_{(\xi)}\|_{L^2}\|W_{(\xi)}\|_{L^2}+\sigma^{-1/2}\|a_{(\xi)}(t)\|_{C^1}\|W_{(\xi)}\|_{L^2})\notag\\
&\lesssim\sum_{\xi\in\Lambda}|g_{(\xi)}(t)|(\|a_{(\xi)}\|_{L^2}+\lambda_{q+1}^{30\alpha-\gamma}(\delta_{q+1}^{1/2}M_L^{1/2}+\gamma_{q+1}^{1/2}))\notag\\
&\lesssim \sum_{\xi\in\Lambda}|g_{(\xi)}(t)|(\|\rho(t)\|_{L^1}^{1/2}+\delta_{q+1}^{1/2}M_L^{1/2}+\gamma_{q+1}^{1/2}),\notag
\end{align}
where we used the conditions on the parameters to deduce $30\alpha-\gamma<0$.
Then we obverse that $g_{(\xi)}$ is $\mathbb{T}/\sigma$-periodic. To again apply Lemma \ref{ihiot}, it is possible to find  $n_1,n_2\in\mathbb{N}$  satisfying $0\leq\frac{n_1-1}{\sigma}<(2\sigma_{q-1})\wedge T_L\leq \frac{n_1}{\sigma},\frac{n_2-1}{\sigma}\leq T_L< \frac{n_2}{\sigma}$, and we define $\rho(t)=\rho(T_L)$ for $t\in(T_L,\frac{n_2}{\sigma}]$. Then by Lemma \ref{ihiot} we obtain
\begin{align*}
\left\|g_{(\xi)}(t)\|\rho(t)\|_{L^1}^{1/2}\right\|&_{L^2_{((2\sigma_{q-1})\wedge T_L,T_L]}}\lesssim \|g_{(\xi)}\|_{L^2_{[\frac{n_1-1}{\sigma},\frac{n_1}{\sigma}]}}\|\rho\|_{C_{[0,T_L],x}^0}^{1/2}+\left\|g_{(\xi)}(t)\|\rho(t)\|_{L^1}^{1/2}\right\|_{L^2_{[\frac{n_1}{\sigma},\frac{n_2}{\sigma}]}}\\
&\lesssim\|g_{(\xi)}\|_{L^2_{[\frac{n_1-1}{\sigma},\frac{n_1}{\sigma}]}}\|\rho\|_{C_{[0,T_L],x}^0}^{1/2}+\|g_{(\xi)}\|_{L^2([0,1])}\|\rho\|_{L^1_{[\frac{n_1}{\sigma},\frac{n_2}{\sigma}]}L^1}^{1/2}\notag\\
&\ \ \ \ +\|g_{(\xi)}\|_{L^2([0,1])}\sigma^{-\frac12}\left\|\|\rho(\cdot)\|_{L^1}^{1/2}\right\|_{C^{0,1}{([\frac{n_1}{\sigma},\frac{n_2}{\sigma}])}}(T_L-(2\sigma_{q-1})\wedge T_L+\frac{1}{\sigma})^{\frac12}.
\end{align*}
Then by the definition of $\rho$, (\ref{3.11}) and (\ref{rho0}) we have

\begin{align*}
   \|\rho\|_{C_{[0,T_L],x}^0}^{1/2}&\lesssim l^{-2}(M_L+qA+\gamma_{q+1})^{1/2},\\
    \|\rho\|^{1/2}_{L^1_{[\frac{n_1}{\sigma},\frac{n_2}{\sigma}]}L^1}&\lesssim (l+\|\mathring{R}_{q}\|_{L^1_{(\sigma_{q-1}\wedge T_L,T_L]}L^1}+\gamma_{q+1})^{\frac12}(T_L-(2\sigma_{q-1})\wedge T_L+\frac{1}{\sigma})^{\frac12}\\&\lesssim (\delta_{q+1}^{1/2}M_L^{1/2}+\gamma_{q+1}^{1/2})(T_L-(2\sigma_{q-1})\wedge T_L+\frac{1}{\sigma})^{\frac12},
\end{align*}
and
$$\partial_t(\|\rho(t)\|_{L^1}^{1/2})\lesssim \|\rho(t)\|_{L^1}^{-1/2}\|\partial_t\rho(t)\|_{L^1}\lesssim l^{-1/2}\|\partial_t\rho(t)\|_{L^1},$$
which together with  (\ref{4.6}) imply that
$$\left\|\|\rho(\cdot)\|_{L^1}^{1/2}\right\|_{C^{0,1}{([\frac{n_1}{\sigma},\frac{n_2}{\sigma}])}}\lesssim l^{-1/2}\|\rho\|_{C_{[0,T_L],x}^1}\lesssim l^{-11/2}(M_L+qA+\gamma_{q+1}).$$
So together with (\ref{4.6}), (\ref{gwnp}) and Remark \ref{rem:gxi} we obtain
\begin{align*}
&\left\|g_{(\xi)}(t)\|\rho(t)\|_{L^1}^{1/2}\right\|_{L^2_{((2\sigma_{q-1})\wedge T_L,T_L]}}\\
&\lesssim\sigma^{-1/2}l^{-2}(M_L+qA+\gamma_{q+1})^{1/2}+(\delta_{q+1}^{1/2}M_L^{1/2}+\gamma_{q+1}^{1/2})(T_L-(2\sigma_{q-1})\wedge T_L+\frac{1}{\sigma})^{\frac12}\notag\\
&\ \ \ \ +\sigma^{-\frac12}l^{-11/2}(M_L+qA+\gamma_{q+1})(T_L-(2\sigma_{q-1})\wedge T_L+\frac{1}{\sigma})^{\frac12}\\
&\lesssim (\delta_{q+1}^{1/2}M_L^{1/2}+\gamma_{q+1}^{1/2})(T_L-(2\sigma_{q-1})\wedge T_L+\frac{1}{\sigma})^{\frac12},
\end{align*}
 where we used the conditions on the parameters to deduce $ M_L+qA+K\leq l^{-1}$ and $13\alpha-\gamma<-\beta$. With these bounds we can deduce by Remark \ref{rem:gxi}
\begin{align}
\|\tilde{w}^{(p)}_{q+1}\|_{L^2_{((2\sigma_{q-1})\wedge T_L,T_L]}L^2}&\lesssim(\delta_{q+1}^{1/2}M_L^{1/2}+\gamma_{q+1}^{1/2})(T_L-(2\sigma_{q-1})\wedge T_L+\frac{2}{\sigma})^{1/2}\notag\\
&\leq \frac{1}{2}M_0(\delta_{q+1}^{1/2}M_L^{1/2}+\gamma_{q+1}^{1/2})(M_L^{1/2}-2\sigma_{q-1})^{1/2}.\label{wp21}
\end{align}
Here $M_0$ is a universal constant. Also we used
\begin{align}\label{tl-sigma}
    T_L-(2\sigma_{q-1})\wedge T_L+\frac{2}{\sigma}\leq M_L^{1/2}-2\sigma_{q-1}.
\end{align}
In fact for $2\sigma_{q-1}\geq T_L$ we have $T_L-(2\sigma_{q-1})\wedge T_L+\frac{2}{\sigma}=\frac{2}{\sigma}\leq M_L^{1/2}-2\sigma_{q-1}$ by the choice of $M_L$, and for $2\sigma_{q-1}< T_L$ we have $T_L-(2\sigma_{q-1})\wedge T_L+\frac2{\sigma}=T_L-2\sigma_{q-1}+\frac{2}{\sigma}\leq M_L^{1/2}-2\sigma_{q-1}$.

And similarly for $t\in(\sigma_{q+1}\wedge T_L, T_L]$ we obtain
\begin{align}
\|\tilde{w}_{q+1}^{(p)}(t)\|_{L^2}\lesssim   \sum_{\xi\in\Lambda} |g_{(\xi)}(t)|(\|\rho(t)\|_{L^1}^{1/2}+(M_L+qA)^{1/2}+\gamma_{q+1}^{1/2}).\notag
\end{align}
By the similar argument above, it is possible to find $n_3,n_4\in\mathbb{N}$ satisfying $\frac{n_3-1}{\sigma}<\sigma_{q+1}\wedge T_L\leq \frac{n_3}{\sigma}$, $\frac{n_4-1}{\sigma}<(2\sigma_{q-1})\wedge T_L\leq \frac{n_4}{\sigma}$, and define $\rho(t)=\rho(T_L)$ for $t\in(T_L,\frac{n_4}{\sigma}]$ if $T_L< 2\sigma_{q-1}$. Then by Lemma \ref{ihiot}, (\ref{4.6}), (\ref{rho0}), (\ref{gwnp}), Remark \ref{rem:gxi} we obtain
\begin{align}
&\ \ \|\tilde{w}^{(p)}_{q+1}\|_{L^2_{(\sigma_{q+1}\wedge T_L,(2\sigma_{q-1})\wedge T_L]}L^2}\leq\|\tilde{w}^{(p)}_{q+1}\|_{L^2_{[\frac{n_3-1}{\sigma},\frac{n_4}{\sigma}]}L^2}\notag\\
&\lesssim((M_L+qA)^{1/2}+\gamma_{q+1}^{1/2})((2\sigma_{q-1})\wedge T_L-\sigma_{q+1}\wedge T_L+\frac{2}{\sigma})^{1/2}+\sum_{\xi\in\Lambda}\|g_{(\xi)}\|_{L^2[0,1]}\|\rho\|_{L^1_{[\frac{n_3-1}{\sigma},\frac{n_4}{\sigma}]}L^1}^{1/2}\notag\\
& +\sum_{\xi\in\Lambda}\|g_{(\xi)}\|_{L^2[0,1]}\sigma^{-\frac12}\left\|\|\rho(\cdot)\|_{L^1}^{1/2}\right\|_{C^{0,1}{[\frac{n_3-1}{\sigma},\frac{n_4}{\sigma}]}}\left((2\sigma_{q-1})\wedge T_L-\sigma_{q+1}\wedge T_L+\frac{2}{\sigma}\right)^{\frac12}\notag\\
&\lesssim ((M_L+qA)^{1/2}+\gamma_{q+1}^{1/2}+\lambda_{q+1}^{11\alpha-\gamma}(M_L+qA+\gamma_{q+1}))(2\sigma_{q-1})^{1/2}\notag\\
&\leq \frac{1}{2}M_0((M_L+qA)^{1/2}+\gamma_{q+1}^{1/2}))(2\sigma_{q-1})^{1/2},\label{wp22}
\end{align}
where we used the conditions on the parameters to deduce $ M_L+qA+K\leq l^{-1}$ and $12\alpha-\gamma<0$. Here $M_0$ is a universal constant. Also we used
\begin{align}
  (2\sigma_{q-1})\wedge T_L-\sigma_{q+1}\wedge T_L+\frac{2}{\sigma}\leq 2\sigma_{q-1}. \label{sigma-sigma}
\end{align}
In fact by the choice of parameters we have $\gamma>\beta$ which yields that $\sigma_{q+1}>\frac{2}{\sigma}$. Thus for $T_L>\sigma_{q+1}$ we have $(2\sigma_{q-1})\wedge T_L-\sigma_{q+1}\wedge T_L+\frac{2}{\sigma}\leq2\sigma_{q-1}-\sigma_{q+1}+\frac{2}{\sigma}\leq 2\sigma_{q-1}$, and for $T_L\leq \sigma_{q+1}$ we obtain $(2\sigma_{q-1})\wedge T_L-\sigma_{q+1}\wedge T_L+\frac{2}{\sigma}=\frac{2}{\sigma}< \sigma_{q+1}\leq2\sigma_{q-1}$.

Thus by (\ref{wp21}) and (\ref{wp22}) we obtain the following bound, which is used below,
\begin{align}
\|\tilde{w}_{q+1}^{(p)}\|_{L^2_{(\sigma_{q+1}\wedge T_L,T_L]}L^2}
&\lesssim ((M_L+qA)^{1/2}+\gamma_{q+1}^{1/2})M_L^{1/4}.\label{wpq+1qtl}
\end{align}

For general $L^m$-norm with $m\in(1,\infty)$, by (\ref{4.8}), (\ref{4.9}) and (\ref{2.4}), for $t\in((2\sigma_{q-1})\wedge T_L,T_L]$ we have
\begin{align}
\|{w}_{q+1}^{(p)}(t)\|_{L^m}&\lesssim \sum_{\xi\in\Lambda}\|a_{(\xi)}(t)\|_{L^\infty}|g_{(\xi)}(t)|\|W_{(\xi)}(t)\|_{L^m}\notag\\
&\lesssim\sum_{\xi\in\Lambda} (M_L^{1/2}\delta_{q+1}^{1/2}+\gamma_{q+1}^{1/2})l^{-8}|g_{(\xi)}(t)|(\nu\mu)^{1/2-1/m},\notag\\
&\lesssim(M_L^{1/2}\delta_{q+1}^{1/2}+\gamma_{q+1}^{1/2})\lambda_{q+1}^{16\alpha+(1/2-1/m)(2-8\gamma)}\sum_{\xi\in\Lambda}|g_{(\xi)}(t)|.\label{4.16}
\end{align}
By (\ref{4.8}), (\ref{4.11}) and (\ref{2.4}) we obtain
\begin{align}
\|{w}_{q+1}^{(c)}(t)\|_{L^m}&\lesssim\sum_{\xi\in\Lambda}|g_{(\xi)}(t)|\left(\|a_{(\xi)}(t)\|_{L^\infty}\|W^{(c)}_{(\xi)}(t)\|_{L^m}+\sigma^{-1}\|\nabla^{\perp} a_{(\xi)}(t)\|_{L^\infty}\|\Psi_{(\xi)}(t)\|_{L^m}\right)\notag\\
&\lesssim \sum_{\xi\in\Lambda}|g_{(\xi)}(t)|(M_L^{1/2}\delta_{q+1}^{1/2}+\gamma_{q+1}^{1/2})l^{-15}(\nu\mu^{-1}(\nu\mu)^{1/2-1/m}+\sigma^{-1}\mu^{-1}(\nu\mu)^{1/2-1/m})\notag\\
&\lesssim(M_L^{1/2}\delta_{q+1}^{1/2}+\gamma_{q+1}^{1/2})\lambda_{q+1}^{30\alpha+(1/2-1/m)(2-8\gamma)-8\gamma}\sum_{\xi\in\Lambda}|g_{(\xi)}(t)|.\label{4.17}
\end{align}
By (\ref{4.8}), (\ref{4.12}), (\ref{hlinfty}) and (\ref{2.4}) we obtain
\begin{align}
\|{w}_{q+1}^{(o)}(t)\|_{L^m}&\lesssim \sigma^{-1}\sum_{\xi\in\Lambda}\|h_{(\xi)}(t)\|_{L^\infty}\|\div(a_{(\xi)}^2(t)\xi\otimes\xi)\|_{L^m}\notag\\
&\lesssim (M_L^{1/2}\delta_{q+1}^{1/2}+\gamma_{q+1}^{1/2})^2\lambda_{q+1}^{46\alpha-2\gamma}.\label{4.18}
\end{align}
By (\ref{4.8}), (\ref{4.13}) and (\ref{2.4}) we obtain
\begin{align}
\|{w}_{q+1}^{(a)}(t)\|_{L^m}&\lesssim\sigma\theta^{-1}\sum_{\xi\in\Lambda}|g_{(\xi)}(t)|\|a_{(\xi)}(t)\|_{L^\infty}^2\|W_{(\xi)}(t)\|_{L^{2m}}^2\notag\\
&\lesssim \sum_{\xi\in\Lambda}(M_L^{1/2}\delta_{q+1}^{1/2}+\gamma_{q+1}^{1/2})^2l^{-16}\sigma\theta^{-1}|g_{(\xi)}(t)|(\nu\mu)^{1-1/m}\notag\\
&\lesssim(M_L^{1/2}\delta_{q+1}^{1/2}+\gamma_{q+1}^{1/2})^2\lambda_{q+1}^{32\alpha+(1-1/m)(2-8\gamma)-3\gamma-1}\sum_{\xi\in\Lambda}|g_{(\xi)}(t)|.\label{4.19}
\end{align}

And for $t\in(\sigma_{q+1}\wedge T_L,T_L]$  the estimates are similar with $(M_L^{1/2}
\delta_{q+1}^{1/2}+\gamma_{q+1}^{1/2})$ replaced by $((M_L+qA)^{1/2}+\gamma_{q+1}^{1/2})$. More specifically for $m\in(1,\infty)$
\begin{align}
\|{w}_{q+1}^{(p)}(t)\|_{L^m}
&\lesssim((M_L+qA)^{1/2}+\gamma_{q+1}^{1/2})\lambda_{q+1}^{16\alpha+(1/2-1/m)(2-8\gamma)}\sum_{\xi\in\Lambda}|g_{(\xi)}(t)|,\label{4.16*}\\
\|{w}_{q+1}^{(c)}(t)\|_{L^m}
&\lesssim((M_L+qA)^{1/2}+\gamma_{q+1}^{1/2})\lambda_{q+1}^{30\alpha+(1/2-1/m)(2-8\gamma)-8\gamma}\sum_{\xi\in\Lambda}|g_{(\xi)}(t)|,\label{4.17*}\\
\|{w}_{q+1}^{(o)}(t)\|_{L^m}
&\lesssim ((M_L+qA)^{1/2}+\gamma_{q+1}^{1/2})^2\lambda_{q+1}^{46\alpha-2\gamma},\label{4.18*}\\
\|{w}_{q+1}^{(a)}(t)\|_{L^m}
&\lesssim((M_L+qA)^{1/2}+\gamma_{q+1}^{1/2})^2\lambda_{q+1}^{32\alpha+(1-1/m)(2-8\gamma)-3\gamma-1}\sum_{\xi\in\Lambda}|g_{(\xi)}(t)|.\label{4.19*}
\end{align}

Combining (\ref{4.17})-(\ref{4.19}) and $\chi(t)\leq1$ we obtain for $t\in((2\sigma _{q-1})\wedge T_L,T_L]$
\begin{align*}
\|\tilde{w}^{(c)}_{q+1}(t)+\tilde{w}^{(t)}_{q+1}(t)\|_{L^2}&\lesssim (M_L^{1/2}\delta_{q+1}^{1/2}+\gamma_{q+1}^{1/2})(\lambda_{q+1}^{30\alpha-8\gamma}
\sum_{\xi\in\Lambda}|g_{(\xi)}(t)|+\lambda_{q+1}^{33\alpha-7\gamma}\sum_{\xi\in\Lambda}|g_{(\xi)}(t)|+\lambda_{q+1}^{47\alpha-2\gamma})\\
&\lesssim(M_L^{1/2}\delta_{q+1}^{1/2}+\gamma_{q+1}^{1/2})(\sum_{\xi\in\Lambda}|g_{(\xi)}(t)|+1),
\end{align*}
where we used the conditions on the parameters to deduce  $M_L+K\leq l^{-1}$ and $47\alpha<2\gamma$. Then integrating in time by Remark \ref{rem:gxi} and (\ref{tl-sigma}) we obtain
\begin{align*}
\|\tilde{w}^{(c)}_{q+1}+\tilde{w}^{(t)}_{q+1}\|_{L^2_{((2\sigma _{q-1})\wedge T_L,T_L]}L^2}\leq\frac 14M_0(M_L^{1/2}\delta_{q+1}^{1/2}+\gamma_{q+1}^{1/2})(M_L^{1/2}-2\sigma_{q-1})^{1/2}.
\end{align*}
Similarly by (\ref{sigma-sigma}), (\ref{4.17*})-(\ref{4.19*}) and Remark \ref{rem:gxi} we obtain
\begin{align*}
\|\tilde{w}^{(c)}_{q+1}+\tilde{w}^{(t)}_{q+1}\|_{L^2_{(\sigma_{q+1}\wedge T_L,(2\sigma _{q-1})\wedge T_L]}L^2}\leq \frac14M_0((M_L+qA)^{1/2}+\gamma_{q+1}^{1/2})(2\sigma_{q-1})^{1/2}.
\end{align*}
Then by the above estimates and (\ref{wp21}), (\ref{wp22}) we obtain
\begin{align}
\|w_{q+1}\|_{L^2_{((2\sigma _{q-1})\wedge T_L,T_L]}L^2}
\leq \frac{3}{4}M_0(\delta_{q+1}^{1/2}M_L^{1/2}+\gamma_{q+1}^{1/2})(M_L^{1/2}-2\sigma_{q-1})^{1/2},\notag\\
\|w_{q+1}\|_{L^2_{(\sigma_{q+1}\wedge T_L,(2\sigma _{q-1})\wedge T_L]}L^2}\leq \frac34M_0((M_L+qA)^{1/2}+\gamma_{q+1}^{1/2})(2\sigma_{q-1})^{1/2}.\label{4.20}
\end{align}
With these bounds, we have all in hand to complete the estimate of $v_{q+1}^2$. We split the details into several subsections.
\subsubsection{Proof of (\ref{3.15})}(\ref{4.20}) together with (\ref{3.9}) yields
\begin{align*}
\|v_{q+1}^2&-v_q^2\|_{L^2_{({(2\sigma _{q-1})}\wedge T_L,T_L]}L^2}\\
&\leq \|v_{q+1}^2-v_l\|_{L^2_{({(2\sigma _{q-1})}\wedge T_L,T_L]}L^2}+\|v_{l}-v_q^2\|_{L^2_{({(2\sigma _{q-1})}\wedge T_L,T_L]}L^2}\notag\\
&\leq\|w_{q+1}\|_{L^2_{({(2\sigma _{q-1})}\wedge T_L,T_L]}L^2}+l\|v_q^2\|_{C_{t,x}^1}(T_L-{(2\sigma _{q-1})}\wedge T_L)^{1/2}\\
&\leq\frac34M_0(M_L^{1/2}\delta_{q+1}^{1/2}+\gamma_{q+1}^{1/2})(M_L^{1/2}-2\sigma_{q-1})^{1/2}+l\lambda_q^4M_L^{1/2}(M_L^{1/2}-2\sigma_{q-1})^{1/2} \\
&\leq M_0(M_L^{1/2}\delta_{q+1}^{1/2}+\gamma_{q+1}^{1/2})(M_L^{1/2}-2\sigma_{q-1})^{1/2},
\end{align*}
where we used the conditions on the parameters to deduce $-\frac{3}{2}\alpha+\frac4b<-\beta$. And similarly we obtain
\begin{align*}
\|v_{q+1}^2-v_q^2\|_{L^2_{(\sigma_{q+1}\wedge T_L,(2\sigma _{q-1})\wedge T_L]}L^2}
&\leq M_0((M_L+qA)^{1/2}+\gamma_{q+1}^{1/2})(2\sigma_{q-1})^{1/2}.
\end{align*}

For $t\in[t_{q+1},\sigma_{q+1}\wedge T_L]$ we have $\chi(t)=0$  which implies $w_{q+1}=0$. By (\ref{3.8}) we obtain $ v_q^2(t)=0$ and hence $ v_l(t)=0$,
 which implies $v_{q+1}^2(t)=w_{q+1}(t)+v_l(t)=0$.
Hence (\ref{3.15}) follows.

\subsubsection{Proof of (\ref{3.8}) on the level $q+1$}
By (\ref{3.15}) we have
\begin{align*}
\|v_{q+1}^2-v_q^2\|_{L^2_{[0,T_L]}L^2}
&\leq M_0(M_L^{1/2}\delta_{q+1}^{1/2}+\gamma_{q+1}^{1/2})(M_L^{1/2}-2\sigma_{q-1})^{1/2}+M_0((M_L+qA)^{1/2}+\gamma_{q+1}^{1/2})(2\sigma_{q-1})^{1/2}\notag\\
&\leq M_0M_L^{3/4}\delta_{q+1}^{1/2}+\sqrt{2}M_0M_L^{1/4}\gamma_{q+1}^{1/2}+M_0(M_L+qA)^{1/2}(2\sigma_{q-1})^{1/2},
\end{align*}
where we used $ \sqrt{a}+\sqrt{b}\leq \sqrt{2(a+b)}$ for $a,b\geq0$. Thus together with (\ref{para2}) we obtain
\begin{align*}
\|v_{q+1}^2\|_{L^2_{[0,T_L]}L^2}&\leq M_0(M_L^{3/4}\sum_{m=1}^{q+1}\delta_{m}^{1/2}+\sqrt{2}M_L^{1/4}\sum_{m=1}^{q+1}\gamma_{m}^{1/2})+\sqrt{2}M_0(M_L+A)^{1/2}\sum_{m=1}^{q}(m\sigma_{m-1})^{1/2}\\
&\leq M_0M_L^{3/4}+\sqrt{2} M_0M_L^{1/4}(K^{1/2}+1)+17M_0(M_L+A)^{1/2},
\end{align*}
where we roughly bounded $\sigma_{q}\leq 2^{-q}$ by (\ref{para1}) and thus $\sum_{m=1}^{q}(m\sigma_{m-1})^{1/2}\leq\sum_{m=1}^{q}m2^{-\frac{m-1}2}<17
$.
And together with the last line in (\ref{3.15}) we obtain (\ref{3.8}) on the level $q+1$.

\subsubsection{Proof of (\ref{3.14}).} By (\ref{4.16*})-(\ref{4.19*}) we have for any $m\in(1,\infty)$
\begin{align}
\|w_{q+1}-\tilde{w}_{q+1}^{(o)}\|_{C_{T_L}L^{m}}\lesssim\lambda_{q+1}^{34\alpha+(1/2-1/m)(2-8\gamma)+8\gamma},  \label{w165}
\end{align}
where we used the condition $M_L+qA+K\leq l^{-1}$. Then by (\ref{3.9}) we have
\begin{align}
\|v_{q+1}^2-v_q^2\|_{C_{T_L}L^{p}}&\lesssim\|w_{q+1}-\tilde{w}_{q+1}^{(o)}\|_{C_{T_L}L^{p}}+\|\tilde{w}_{q+1}^{(o)}\|_{C_{T_L}L^{p}}+\|v_l-v_q^2\|_{C_{T_L}L^{p}}\notag\\
&\lesssim \lambda_{q+1}^{34\alpha+(1/2-1/p)(2-8\gamma)+8\gamma}+\lambda_{q+1}^{46\alpha-2\gamma}+l\lambda_q^4M_L^{1/2}\notag\\
&\lesssim (\lambda_{q+1}^{-\gamma}+\lambda_{q+1}^{-\alpha})M_L^{1/2}\leq  M_L^{1/2}\delta_{q+1}^{1/2},\label{4.20*}
\end{align}
where we used the conditions on the parameters to deduce  $(1/2-1/p)(2-8\gamma)<-10\gamma$,   $46\alpha<\gamma$ and $l\lambda_{q}^4\leq\lambda_{q+1}^{-\alpha}$ in the last second step, and $\gamma>\alpha>\beta$ in the last step.
Then (\ref{3.14}) holds and (\ref{3.10}) at the level $q+1$ follows.

\subsubsection{Proof of (\ref{3.9}) at the level $q+1$}
Next we estimate $C_{t,x}^1$-norm,
by (\ref{4.8}), (\ref{wp+wc}) and (\ref{gwnp}) and (\ref{2.4}) for $t\in[0,T_L]$
\begin{align}
\|w_{q+1}^{(p)}+w_{q+1}^{(c)}\|_{C_{t,x}^1}&=\|\sigma^{-1}\sum_{\xi\in\Lambda}\nabla^{\perp}[a_{(\xi)}g_{(\xi)}\Psi_{(\xi)}]\|_{C_{t,x}^1}\lesssim\sigma^{-1}\sum_{\xi\in\Lambda}\|a_{(\xi)}\|_{C_{t,x}^2}\|\nabla^{\perp}[g_{(\xi)}\Psi_{(\xi)}]\|_{C_{t,x}^1}\notag\\
&\lesssim ((M_L+qA)^{1/2}+\gamma_{q+1}^{1/2})l^{-22}(\eta^{3/2}\sigma+\sigma\mu\eta^{1/2}+\theta\eta\mu)(\nu\mu)^{1/2}\lesssim\lambda_{q+1}^{45\alpha+3+17\gamma},\label{4.21}
\end{align}
where we used (\ref{2.4}) and (\ref{B.10}) to bound $\Psi_{(\xi)}$. For the last inequality we used the condition $M_L+qA+K\leq l^{-1}$.

By (\ref{gwnp}) and (\ref{parth}) we have for $t\in[0,T_L]$
\begin{align}
\|w_{q+1}^{(o)}\|_{C_{t,x}^1}&\lesssim \sigma^{-1}\sum_{\xi\in\Lambda}\|\partial_th_{(\xi)}\|_{L^\infty}\left(\|\div(a_{(\xi)}^2\xi\otimes\xi)\|_{C_tW^{1+\alpha,\infty}}+\|\div(a_{(\xi)}^2\xi\otimes\xi)\|_{C_t^1W^{\alpha,\infty}}\right)\notag\\
&\lesssim ((M_L+qA)^{1/2}+\gamma_{q+1}^{1/2})^2l^{-31}\eta\lesssim\lambda_{q+1}^{64\alpha+16\gamma},\label{4.23}
\end{align}
where we used the condition $M_L+qA+K\leq l^{-1}$.

By (\ref{4.8}), (\ref{gwnp}), (\ref{2.4}) and (\ref{B.9}) we have for $t\in[0,T_L]$
\begin{align}
\|w_{q+1}^{(a)}\|_{C_{t,x}^1}&\lesssim\sigma\theta^{-1}\sum_{\xi\in\Lambda}\|g_{(\xi)}\|_{W^{1,\infty}}\|a_{(\xi)}^2\|_{C_{t,x}^{1+\alpha}}\left(\|\ |W_{(\xi)}|^2\xi\|_{C_tW^{1+\alpha,\infty}}+\|\partial_t|W_{(\xi)}|^2\xi\|_{C_tW^{\alpha,\infty}}\right)\notag\\
&\lesssim ((M_L+qA)^{1/2}+\gamma_{q+1}^{1/2})^2l^{-24}\sigma^2\theta^{-1}\eta^{3/2}(\sigma\nu\mu^2+\nu\mu^2\theta\eta^{1/2})(\sigma\mu)^\alpha\lesssim\lambda_{q+1}^{52\alpha+3+28\gamma}.\label{4.24}
\end{align}
Here we used the condition $M_L+qA+K\leq l^{-1}$ and we have a extra $\alpha$ since $\mathbb{P}_{\rm{H}}\mathbb{P}_{\neq0}$ is not a bounded operator on $C^0$. In particular, we see that the fact that the time derivative of $\chi$ behaves like $\sigma^{-1}_{q+1}\lesssim l^{-1}$ does not pose any problems as the $C_{t,x}^0$-norms of $\tilde{w}_{q+1}^{(p)},\tilde{w}_{q+1}^{(c)},\tilde{w}_{q+1}^{(t)}$ always contain smaller powers of $l^{-1}$.
So we obtain (\ref{3.9}) at the level $q+1$ by (\ref{4.21})-(\ref{4.24}): for $t\in[0,T_L]$
\begin{align}
\|v_{q+1}^2\|_{C_{t,x}^1}&\leq \|v_{l}\|_{C_{t,x}^1}+\|w_{q+1}\|_{C_{t,x}^1}\leq M_L^{1/2}\lambda_q^4+\lambda_{q+1}^{64\alpha+3+28\gamma}\notag\\
&\leq M_L^{1/2}\lambda_{q+1}^\alpha+ \lambda_{q+1}^{64\alpha+3+28\gamma} \leq  M_L^{1/2}\lambda_{q+1}^{4}, \notag
\end{align}
where we used the conditions $\alpha<\frac{1}{128},\gamma<\frac{1}{56}$.

\subsubsection{Proof of (\ref{3.16})}
We conclude this part with further $W^{1,m}$-norm for $m\in(1,\infty)$. By (\ref{4.8}), (\ref{wp+wc}) and (\ref{2.4}) we have for $t\in[0,T_L]$
\begin{align}
\|\tilde{w}_{q+1}^{(p)}(t)+\tilde{w}_{q+1}^{(c)}(t)\|_{W^{1,m}}&\lesssim \sigma^{-1}\sum_{\xi\in\Lambda}\|\nabla^{\perp}[a_{(\xi)}g_{(\xi)}\Psi_{(\xi)}](t)\|_{W^{1,m}}\notag\\
&\lesssim\sigma^{-1}\sum_{\xi\in\Lambda}\|\nabla^{\perp}a_{(\xi)}\|_{C_{t,x}^1}|g_{(\xi)}(t)|\|\nabla^{\perp}\Psi_{(\xi)}\|_{C_tW^{1,m}}\notag\\
&\lesssim ((M_L+qA)^{1/2}+\gamma_{q+1}^{1/2})l^{-22}\sigma\mu(\nu\mu)^{1/2-1/m}\sum_{\xi\in\Lambda}|g_{(\xi)}(t)|
,\notag\\
&\lesssim\lambda_{q+1}^{45\alpha+(1/2-1/m)(2-8\gamma)+1+2\gamma}\sum_{\xi\in\Lambda}|g_{(\xi)}(t)|,\label{4.25}
\end{align}
where we used the condition $M_L+qA+K\leq l^{-1}$. By (\ref{4.8}) and (\ref{hlinfty}) we have
\begin{align}
\|\tilde{w}_{q+1}^{(o)}(t)\|_{W^{1,m}}&\lesssim \sigma^{-1}\sum_{\xi\in\Lambda}\|h_{(\xi)}\|_{L^\infty}\|\div(a_{(\xi)}^2\xi\otimes\xi)\|_{C_tW^{1,m}}\notag\\
&\lesssim ((M_L+qA)^{1/2}+\gamma_{q+1}^{1/2})^2\sigma^{-1}\lambda_{q+1}^{60\alpha}\lesssim\lambda_{q+1}^{62\alpha-2\gamma},\label{4.27}
\end{align}
and
\begin{align}
\|\tilde{w}_{q+1}^{(a)}(t)\|_{W^{1,m}}&\lesssim\sigma\theta^{-1}\sum_{\xi\in\Lambda}|g_{(\xi)}|\|a_{(\xi)}\|_{C_{t,x}^0}\|a_{(\xi)}\|_{C_{t,x}^1}\|W_{(\xi)}\|_{C_tW^{1,2m}}\|W_{(\xi)}\|_{C_tL^{2m}}\notag\\
&\lesssim ((M_L+qA)^{1/2}+\gamma_{q+1}^{1/2})^2\lambda_{q+1}^{46\alpha}\sigma\theta^{-1}\sigma\mu(\nu\mu)^{1-1/m}\sum_{\xi\in\Lambda}|g_{(\xi)}(t)|\notag\\
&\lesssim\lambda_{q+1}^{48\alpha+(1-1/m)(2-8\gamma)-\gamma}\sum_{\xi\in\Lambda}|g_{(\xi)}(t)|.\label{4.28}
\end{align}

Moreover together with (\ref{gwnp}) we have
\begin{align}
\|{w}_{q+1}-\tilde{w}_{q+1}^{(o)}\|_{C_{T_L}W^{1,6/5}}\lesssim
\lambda_{q+1}^{48\alpha+\frac{1}{3}+\frac{38}{3}\gamma}\lesssim\lambda_{q+1}^{\frac{1}{3}+13\gamma},\notag
\end{align}
where we used the conditions to deduce $48\alpha<\frac{1}{3}\gamma$.

By (\ref{w165}) we have
\begin{align*}
    \|{w}_{q+1}-\tilde{w}_{q+1}^{(o)}\|_{C_{T_L}L^{6/5}}\lesssim \lambda_{q+1}^{34\alpha-\frac{2}{3}+\frac{32}{3}\gamma}\lesssim\lambda_{q+1}^{-\frac{2}{3}+11\gamma},
\end{align*}
where we used the conditions to deduce $34\alpha<\frac{1}{3}\gamma$.

Then together the bounds above and   (\ref{3.9}), (\ref{4.18*}) and (\ref{4.27}) we obtain by interpolation that
\begin{align}
\|v_{q+1}^2&-v_q^2\|_{C_{T_L}W^{1/2,6/5}}\notag\\
&\leq \|v_l^2-v_q^2\|_{C_{T_L}W^{1/2,6/5}}+ \|{w}_{q+1}-\tilde{w}_{q+1}^{(o)}\|_{C_{T_L}W^{1/2,6/5}}+\|\tilde{w}_{q+1}^{(o)}\|_{C_{T_L}W^{1/2,6/5}}\notag\\
&\lesssim l^{1/2}\|v_q^2\|_{C_{[0,T_L],x}^1}+\|{w}_{q+1}-\tilde{w}_{q+1}^{(o)}\|_{C_{T_L}L^{6/5}}^{1/2}\|{w}_{q+1}-\tilde{w}_{q+1}^{(o)}\|_{C_{T_L}W^{1,6/5}}^{1/2}+\|\tilde{w}_{q+1}^{(o)}\|_{C_{T_L}W^{1/2,6/5}}\notag\\
&\lesssim l^{1/2}\lambda_q^4M_L^{1/2}+\lambda_{q+1}^{-\frac{1}{6}+12\gamma}+\lambda_{q+1}^{55\alpha-2\gamma}\leq M_L^{1/2}\lambda_{q+1}^{-\alpha}\leq M_L^{1/2}\delta_{q+1}^{1/2},\notag
\end{align}
which implies (\ref{3.16}). Here we used the condition on the parameters to deduce $l^{1/2}\lambda_q^4<\lambda_{q+1}^{-\alpha}$, $-\frac{1}{6}+12\gamma<-\alpha<-\beta$, $55\alpha-2\gamma<-\alpha<-\beta$.

\subsection{Proof of (\ref{3.17})}\label{subsub5}
For $t\in(2\wedge T_L,T_L]$ we have $\chi(t)=1$ by (\ref{def:chi}) . Then we find
\begin{align}
| \|v_{q+1}^2&\|_{L^2_{(2\wedge T_L,T_L]}L^2}^2-\|v_q^2\|_{L^2_{(2\wedge T_L,T_L]}L^2}^2-2\gamma_{q+1}(T_L-2\wedge T_L)|\notag\\
&\leq| \|w_{q+1}^{(p)}\|_{L^2_{(2\wedge T_L,T_L]}L^2}^2-2\gamma_{q+1}(T_L-2\wedge T_L)|+\|w_{q+1}^{(c)}+w_{q+1}^{(t)}\|_{L^2_{(2\wedge T_L,T_L]}L^2}^2\notag\\
&+2\|v_l(w_{q+1}^{(c)}+w_{q+1}^{(t)})\|_{L^1_{(2\wedge T_L,T_L]}L^1}+2\|v_lw_{q+1}^{(p)}\|_{L^1_{(2\wedge T_L,T_L]}L^1}\notag\\
&+2\|w_{q+1}^{(p)}(w_{q+1}^{(c)}+w_{q+1}^{(t)})\|_{L^1_{(2\wedge T_L,T_L]}L^1}+|\|v_l\|_{L^2_{(2\wedge T_L,T_L]}L^2}^2-\|v_q^2\|_{L^2_{(2\wedge T_L,T_L]}L^2}^2|.\label{4.29}
\end{align}
We use (\ref{4.10}) and the fact that $\mathring{R}_l$ is traceless to deduce
\begin{align*}
|w_{q+1}^{(p)}|^2-2\gamma_{q+1}&=4\sqrt{l^2+|\mathring{R}_l|^2}
+\sum_{\xi\in\Lambda}a_{(\xi)}^2g_{(\xi)}^2\mathbb{P}_{\neq0}|W_{(\xi)}|^2+\sum_{\xi\in\Lambda}a_{(\xi)}^2(g_{(\xi)}^2-1),
\end{align*}
hence
\begin{align}
| \|{w}_{q+1}^{(p)}\|_{L^2_{(2\wedge T_L,T_L]}L^2}^2&-2\gamma_{q+1}(T_L-2\wedge T_L)|\leq 4lT_L+4\|\mathring{R}_l\|_{L^1_{(2\wedge T_L,T_L]}L^1}\notag\\
&+\sum_{\xi\in\Lambda}\|a_{(\xi)}^2g_{(\xi)}^2\mathbb{P}_{\neq0}|W_{(\xi)}|^2\|_{L^1_{(2\wedge T_L,T_L]}L^1}+\sum_{\xi\in\Lambda}\left|\int_{2\wedge T_L}^{T_L}(g_{(\xi)}^2-1)\|a_{(\xi)}\|_{L^2}^2\dif t\right|\notag.
\end{align}
Now we estimate each term separately. Using (\ref{3.11}) and $\supp \phi_l\subset[0,l]$, $(2\wedge T_L,T_L]\subset((2\sigma_{q-1})\wedge T_L,T_L]$ we have
$$\|\mathring{R}_l\|_{L^1_{(2\wedge T_L,T_L]}L^1}\leq M_L\delta_{q+1}.$$
Moreover,
$$4lT_L\leq4\lambda_{q+1}^{-\frac{3\alpha}{2}}T_L\leq\frac{1}{8}\delta_{q+1}M_L,$$
which requires $2\beta<\frac{3\alpha}{2}$ and we choose $a$ large enough to absorb the constant.

We obverse that $W_{(\xi)}$ is $(\mathbb{T}/\sigma)^2$-periodic so
$$\mathbb{P}_{\neq0}(W_{(\xi)}\otimes W_{(\xi)})=\mathbb{P}_{\geq\frac{\sigma}{2}}(W_{(\xi)}\otimes W_{(\xi)}),$$
where $\mathbb{P}_{\geq r}=\Id-\mathbb{P}_{<r}$, and $\mathbb{P}_{<r}$ denotes the Fourier multiplier operator, which projects a function
onto its Fourier frequencies $< r$ in absolute value. By \cite[Proposition C.1]{BDLPS15}, (\ref{4.8}) and (\ref{2.4}) we obtain for $t\in(2\wedge T_L,T_L]$
\begin{align}
\sum_{\xi\in\Lambda}|\int_{\mathbb{T}^2}a_{(\xi)}^2\mathbb{P}_{\neq0}|W_{(\xi)}|^2
\dif x|&=\sum_{\xi\in\Lambda}|\int_{\mathbb{T}^2}a_{(\xi)}^2\mathbb{P}_{\geq\frac{\sigma}{2}}|W_{(\xi)}|^2
\dif x|=\sum_{\xi\in\Lambda}|\int_{\mathbb{T}^2}|\nabla|^{\gamma^*}a_{(\xi)}^2|\nabla|^{-\gamma^*}\mathbb{P}_{\geq\frac{\sigma}{2}}|W_{(\xi)}|^2
\dif x|\notag\\
&\lesssim\|a_{(\xi)}^2\|_{C^{\gamma^*}}\sigma^{-\gamma^*}\||W_{(\xi)}|^2\|_{L^2}\notag\lesssim(M_L+K)\lambda_{q+1}^{(14\gamma^*+32)\alpha-2\gamma \gamma^*+1-4\gamma}\notag\\
&\lesssim\lambda_{q+1}^{(14\gamma^*+34)\alpha-1-4\gamma}\lesssim \lambda_{q+1}^{-4\gamma},\notag
\end{align}
where we recall $\gamma^*=\frac{1}{\gamma}$. Here we used the conditions on the the parameters to deduce $M_L+K\leq l^{-1}$ and
$(14\gamma^*+34)\alpha=\frac{14\gamma^*+34}{145\gamma^*}<1$, which togther with $2\gamma>\beta$ imply that
$$\sum_{\xi\in\Lambda}\|a_{(\xi)}^2g_{(\xi)}^2\mathbb{P}_{\neq0}|W_{(\xi)}|^2\|_{L^1_{(2\wedge T_L,T_L]}L^1}\leq \frac{1}{8}\delta_{q+1}M_L,$$
where we chose $a$ large enough to absorb the universal constant.

For the last term, to apply Theorem \ref{ihiot} it is possible to find $n_1,n_2\in\mathbb{N}_0$ such that $\frac{n_1}\sigma<2\wedge T_L\leq\frac{n_1+1}\sigma, \frac{n_2}{\sigma} \leq T_L<\frac{ n_2+1}\sigma$ and define $a_{(\xi)}(t)=a_{(\xi)}(T_L)$ for $t\in(T_L, \frac{n_2+1}\sigma]$. Note that $g^2_{(\xi)}-1$ is mean-zero on $[0,1]$. Together with (\ref{4.8}), (\ref{gwnp}) we apply Theorem \ref{ihiot} and Remark \ref{rem:gxi} to obtain
\begin{align*}
&\ \ \sum_{\xi\in\Lambda}|\int_{2\wedge T_L}^{T_L}(g_{(\xi)}^2-1)\|a_{(\xi)}\|_{L^2}^2\dif t|\notag\\
&\leq\sum_{\xi\in\Lambda}\left(|\int_{\frac{n_1+1}\sigma}^{\frac{n_2}\sigma}(g_{(\xi)}^2-1)\|a_{(\xi)}\|_{L^2}^2\dif t|+\int_{2\wedge T_L}^{\frac{n_1+1}\sigma}(g_{(\xi)}^2+1)\|a_{(\xi)}\|_{L^2}^2\dif t+\int_{\frac{n_2}\sigma}^{T_L}(g_{(\xi)}^2+1)\|a_{(\xi)}\|_{L^2}^2\dif t\right)\\
&\lesssim \frac{n_2-n_1-1}{\sigma^2}\|g_{(\xi)}^2+1\|_{L_{[0,1]}^1}\left\|\|a_{(\xi)}\|_{L^2}^2\right\|_{C^{0,1}_{[\frac{n_1+1}{\sigma},\frac{n_2}{\sigma}]}}+\|g_{(\xi)}^2+1\|_{L^1([\frac{n_1}\sigma,\frac{n_1+1}\sigma]\cup[\frac{n_2}\sigma,\frac{n_2+1}\sigma])}\|a_{(\xi)}\|_{C_{[0,T_L],x}^0}^2\\
&\lesssim \sigma^{-1}T_Ll^{-23}((M_L+qA)^{1/2}+\gamma_{q+1}^{1/2})^2+\sigma^{-1} l^{-16}((M_L+qA)^{1/2}+\gamma_{q+1}^{1/2})^2\\
&\lesssim\lambda_{q+1}^{48\alpha-2\gamma}M_L^{1/2}+\lambda_{q+1}^{34\alpha-2\gamma}\leq \frac1{8}\delta_{q+1}M_L,
\notag
\end{align*}
where we used the conditions on the parameters to deduce $M_L+qA+K\leq l^{-1}$, $48\alpha-2\gamma<-2\beta $.

Go back to (\ref{4.29}), we control the remaining terms as follows. Using (\ref{4.17})-(\ref{4.19}) and Remark \ref{rem:gxi} we obtain
\begin{align}
\|w_{q+1}^{(c)}+{w}_{q+1}^{(t)}\|_{L^2_{(2\wedge T_L,T_L]}L^2}^2
&\lesssim \sum_{\xi\in\Lambda}\|g_{(\xi)}^2+1\|_{L^1_{(2\wedge T_L,T_L]}}(\lambda_{q+1}^{31\alpha-8\gamma}
+\lambda_{q+1}^{48\alpha-2\gamma}+\lambda_{q+1}^{34\alpha-7\gamma})^2\notag\\
&\lesssim (T_L-2\wedge T_L+\frac{2}{\sigma})\lambda_{q+1}^{96\alpha-4\gamma}\leq\frac{1}{8}\delta_{q+1}M_L,\notag
\end{align}
where we used the conditions on the parameters to deduce $M_L+K\leq l^{-1}$, and in the last inequality we used $T_L-2\wedge T_L+\frac{2}{\sigma}\lesssim M_L^{1/2}$, $48\alpha-2\gamma<-\beta$, and chose $a$ large enough to absorb the constant.

By (\ref{3.8}), (\ref{wp21}) and (\ref{4.17*})-(\ref{4.19*}) we obtain
\begin{align}
2\|v_l&({w}_{q+1}^{(c)}+{w}_{q+1}^{(t)})\|_{L^1_{(2\wedge T_L,T_L]}L^1}
+2\|{w}_{q+1}^{(p)}({w}_{q+1}^{(c)}+{w}_{q+1}^{(t)})\|_{L^1_{(2\wedge T_L,T_L]}L^1}\notag\\
&\lesssim (\|{v}_l\|_{L^2_{(2\wedge T_L,T_L]}L^2}+\|{w}_{q+1}^{(p)}\|_{L^2_{(2\wedge T_L,T_L]}L^2})\|{w}_{q+1}^{(c)}+{w}_{q+1}^{(t)}\|_{L^2_{(2\wedge T_L,T_L]}L^2}\notag\\
&\lesssim M_0(M_L^{3/4}+A^{1/2}+M_L^{1/4}K^{1/2})(\lambda_{q+1}^{31\alpha-8\gamma}
+\lambda_{q+1}^{48\alpha-2\gamma}+\lambda_{q+1}^{34\alpha-7\gamma})M_L^{1/4}\notag\\
&\lesssim \lambda_{q+1}^{49\alpha-2\gamma}M_L^{1/2}\leq\frac{1}{4}\delta_{q+1}M_L.\notag
\end{align}
where we used $\|{v}_l\|_{L^2_{(2\wedge T_L,T_L]}L^2}\leq\|v_q^2\|_{L^2_{[t_q,T_L]}L^2}$ by $l\leq 2$, and we used the conditions on the parameters to deduce $M_L+A+K\leq l^{-1}$ and $49\alpha-2\gamma<-2\beta$ and we possibly increase $a$ to absorb the constant.

By (\ref{3.9}), (\ref{4.16}) and (\ref{gwnp}) we obtain
\begin{align}
2\|v_l{w}_{q+1}^{(p)}\|_{L^1_{(2\wedge T_L,T_L]}L^1}&\lesssim\|v_l\|_{C_{[0,T_L],x}^0}\|w_{q+1}^{(p)}\|_{L^1_{(2\wedge T_L,T_L]}L^1}\notag\\
&\lesssim\lambda_q^{4}M_L^{1/2}(M_L^{1/2}\delta_{q+1}^{1/2}+K^{1/2})\lambda_{q+1}^{16\alpha-1-4\gamma}M_L^{1/2}\leq\frac{1}{8}\delta_{q+1}M_L,\notag
\end{align}
where we used the conditions on the parameters to deduce $M_L+K\leq l^{-1}$, $\lambda_q^4\leq\lambda_{q+1}^{\alpha}$ and $18\alpha-1-4\gamma<-2\beta$.

For the last term, by (\ref{3.8}) and (\ref{3.9})
\begin{align}
|\|v_l&\|_{L^2_{(2\wedge T_L,T_L]}L^2}^2-\|v_q^2\|_{L^2_{(2\wedge T_L,T_L]}L^2}^2|\notag\\
&\leq\|v_l-v_q^2\|_{L^2_{(2\wedge T_L,T_L]}L^2}(\|v_l\|_{L^2_{(2\wedge T_L,T_L]}L^2}+\|v_q^2\|_{L^2_{(2\wedge T_L,T_L]}L^2})\notag\\
&\leq l\|v_q^2\|_{C_{[0,T_L],x}^1}M_L^{1/4}\|v_q^2\|_{L^2_{[t_q,T_L]}L^2}\notag\\
&\lesssim l\lambda_{q}^{4}M_L^{1/4}\times M_0(M_L^{3/4}+A^{1/2}+M_L^{1/4}K^{1/2})\notag\\
&\lesssim\lambda_{q+1}^{-\alpha}M_0M_L^{1/2}(M_L+A+K)^{1/2}\leq\frac{1}{8}\delta_{q+1}M_L.\notag
\end{align}
where we used $\|{v}_l\|_{L^2_{(2\wedge T_L,T_L]}L^2}\leq\|v_q^2\|_{L^2_{[t_q,T_L]}L^2}$ by $l\leq 2$, and used the conditions on the parameters to deduce $l\lambda_{q}^{4}\leq\lambda_{q+1}^{-\alpha}$ and $M_0(M_L+A+K)^{1/2}\ll\lambda_{q+1}^{\alpha-2\beta}$ and then chose $a$ large enough to absorb the constant.

Combining the above estimates (\ref{3.17}) follows.

\subsection{Definition of the Reynolds stress}\label{dotrs}
Combine (\ref{3.6}) at the level $q+1$, (\ref{4.5}) and the definition of $w_{q+1}$ we get
\begin{align}
\div\mathring{R}_{q+1}&-\nabla p_{q+1}^2+\nabla p_l\notag\\
&=-\Delta w_{q+1}+\partial_t(\tilde{w}_{q+1}^{(p)}+\tilde{w}_{q+1}^{(c)})+\div((v_l+v_{q+1}^1)\otimes w_{q+1}+w_{q+1}\otimes (v_l+v_{q+1}^1))\notag\\
&\ \ \ \ \ \ \ \ \ \ \ \ \ \ \ \ \ \ \ \ \ \ \ \ \ \ \ \ \ \ \ \ \ \ \ \ \ \ \ \ \ \ \ \ \ \ \ \ \ \ \ \ \ \ \  \ \ \ \ \ \ \ \ \ \ \ \ \ (:=\div R_{lin}+\nabla p_{lin})\notag\\
&+\div((\tilde{w}_{q+1}^{(c)}+\tilde{w}_{q+1}^{(t)})\otimes w_{q+1}+\tilde{w}_{q+1}^{(p)}\otimes (\tilde{w}_{q+1}^{(c)}+\tilde{w}_{q+1}^{(t)}))\ (:=\div R_{cor}+\nabla p_{cor})\notag\\
&+\partial_t\tilde{w}_{q+1}^{(t)}+\div(\tilde{w}_{q+1}^{(p)}\otimes \tilde{w}_{q+1}^{(p)}+\mathring{R}_l)\ (:=\div R_{osc}+\nabla p_{osc})\notag\\
&+\div((v_{q+1}-v_q)\otimes\Delta_{\leq R}z+\Delta_{\leq R}z\otimes(v_{q+1}-v_q)\notag\\
&\ \ \ \ \ \ \ \ \ \ \  \ \ \ +(v_{q+1}-v_q)\succcurlyeq\!\!\!\!\!\!\!\bigcirc\Delta_{>R}z+\Delta_{>R}z\preccurlyeq\!\!\!\!\!\!\!\bigcirc(v_{q+1}-v_q))(:=\div R_{com1}+\nabla p_{com1})\notag\\
&+\div((v_l+v_{q+1}^1)\otimes(v_l+v_{q+1}^1)-(v_q^1+v_q^2)\otimes(v_q^1+v_q^2))(:=\div R_{com2}+\nabla p_{com2})\notag\\
&+\div((v_{q+1}-v_q)\otimes z^{in}+z^{in}\otimes(v_{q+1}-v_q))(:=\div R_{com3}+\nabla p_{com3} )\notag\\
&+\div{R}_{com},\notag
\end{align}
where $R_{com}$ is introduced in Section \ref{moll}.

Applying the inverse divergence operator $\mathcal{R}$, by (\ref{rdeltav}) we could define
\begin{align*}
R_{lin}:&=-\nabla w_{q+1}-\nabla^Tw_{q+1}+ \mathcal{R}(\partial_t(\tilde{w}_{q+1}^{(p)}+\tilde{w}_{q+1}^{(c)}))+(v_l+v_{q+1}^1)\mathring{\otimes} w_{q+1}+w_{q+1}\mathring{\otimes} (v_l+v_{q+1}^1),\\
R_{cor}:&=(\tilde{w}_{q+1}^{(c)}+\tilde{w}_{q+1}^{(t)})\mathring{\otimes} w_{q+1}+\tilde{w}_{q+1}^{(p)}\mathring{\otimes} (\tilde{w}_{q+1}^{(c)}+\tilde{w}_{q+1}^{(t)}),\\
R_{com2}:&=(v_l+v_{q+1}^1)\mathring{\otimes}(v_l+v_{q+1}^1)-(v_q^1+v_q^2)\mathring{\otimes}(v_q^1+v_q^2),\\
R_{com3}:&=(v_{q+1}-v_q)\mathring{\otimes} z^{in}+z^{in}\mathring{\otimes}(v_{q+1}-v_q),
\end{align*}
$R_{com1}$ is the trace-free part of the matrix
\begin{align*}
(v_{q+1}-v_q){\otimes}\Delta_{\leq R}z+\Delta_{\leq R}z{\otimes}(v_{q+1}-v_q) +(v_{q+1}-v_q)\ {\succcurlyeq\!\!\!\!\!\!\!\bigcirc}\Delta_{>R}z+\Delta_{>R}z\ {\preccurlyeq\!\!\!\!\!\!\!\bigcirc}(v_{q+1}-v_q).
\end{align*}

In order to define the remaining oscillation error from the forth line, we apply (\ref{4.10}), (\ref{t=0+a}) and direct calculation to obtain
\begin{align}
\partial_t&\tilde{w}_{q+1}^{(t)}+\div(\tilde{w}_{q+1}^{(p)}\otimes \tilde{w}_{q+1}^{(p)}+\mathring{R}_l)\notag\\
=&\sum_{\xi\in\Lambda}\chi^2g^2_{(\xi)}{\mathbb{P}_{\neq0}\left(\nabla (a_{(\xi)}^2)\mathbb{P}_{\neq0}(W_{(\xi)}\otimes W_{(\xi)})\right)}+\chi^2(\Id-\mathbb{P}_{\rm{H}})\sum_{\xi\in\Lambda}\mathbb{P}_{\neq0}\left(a_{(\xi)}^2g^2_{(\xi)}\div(W_{(\xi)}\otimes W_{(\xi)})\right)\notag\\
&\ \ \ \ \ \ \ \ \ \ \ \ \ \ \ \ \ \ \ \ \ \ \ \ \ \ \ \ \ \ \ \ \ \ \ \ \ \ \ \ \ \ \ \ \ \ \ \ \ \ \ \ \ \ \ \ \ \ \ \ \ \ \ \ \ \ \ \ \ \ \ \ \ \ \ \ \ \ \ \ \ \ \ \ (:=\chi^2R_{osc,x}+\chi^2\nabla p_{osc,x})\notag\\
&+\chi^2\partial_t{w}_{q+1}^{(a)}+\chi^2\mathbb{P}_{\rm{H}}\sum_{\xi\in\Lambda}\mathbb{P}_{\neq0}\left(a_{(\xi)}^2g^2_{(\xi)}\div(W_{(\xi)}\otimes W_{(\xi)})\right)(:=\chi^2R_{osc,a}+\chi^2\nabla p_{osc,a})\notag\\
&+\chi^2\partial_t{w}_{q+1}^{(o)}+\sum_{\xi\in\Lambda}\chi^2(g^2_{(\xi)}-1)\div(a_{(\xi)}^2\xi\otimes\xi){(:=\chi^2R_{osc,t}+\chi^2\nabla p_{osc,t})}\notag\\
&+(\chi^2)'w_{q+1}^{(t)}+\div( (1-\chi^2)\mathring{R}_{l})+\nabla \rho.\notag
\end{align}

Applying the inverse divergence operator $\mathcal{B}$, we could define
\begin{align*}
    R_{osc,x}:&=\sum_{\xi\in\Lambda}g^2_{(\xi)}\mathcal{B}\left(\nabla a_{(\xi)}^2,\mathbb{P}_{\neq0}(W_{(\xi)}\otimes W_{(\xi)})\right),\\
    \nabla p_{osc,x}:&=(\Id-\mathbb{P}_{\rm{H}})\sum_{\xi\in\Lambda}\mathbb{P}_{\neq0}\left(a_{(\xi)}^2g^2_{(\xi)}\div(W_{(\xi)}\otimes W_{(\xi)})\right).
\end{align*}
For $R_{osc,a}$ by (\ref{4.13}), (\ref{B.9}) we obtain
\begin{align*}
    \partial_t{w}_{q+1}^{(a)}&+\mathbb{P}_{\rm{H}}\sum_{\xi\in\Lambda}\mathbb{P}_{\neq0}\left(a_{(\xi)}^2g^2_{(\xi)}\div(W_{(\xi)}\otimes W_{(\xi)})\right)=-\theta^{-1}\sigma\mathbb{P}_{\rm{H}}\mathbb{P}_{\neq0}\sum_{\xi\in\Lambda}\partial_t(a_{(\xi)}^2g_{(\xi)})|W_{(\xi)}|^2\xi.
\end{align*}
Thus we define
\begin{align*}
  R_{osc,a}:&=-\theta^{-1}\sigma\sum_{\xi\in\Lambda}\mathcal{B}\left(\partial_t(a_{(\xi)}^2g_{(\xi)}),|W_{(\xi)}|^2\xi\right),\\
  \nabla p_{osc,a}:&=\theta^{-1}\sigma(\Id-\mathbb{P}_{\rm{H}})\mathbb{P}_{\neq0}\sum_{\xi\in\Lambda}\partial_t(a_{(\xi)}^2g_{(\xi)})|W_{(\xi)}|^2\xi.
\end{align*}
For $R_{osc,t}$ by (\ref{4.12}) and (\ref{parth}) we obtain
\begin{align*}
\partial_t{w}_{q+1}^{(o)}&+\sum_{\xi\in\Lambda}(g^2_{(\xi)}-1)\div(a_{(\xi)}^2\xi\otimes\xi)\\
&=-\sigma^{-1}\mathbb{P}_{\rm{H}}\mathbb{P}_{\neq0}\sum_{\xi\in\Lambda}h_{(\xi)}\div\partial_t(a_{(\xi)}^2\xi\otimes\xi)+(\Id-\mathbb{P}_{\rm{H}})\mathbb{P}_{\neq0}\sum_{\xi\in\Lambda}(g^2_{(\xi)}-1)\div(a_{(\xi)}^2\xi\otimes\xi),\notag
\end{align*}
then we define
\begin{align*}
R_{osc,t}:&=-\sigma^{-1}\sum_{\xi\in\Lambda}\mathcal{R}\left(h_{(\xi)}\div\partial_t(a_{(\xi)}^2\xi\otimes\xi)\right),\\
\nabla p_{osc,t}:&=\sigma^{-1}(\Id-\mathbb{P}_{\rm{H}})\mathbb{P}_{\neq0}\sum_{\xi\in\Lambda}h_{(\xi)}\div\partial_t(a_{(\xi)}^2\xi\otimes\xi)+(\Id-\mathbb{P}_{\rm{H}})\mathbb{P}_{\neq0}\sum_{\xi\in\Lambda}(g^2_{(\xi)}-1)\div(a_{(\xi)}^2\xi\otimes\xi).
\end{align*}
Together with the above we could define $$R_{osc}:=\chi^2(R_{osc,x}+R_{osc,a}+R_{osc,t})+(\chi^2)'\mathcal{R}w_{q+1}^{(t)}+ (1-\chi^2)\mathring{R}_{l}.$$ Finally we define the Reynolds stress on the level $q+1$ by
$$\mathring{R}_{q+1}=R_{lin}+R_{cor}+R_{osc}+
R_{com}+R_{com1}+R_{com2}+R_{com3}.$$
We observe that by construction, $\mathring{R}_{q+1}$ is $(\mathcal{F}_t)_{t\geq 0}$-adapted.

\subsection{Verification of the inductive estimate for $\mathring{R}_{q+1}$}\label{votiefr}
To conclude the proof of Proposition \ref{prop:3.1} we shall verify (\ref{3.11}) at the level $q+1$.
In order to establish the iterative estimate, we distinguish four cases corresponding to the four time intervals.

Case \uppercase\expandafter{\romannumeral1}. Let $t\in(\sigma _q\wedge T_L,T_L]$. Note that if $T_L\leq \sigma_q$ then there is nothing to estimate here, hence
we assume that $\sigma_q < T_L$ and $t\in(\sigma_q, T_L]$.
Using (\ref{para1}) and  (\ref{para1}) we conclude that $\chi(t)=1$ on the interval $(\sigma_q,T_L]$. We will estimate each term present in the definition of $\mathring{R}_{q+1}$ in Section \ref{dotrs} separately.

For the linear error $R_{lin}$, we estimate each term separately. By (\ref{wp+wc}) and (\ref{B.10}) we obtain
\begin{align*}
    \partial_t(w_{q+1}^{(p)}+w_{q+1}^{(c)})&=\sigma^{-1}\sum_{\xi\in\Lambda}\nabla^{\perp}\partial_t(a_{(\xi)}g_{(\xi)}\Psi_{(\xi)})\\
    &=\sigma^{-1}\sum_{\xi\in\Lambda}\nabla^{\perp}[\partial_t(a_{(\xi)}g_{(\xi)})\Psi_{(\xi)}]+\sigma^{-2}\theta\sum_{\xi\in\Lambda}\nabla^{\perp}[a_{(\xi)}g_{(\xi)}^2(\xi\cdot\nabla)\Psi_{(\xi)}].
\end{align*}
Now by the fact that $\mathcal{R}\nabla^{\perp}$ is $L^p\to L^p$ bounded for any $p>1$, together with (\ref{4.8}), (\ref{2.4}) and (\ref{2.4*}) we obtain for $\epsilon>0$ small enough
\begin{align}
&\|\mathcal{R}(\partial_t(w_{q+1}^{(p)}+w_{q+1}^{(c)}))(t)\|_{L^1}\lesssim\|\mathcal{R}(\partial_t(w_{q+1}^{(p)}+w_{q+1}^{(c)}))(t)\|_{L^{1+\epsilon}}\notag\\
&\lesssim \sigma^{-1}\sum_{\xi\in\Lambda}\|\partial_t(a_{(\xi)}g_{(\xi)})\Psi_{(\xi)}(t)\|_{L^{1+\epsilon}}+\sigma^{-2}\theta\sum_{\xi\in\Lambda}\|a_{(\xi)}g_{(\xi)}^2(\xi\cdot\nabla)\Psi_{(\xi)}\|_{L^{1+\epsilon}}\notag\\
&\lesssim\sum_{\xi\in\Lambda}\|a_{(\xi)}\|_{C_{t,x}^1}\left(\sigma^{-1}\|\Psi_{(\xi)}\|_{C_tL^{1+\epsilon}}(|\partial_tg_{(\xi)}(t)|+|g_{(\xi)}(t)|)+\sigma^{-2}\theta\|(\xi\cdot\nabla)\Psi_{(\xi)}\|_{C_tL^{1+\epsilon}}g^2_{(\xi)}(t)\right)\notag\\
&\lesssim ((M_L+qA)^{\frac12}+K^{\frac12})l^{-15}\sum_{\xi\in\Lambda}(\sigma^{-1
}\nu^{-\frac12}\mu^{-\frac32}(|\partial_tg_{(\xi)}(t)|+|g_{(\xi)}(t)|)+\sigma^{-1}\theta\nu^{\frac12}\mu^{-\frac32}g^2_{(\xi)}(t)
)(\nu\mu)^{\epsilon}\notag\\
&\lesssim \lambda_{q+1}^{31\alpha-2+2\gamma+2\epsilon} \sum_{\xi\in\Lambda}(|\partial_tg_{(\xi)}(t)|+|g_{(\xi)}(t)|)+\lambda_{q+1}^{31\alpha-\gamma+2\epsilon}\sum_{\xi\in\Lambda}g^2_{(\xi)}(t)\notag\\
&\lesssim \lambda_{q+1}^{-1}\sum_{\xi\in\Lambda} (|\partial_tg_{(\xi)}(t)|+|g_{(\xi)}(t)|)+\lambda_{q+1}^{-\frac12\gamma}\sum_{\xi\in\Lambda}g^2_{(\xi)}(t),\notag
\end{align}
where we chose $2\epsilon<\alpha$ and used the conditions on the parameters to deduce $32\alpha<\frac12\gamma,\frac{5}{2}\gamma<1$, $M_L+qA+K\leq l^{-1}$.

By (\ref{4.25})-(\ref{4.28})
we obtain for $\epsilon>0$ small enough
\begin{align}
\|&\nabla w_{q+1}(t)+\nabla^T w_{q+1}(t)\|_{L^{1}}\lesssim\|w_{q+1}\|_{C_tW^{1,1+\epsilon}}\notag\\
&\lesssim
 \lambda_{q+1}^{48\alpha+6\gamma+2\epsilon}\sum_{\xi\in\Lambda}|g_{(\xi)}(t)|+\lambda_{q+1}^{62\alpha-2\gamma}\lesssim \lambda_{q+1}^{7\gamma}\sum_{\xi\in\Lambda}|g_{(\xi)}(t)|+\lambda_{q+1}^{-\gamma}.\notag
\end{align}
where we $2\epsilon<\alpha$ and used the conditions {on parameters} to deduce $62\alpha<\gamma$.

By (\ref{3.9}), (\ref{4.4}) and (\ref{4.20*}) we obtain
\begin{align}
\|(v_l+v_{q+1}^1)&\mathring{\otimes} w_{q+1}+w_{q+1}\mathring{\otimes}(v_l+v_{q+1}^1) \|_{L^1}\lesssim
\|v_l+v_{q+1}^1\|_{C_t L^\infty}\|w_{q+1}\|_{C_tL^p}\notag\\
&\lesssim\lambda_{q+1}^{-\gamma}M_L^{1/2}(\lambda_q^{4}+L^2\lambda_{q+1}^{\alpha})\lesssim \lambda_{q+1}^{-\frac12\gamma}M_L^{1/2},\notag
\end{align}
where we used the conditions {on parameters} to deduce $L^2\leq l^{-1}$ and $\lambda_q^4<\lambda_{q+1}^{\alpha}$ and $6\alpha<\gamma$.

Together with all the estimates above we obtain
\begin{align}
\|{R}_{lin}(t)\|_{L^1}\lesssim\lambda_{q+1}^{-1} \sum_{\xi\in\Lambda}|\partial_tg_{(\xi)}(t)|+\lambda_{q+1}^{7\gamma}\sum_{\xi\in\Lambda}|g_{(\xi)}(t)|+\lambda_{q+1}^{-\frac12\gamma}\sum_{\xi\in\Lambda}g^2_{(\xi)}(t)+\lambda_{q+1}^{-\frac12\gamma}M_L^{1/2}.\label{rlin}
\end{align}

The corrector error ${R}_{cor}$ is estimated using (\ref{4.17*})-(\ref{4.19*}) for $t\in(\sigma_q,T_L]$
\begin{align}
\|{R}_{cor}(t)\|_{L^1}&\leq\|w_{q+1}^{(c)}(t)+w_{q+1}^{(t)}(t)\|_{L^{2}}(\|w_{q+1}(t)\|_{L^{2}}
+\|w_{q+1}^{(p)}(t)\|_{L^{2}})\notag\\
&\lesssim
\lambda_{q+1}^{48\alpha-2\gamma}(\sum_{\xi\in\Lambda}|g_{(\xi)}(t)|+1)(\|w_{q+1}(t)\|_{L^{2}}
+\|w_{q+1}^{(p)}(t)\|_{L^{2}})\notag\\
&\lesssim
\lambda_{q+1}^{-\gamma}(\sum_{\xi\in\Lambda}|g_{(\xi)}(t)|+1)(\|w_{q+1}(t)\|_{L^{2}}
+\|w_{q+1}^{(p)}(t)\|_{L^{2}}),\label{2.30}
\end{align}
where we used the conditions on the paramters to deduce $M_L+qA+K\leq l^{-1},48\alpha<\gamma$.

We continue with oscillation error ${R}_{osc}$. By (\ref{def:chi}) we have $\chi(t)=1$ on the interval $(\sigma_q,T_L]$, which implies $R_{osc}=R_{osc,x}+R_{osc,a}+R_{osc,t}$. In order to bound the first term, we note that $W_{(\xi)}$ is $(\mathbb{T}/\sigma)^2$-periodic. Then we apply {(\ref{bb1}),} (\ref{bb}), (\ref{4.8}) and (\ref{2.4}) to obtain for all $t\in[0,T_L]$
\begin{align}
\|{R}_{osc,x}(t)\|_{L^1}&\lesssim\sum_{\xi\in\Lambda}\|\mathcal{B}(\nabla a_{(\xi)}^2,\mathbb{P}_{\neq0}(W_{(\xi)}\otimes W_{(\xi)}))\|_{C_tL^1}g^2_{(\xi)}(t)\notag\\
&\lesssim \sum_{\xi\in\Lambda}\|\nabla a_{(\xi)}^2\|_{C_tC^1}\|\mathcal{R}(W_{(\xi)}\otimes W_{(\xi)})\|_{C_tL^{1}}g^2_{(\xi)}(t)\notag\\
&\lesssim\sum_{\xi\in\Lambda} l^{-30}(M_L+qA+\gamma_{q+1})\sigma^{-1}\|W_{(\xi)}\otimes W_{(\xi)}\|_{C_tL^{1}}g^2_{(\xi)}(t)\notag\\
&\lesssim \lambda_{q+1}^{62\alpha-2\gamma}\sum_{\xi\in\Lambda}g^2_{(\xi)}(t)\lesssim\lambda_{q+1}^{-\gamma}\sum_{\xi\in\Lambda}g^2_{(\xi)}(t),\notag
\end{align}
where we used the conditions on parameters to deduce $M_L+qA+K\leq l^{-1},62\alpha<\gamma$.

We apply {(\ref{bb1})} (\ref{bb}), (\ref{4.8})
and (\ref{2.4}) to obtain for all $t\in[0,T_L]$
\begin{align}
\|{R}_{osc,a}(t)\|_{L^1}&\lesssim \theta^{-1}\sigma\sum_{\xi\in\Lambda}\|\mathcal{B}\left(\partial_t(a_{(\xi)}^2g_{(\xi)}),|W_{(\xi)}|^2\xi\right)\|_{L^{1}}\notag\\
&\lesssim\theta^{-1}\sigma\sum_{\xi\in\Lambda}\|\partial_t(a_{(\xi)}^2g_{(\xi)})\|_{C^1}\|\mathcal{R}(|W_{(\xi)}|^2\xi)\|_{C_tL^{1}}\notag\\
&\lesssim \theta^{-1}\|a_{(\xi)}^2\|_{C_{t,x}^2}\sum_{\xi\in\Lambda}(| g_{(\xi)}(t)|+|\partial_t g_{(\xi)}(t)|)\|W_{(\xi)}\|_{C_tL^{2}}^2
,\notag\\
&\lesssim \theta^{-1}l^{-30}(M_L+qA+\gamma_{q+1})
\sum_{\xi\in\Lambda}(| g_{(\xi)}(t)|+|\partial_t g_{(\xi)}(t)|)\notag\\
&\lesssim \lambda_{q+1}^{62\alpha-1-5\gamma}\sum_{\xi\in\Lambda}(| g_{(\xi)}(t)|+|\partial_t g_{(\xi)}(t)|)\lesssim \lambda_{q+1}^{-1}\sum_{\xi\in\Lambda}(| g_{(\xi)}(t)|+|\partial_t g_{(\xi)}(t)|),\notag
\end{align}
where we used the conditions on parameters to deduce $M_L+qA+K\leq l^{-1}$ and $62\alpha<5\gamma$.

By (\ref{4.8}) and (\ref{hlinfty}) we have for all $t\in[0,T_L]$
\begin{align}
\|{R}_{osc,t}(t)\|_{L^1}&\lesssim\sigma^{-1}\sum_{\xi\in\Lambda}\|a_{(\xi)}^2\|_{C_{t,x}^1}\|h_{(\xi)}\|_{L^\infty}\lesssim\sigma^{-1}l^{-23}(M_L+qA+\gamma_{q+1})
\lesssim\lambda_{q+1}^{-\gamma},\notag
\end{align}
where we used the conditions to deduce $M_L+qA+K\leq l^{-1},48\alpha<\gamma$. Therefore we obtain
\begin{align}
   \|{R}_{osc}(t)\|_{L^1} \lesssim\lambda_{q+1}^{-\gamma}(\sum_{\xi\in\Lambda}g^2_{(\xi)}(t)+1)+\lambda_{q+1}^{-1}\sum_{\xi\in\Lambda}(| g_{(\xi)}(t)|+|\partial_t g_{(\xi)}(t)|).
\end{align}

For ${R}_{com1}$, by Lemma \ref{lem:2.3} with $1-2\kappa-\kappa_0>0$, (\ref{3.16}), (\ref{3.12}) and (\ref{est:deltaz}) we obtain for $t\in[0,T_L]$
\begin{align}
\|{R}_{com1}(t)\|_{L^1}&\lesssim
(\|v_{q+1}^2-v_q^2\|_{C_tW^{1/2,6/5}}+\|v_{q+1}^1-v_q^1\|_{C_tB^{1-\kappa-\kappa_0}_{p,\infty}})(\|\Delta_{\leq R}z\|_{C_tL^\infty}+\|z\|_{C_tC^{-\kappa}})\notag\\
&\lesssim L^{5/2}M_L^{1/2}(\lambda_{q+1}^{-\alpha}+L^{5/4}\lambda_{q}^{-\frac{\alpha}{8}\kappa_0})\lesssim M_L(\lambda_{q+1}^{-\alpha}+L^{5/4}\lambda_{q}^{-\frac{\alpha}{8}\kappa_0}),\label{est:com1}
\end{align}
where we used the condition on the parameters to deduce $M_L\geq L^5$.

For ${R}_{com2}$, by
(\ref{3.9}), (\ref{3.12}) and (\ref{3.13}) we obtain for $t\in[0,T_L]$
\begin{align}
\|{R}_{com2}(t)\|_{L^1}&\lesssim (\|v_{q+1}^1(t)\|_{L^2}+\|v_{q}^1(t)\|_{L^2}+\|v_q^2(t)\|_{L^2}+\|v_l(t)\|_{L^2})\notag\\
&\ \ \ \ \times(\|v_{q+1}^1-v_q^1(t)\|_{C_tL^2}+\|v_l-v_q^2\|_{C_tL^2})\notag\\
&\lesssim (M_L^{1/2}+\|v_q^2(t)\|_{L^2}+\|v_l(t)\|_{L^2}) M_L^{1/2}(\lambda_{q+1}^{-\alpha}+L^{5/4}\lambda_{q}^{-\frac{\alpha}{8}\kappa_0}+l\lambda_q^4)\notag\\
&\lesssim (M_L^{1/2}+\|v_q^2(t)\|_{L^2}+\|v_l(t)\|_{L^2})M_L^{1/2}(\lambda_{q+1}^{-\alpha}+L^{5/4}\lambda_{q}^{-\frac{\alpha}{8}\kappa_0}),\label{est:rcom2}
\end{align}
where we used the conditions on the parameters to deduce $l\lambda_q^4\leq\lambda_{q+1}^{-\alpha}$.

For $t\in[0,T_L]$, it holds that by H\"older's inequality
\begin{align}
\|{R}_{com3}(t)\|_{L^1}\lesssim(\|v_{q+1}^1-v_q^1\|_{C_tL^p}+\|v_{q+1}^2-v_q^2\|_{C_tL^{p}})\|z^{in}(t)\|_{L^{\frac{p}{p-1}}}.\notag
\end{align}
By \cite[Lemma 9]{DV} and (\ref{omegan}) we have
$$\|z^{in}(t)\|_{L^{\frac{p}{p-1}}}\lesssim(1+t^{\frac{p-2}{p}})\|u_0\|_{L^p}\lesssim(1+t^{\frac{p-2}{p}})N,$$
which together with (\ref{3.14}), (\ref{3.12}) imply
\begin{align}
\|{R}_{com3}(t)\|_{L^1}&\lesssim (\lambda_{q+1}^{-\alpha}+L^{5/4}\lambda_{q}^{-\frac{\alpha}{8}\kappa_0})M_L^{1/2}
(1+t^{\frac{p-2}{p}})N.\label{est:com3}
\end{align}

Finally let us now consider each term in ${R}_{com}$  separately, which is defined in Section \ref{moll}. For $I_1=(v_q^1+v_q^2+z^{in})\otimes\Delta_{\leq R}z+\Delta_{\leq R}z\otimes(v_q^1+v_q^2+z^{in})$ and $t\in(\sigma_{q+1},T_L]$ we have
\begin{align}
\|(I_1-I_1*_x\phi_l*_t\varphi_l)(t)\|_{L^1}&\lesssim
l^{\frac{\kappa_0}{2}}(\|v_q\|_{C_{[\frac{\sigma_{q+1}}2,t]}^{\kappa_0/2}L^2}+\|z^{in}\|_{C_{[\frac{\sigma_{q+1}}2,t]}^{\kappa_0/2}L^2})\|\Delta_{\leq R}z\|_{C_t^{\kappa_0/2}L^\infty}\notag\\
&+l^{\frac{\kappa_0}{2}}(\|v_q\|_{C_{{[\frac{\sigma_{q+1}}2,t]}}H^{1/3}}+\|z^{in}\|_{C_{[\frac{\sigma_{q+1}}2,t]}H^1})\|\Delta_{\leq R}z\|_{C_tC^{\kappa_0}}\notag\\
&\lesssim l^{\frac{\kappa_0}{2}}(\lambda_q^4M_L^{1/2}+M_L^{1/2}\lambda_q^\alpha+\sigma_{q+1}^{-1/p} M_L^{1/2})L^{\frac92},
\notag\\
&\lesssim l^{\frac{\kappa_0}{2}}(\lambda_q^4+\sigma_{q+1}^{-1/p})M_L,\notag
\end{align}
where we used
(\ref{deltaf}) and (\ref{3.7}) to deduce
$$\|\Delta_{\leq R}z\|_{C_tC^{\kappa_0}}\lesssim 2^{(\kappa+\kappa_0) R}\|z\|_{C_tC^{-\kappa}}\lesssim 2^{2\kappa_0 R}\|z\|_{C_tC^{-\kappa}}\lesssim L^{\frac{17}{4}},$$
$$\|\Delta_{\leq R}z\|_{C_t^{\kappa_0/2}L^\infty}\lesssim\|\Delta_{\leq R}z\|_{C_t^{\kappa_0/2}C^{\kappa_0/8}}\lesssim 2^{(\frac{9}8\kappa_0+\kappa) R}\|z\|_{C_t^{\kappa_0/2}C^{-\kappa-\kappa_0}}\lesssim 2^{\frac{17}8\kappa_0 R}\|z\|_{C_t^{\kappa_0/2}C^{-\kappa-\kappa_0}}\lesssim L^{\frac{9}{2}},$$
and by \cite[Lemma 9]{DV} and interpolation we have
\begin{align}
\|z^{in}\|_{C_{[\frac{\sigma_{q+1}}2,t]}^{\kappa_0/2}L^2}+\|z^{in}\|_{C_{[\frac{\sigma_{q+1}}2,t]}^{\frac{1}{4}(1-2\kappa-\kappa_0)}H^{1/2-\kappa_0/2}}&\lesssim\|z^{in}\|_{C_{[\frac{\sigma_{q+1}}2,t]}H^1}+\|z^{in}\|_{C_{[\frac{\sigma_{q+1}}2,t]}^{1/2}L^2}\notag\\
&\lesssim\sigma_{q+1}^{-\frac1p}\|u_0\|_{L^p}\lesssim\sigma_{q+1}^{-\frac1p}M_L^{1/2}.   \label{zin}
\end{align}
Here we bounded $v_q^1$ using equation (\ref{4.4}) and $v_q^2$ using equation (\ref{3.9}). Additionally, we used $M_L\geq L^{9}$ in the last inequality.

For $I_2=(v_q^1+v_q^2+z^{in})\otimes(v_q^1+v_q^2+z^{in})$ and $t\in(\sigma_{q+1},T_L]$, by (\ref{3.9}), (\ref{3.13}), (\ref{4.4}) and (\ref{zin}) we have
\begin{align}
\|(I_2-I_2&*_x\phi_l*_t\varphi_l)(t)\|_{L^1}\notag\\
&\lesssim
l^{\frac{\kappa_0}{2}}(\|v_q\|_{C_{[\frac{\sigma_{q+1}}2,t]}^{\kappa_0/2}L^2}+\|z^{in}\|_{C_{[\frac{\sigma_{q+1}}2,t]}^{\kappa_0/2}L^2}+\|v_q\|_{C_{{[\frac{\sigma_{q+1}}2,t]}}H^{1/3}}+\|z^{in}\|_{C_{[\frac{\sigma_{q+1}}2,t]}H^1})\notag\\
&\ \ \ \ \ \ \ \times(\|v_q^1\|_{C_tL^2}+\|v_q^2\|_{C_{t,x}^1}+\|z^{in}\|_{C_{[\frac{\sigma_{q+1}}2,t]}L^2})\notag\\
&\lesssim l^{\frac{\kappa_0}{2}}(\lambda_q^4M_L^{1/2}+M_L^{1/2}\lambda_q^\alpha+\sigma_{q+1}^{-1/p} M_L^{1/2})(\lambda_q^4M_L^{1/2}+\sigma_{q+1}^{-1/p}M_L^{1/2})\notag\\
&\lesssim l^{\frac{\kappa_0}{2}}(\lambda_q^4+\sigma_{q+1}^{-1/p})^2M_L.\notag
\end{align}

For $I_3=(v_q^1+v_q^2+z^{in})\succcurlyeq\!\!\!\!\!\!\!\bigcirc\Delta_{>R}z+\Delta_{>R}z\preccurlyeq\!\!\!\!\!\!\!\bigcirc(v_q^1+v_q^2+z^{in})$ and $t\in(\sigma_{q+1},T_L]$ by   (\ref{3.7}), (\ref{3.9}), (\ref{3.13}), (\ref{4.2}) and (\ref{zin}) we obtain
\begin{align}
&\|(I_3-I_3*_x\phi_l*_t\varphi_l)(t)\|_{L^1}\notag\\
&\lesssim l^{\frac{1}{4}(1-2\kappa-\kappa_0)}\Big{(}(\|v_q^1\|_{C_tB^{1-\kappa-\kappa_0}_{p,\infty}}+\|v_q^2\|_{C_{[\frac{\sigma_{q+1}}2,t]}H^{1-\kappa-\kappa_0}}+\|z^{in}\|_{C_{[\frac{\sigma_{q+1}}2,t]}H^{1-\kappa-\kappa_0}})\|z\|_{C_t C^{-\kappa}}\notag\\
&\ \ +(\|v_q^1\|_{C_t^{\frac{1}{4}(1-2\kappa-\kappa_0)}B^{1/2-\kappa_0/2}_{p,\infty}}+\|v_q^2\|_{C_{[\frac{\sigma_{q+1}}2,t]}^{\frac{1}{4}(1-2\kappa-\kappa_0)}H^{1/2-\kappa_0/2}}+\|z^{in}\|_{C_{[\frac{\sigma_{q+1}}2,t]}^{\frac{1}{4}(1-2\kappa-\kappa_0)}H^{1/2-\kappa_0/2}})\|z\|_{C_t C^{-\kappa}}\notag\\
&\ \  +(\|v_q^1\|_{C_tB^{1-\kappa-\kappa_0}_{p,\infty}}+\|v_q^2\|_{C_{[\frac{\sigma_{q+1}}2,t]}H^{1-\kappa-\kappa_0}}+\|z^{in}\|_{C_{[\frac{\sigma_{q+1}}2,t]}H^{1-\kappa-\kappa_0}})\|z\|_{C_t^{\frac{1}{4}(1-2\kappa-\kappa_0)} C^{-\frac12+\frac{\kappa_0}2}}\Big{)}\notag\\
&\lesssim l^{\frac{1}{4}(1-3\kappa_0)}(\lambda_q^4M_L^{1/2}+M_L^{1/2}+\sigma_{q+1}^{-1/p} M_L^{1/2})L^{1/4}\lesssim l^{\frac{1}{4}(1-3\kappa_0)}(\lambda_q^4+\sigma_{q+1}^{-1/p})M_L.\notag
\end{align}

Therefore, we have for $t\in(\sigma_{q+1},T_L]$
\begin{align}
\|{R}_{com}(t)\|_{L^1}
&\lesssim l^{\frac{\kappa_0}{2}}(\lambda_q^4+\sigma_{q+1}^{-1/p})^2M_L+l^{\frac{1}{4}(1-3\kappa_0)}(\lambda_q^4+\sigma_{q+1}^{-1/p})M_L\lesssim l^{\kappa_1} M_L,\label{2.37}
\end{align}
where we recall $\kappa_1=\frac{\kappa_0}{4}\wedge(\frac{1}{6}-\frac{\kappa_0}2)$ and in the first inequality we used the conditions on the parameters to deduce $l^{\kappa_1/2}\lambda_q^4\leq1, l^{\kappa_1/2}\leq \sigma_{q+1}^{1/p}$ by choosing $\alpha b>\frac{16}{3\kappa_1}$ and $\kappa_1\alpha>\frac{8}{3}\beta$.

Summarizing (\ref{rlin})-(\ref{est:com3}) and (\ref{2.37}) we obtain
\begin{align*}\|\mathring{R}_{q+1}(t)\|_{L^1}&\lesssim \sum_{\xi\in\Lambda}(\lambda_{q+1}^{-1} |\partial_t g_{(\xi)}(t)|+\lambda_{q+1}^{7\gamma}|g_{(\xi)}(t)|+\lambda_{q+1}^{-\frac12\gamma}g^2_{(\xi)}(t)+\lambda_{q+1}^{-\frac12\gamma}M_L^{1/2})\\&+\lambda_{q+1}^{-\gamma}(\sum_{\xi\in\Lambda}|g_{(\xi)}(t)|+1)(\|w_{q+1}(t)\|_{L^{2}}+\|{w}_{q+1}^{(p)}(t)\|_{L^{2}}),\\
&+(M_L^{1/2}+\|v_q^2(t)\|_{L^2}+\|v_l(t)\|_{L^2}+Nt^{\frac{p-2}{p}})M_L^{1/2}(\lambda_{q+1}^{-\alpha}+L^{5/4}\lambda_{q}^{-\frac{\alpha}{8}\kappa_0})+l^{\kappa_1} M_L.\end{align*}
Then we integrate over the interval $(\sigma_q,T_L]$ to obtain
\begin{align*}
\|\mathring{R}_{q+1}\|_{L^1_{(\sigma_q,T_L]}L^1}&\lesssim \sum_{\xi\in\Lambda}(\lambda_{q+1}^{-1}\|\partial_tg_{(\xi)}\|_{L_{(\sigma_q,T_L]}^{1}}+\lambda_{q+1}^{7\gamma}\|g_{(\xi)}\|_{L_{(\sigma_q,T_L]}^{1}}+\lambda_{q+1}^{-\frac12\gamma}\|g_{(\xi)}\|_{L_{(\sigma_q,T_L]}^{2}}^2+\lambda_{q+1}^{-\frac12\gamma}M_L^{1/2}L)\\
&+\lambda_{q+1}^{-\gamma}(\sum_{\xi\in\Lambda}\|g_{(\xi)}\|_{L_{(\sigma_{q},T_L]}^{2}}+L)(\|w_{q+1}\|_{L^2_{(\sigma_{q+1},T_L]}L^{2}}
+\|\tilde{w}_{q+1}^{(p)}\|_{L^2_{(\sigma_{q+1},T_L]}L^{2}})\\
&+(M_L^{1/2}L+\|v_q^2\|_{L^2_{[t_q,T_L]}L^2}+
NL^{2-2/p})M_L^{1/2}(\lambda_{q+1}^{-\alpha}+L^{5/4}\lambda_{q}^{-\frac{\alpha}{8}\kappa_0})+l^{\kappa_1} M_LL,
\end{align*}
where we used $\|v_q^2(t)\|_{L^2_{(\sigma_q,T_L]}L^2}+\|v_l(t)\|_{L^2_{(\sigma_q,T_L]}L^2}\lesssim\|v_q^2\|_{L^2_{[t_q,T_L]}L^2}$ by $l\leq \sigma_q$.
 Then we bound $g_{(\xi)}$ by (\ref{gwnp}), Remark \ref{rem:gxi} and that $T_L-\sigma_q+\frac{2}{\sigma}\leq L+2\leq2L$. Together with (\ref{3.8}), (\ref{wpq+1qtl}) and (\ref{4.20}) we obtain
\begin{align*}
\|\mathring{R}_{q+1}\|_{L^1_{(\sigma_q,T_L]}L^1}
&\lesssim (\lambda_{q+1}^{-1+10\gamma} +\lambda_{q+1}^{-\gamma}+\lambda_{q+1}^{-\frac12\gamma})M_L+\lambda_{q+1}^{-\gamma}((M_L+qA)^{1/2}+K^{1/2})M_L^{1/2},\\
&\ \ \ \ +(M_L^{3/4}+A^{1/2}+M_L^{1/4}K^{1/2})M_L^{1/2}(\lambda_{q+1}^{-\alpha}+L^{5/4}\lambda_{q}^{-\frac{\alpha}{8}\kappa_0})+l^{\kappa_1} M_LL\\
&\lesssim \lambda_{q+1}^{-\frac12\gamma}M_L+(\lambda_{q+1}^{-\alpha}+L^{5/4}\lambda_{q}^{-\frac{\alpha}{8}\kappa_0})M_L^{3/4}(M_L^{1/2}+(qA)^{1/2}+K^{1/2})+l^{\kappa_1} M_LL\\
&\leq M_L^{1/2}(M_L^{1/2}-\sigma_q^{2-2/p})\delta_{q+2}\leq M_L\delta_{q+2}.
\end{align*}
where we used the conditions  on the parameters to deduce  $\frac12\gamma>\alpha>\max\{\frac{16\beta b^2}{\kappa_0},\frac{4\beta b}{3\kappa_1}\},$  $(\lambda_{q+1}^{-\alpha+2\beta b}+L^{5/4}\lambda_{q}^{-\frac{\alpha}{8}\kappa_0+2\beta b^2})(M_L+qA+K)^{1/2}+ L\lambda_{q+1}^{-\frac{3}2\kappa_1\alpha+2\beta b}\ll1$, then we chose $a$ large enough to absorb the universal constant. Here the extra $\sigma_q^{2-2/p}\leq1$ in the last second inequality is used to absorb the bound in Case \uppercase\expandafter{\romannumeral2}. Thus we obtained the first line of (\ref{3.11}) at the level $q+1$.

Case \uppercase\expandafter{\romannumeral2}.
Let $t\in(\sigma_{q+1}\wedge T_L,\sigma_q\wedge T_L]$. If $T_L\leq \sigma_{q+1}$
then there is nothing to estimate, hence we may assume $\sigma_{q+1} < T_L$ and $t\in(\sigma_{q+1}, \sigma_q \wedge T_L]$. In this case, $R_{cor}, R_{com1}, R_{com3}$ and $R_{com}$ in the definition of $\mathring{R}_{q+1}$ are similar to those in Case \uppercase\expandafter{\romannumeral1}, with the added factor of $\chi(t)\leq1$, so we can estimate them in the same way as before. For $R_{com2}$ we have $v_q^2(t)=v_l(t)=0$ by (\ref{3.8}) and $\sigma_{q+1}\leq t\leq \sigma_q\wedge T_L$. Using a similar approach as  (\ref{est:rcom2}), we obtain
\begin{align}
\|{R}_{com2}(t)\|_{L^1}
&\lesssim M_L(\lambda_{q+1}^{-\alpha}+L^{5/4}\lambda_{q}^{-\frac{\alpha}{8}\kappa_0}).\notag
\end{align}
For $R_{lin}$, there is an additional term  $\chi'(t)(w_{q+1}^{(p)}(t)+w_{q+1}^{(c)}(t))$. For $R_{osc}$, there are additional terms $(\chi^2)'\mathcal{R}w_{q+1}^{(t)}$ and $(1-\chi^2)\mathring{R}_{l}$ by definition in Section \ref{dotrs}. The rest parts can be estimated similarly. Thus we only have to consider $(1-\chi^2)\mathring{R}_l$ and
\begin{align}
{R}_{cut}:=\chi'(t)(\mathcal{R}w_{q+1}^{(p)}(t)+\mathcal{R}w_{q+1}^{(c)}(t))+(\chi^2)'(t)\mathcal{R}w_{q+1}^{(t)}(t).\notag
\end{align}
As for ${R}_{cut}$ by (\ref{4.16*})-(\ref{4.19*}) we have for $\epsilon>0$ small enough
\begin{align*}
\|{R}_{cut}(t)\|_{L^1}&\leq\|\chi'(t)(w_{q+1}^{(p)}+w_{q+1}^{(c)})\|_{L^{1+\epsilon}}+\|(\chi^2)'w_{q+1}^{(t)}\|_{L^{1+\epsilon}}\notag\\
&\lesssim \sigma_{q+1}^{-1}(\|w_{q+1}^{(p)}(t)\|_{L^{1+\epsilon}}+\|w_{q+1}^{(c)}(t)\|_{L^{1+\epsilon}}+\|w_{q+1}^{(t)}(t)\|_{L^{1+\epsilon}})\\
&\leq \sigma_{q+1}^{-1}((M_L+qA)^{1/2}+\gamma_{q+1}^{1/2})^2\lambda_{q+1}^{46\alpha-2\gamma+2\epsilon}\lesssim\lambda_{q+1}^{-\gamma},
\end{align*}
where we choose ${2\epsilon}<\alpha$ and used the conditions on the parameters to deduce $\|\chi(t)\|_{C^1}\lesssim{\sigma_{q+1}^{-1}}\leq l^{-1}$, $M_L+qA+K\leq l^{-1}$ and $51\alpha<\gamma$.

Therefore, using similar calculations as in Case \uppercase\expandafter{\romannumeral1} we obtain for $\sigma_{q+1}<t\leq\sigma_q\wedge T_L\leq1$
\begin{align*}\|\mathring{R}_{q+1}(t)\|_{L^1}&\lesssim \lambda_{q+1}^{-1} \sum_{\xi\in\Lambda}|\partial_t g_{(\xi)}(t)|+\lambda_{q+1}^{7\gamma}\sum_{\xi\in\Lambda}|g_{(\xi)}(t)|+\lambda_{q+1}^{-\frac12\gamma}(\sum_{\xi\in\Lambda}g^2_{(\xi)}(t)+M_L^{1/2})\\&+\lambda_{q+1}^{-\gamma}(\sum_{\xi\in\Lambda}|g_{(\xi)}(t)|+1)(\|w_{q+1}(t)\|_{L^{2}}+\|\tilde{w}_{q+1}^{(p)}(t)\|_{L^{2}}),\\
&+t^{\frac{p-2}{p}}M_L(\lambda_{q+1}^{-\alpha}+L^{5/4}\lambda_{q}^{-\frac{\alpha}{8}\kappa_0})+l^{\kappa_1} M_L+(1-\chi^2)\|\mathring{R}_l(t)\|_{L^1}.\end{align*}
Here we used $t^{1-2/p}\geq1$ since $1-2/p<0$.
Then by the similar calculation as before, together with (\ref{def:chi}) we obtain
\begin{align*}
\|\mathring{R}_{q+1}\|_{L^1_{(\sigma_{q+1},\sigma_q\wedge T_L]}L^1}
&\lesssim\lambda_{q+1}^{-\frac12\gamma}M_L^{1/2}\sigma_q+\lambda_{q+1}^{-\gamma}\sigma_q^{1/2}(\|w_{q+1}\|_{L^2_{(\sigma_{q+1}, T_L]}L^{2}}+\|\tilde{w}_{q+1}^{(p)}\|_{L^2_{(\sigma_{q+1},T_L]}L^{2}})\\
&+\sigma_q^{2-2/p}M_L(\lambda_{q+1}^{-\alpha}+L^{5/4}\lambda_{q}^{-\frac{\alpha}{8}\kappa_0})+l^{\kappa_1} M_L\sigma_q+\|\mathring{R}_l\|_{L^1_{(\sigma_{q+1},(2\sigma_{q+1})\wedge T_L]}L^1}\\
&\lesssim\lambda_{q+1}^{-\frac12\gamma}M_L^{1/2}\sigma_q+\lambda_{q+1}^{-\frac12\gamma}\sigma_q((M_L+qA)^{1/2}+K^{1/2})M_L^{1/4}\\
&+\sigma_q^{2-2/p}M_L(\lambda_{q+1}^{-\alpha}+L^{5/4}\lambda_{q}^{-\frac{\alpha}{8}\kappa_0})+l^{\kappa_1} M_L\sigma_q+\|\mathring{R}_l\|_{L^1_{(\sigma_{q+1},(2\sigma_{q+1})\wedge T_L]}L^1}\\
&\leq M_L^{1/2}\sigma_q^{2-2/p}\delta_{q+2}+\|\mathring{R}_l\|_{L^1_{(\sigma_{q+1},(2\sigma_{q+1})\wedge T_L]}L^1},\end{align*}
where we bound $g_{(\xi)}$ by (\ref{gwnp}), Remark \ref{rem:gxi} and note that $\sigma_q\wedge T_L-\sigma_{q+1}+\frac{2}{\sigma}\lesssim \sigma_q$. In the last second inequality we used (\ref{4.20}) and the conditions on the parameters to deduce $M_L+qA+K\leq l^{-1}$, $\lambda_{q+1}^{-\frac12\gamma}\leq\sigma_q^{\frac12}$ by $\gamma>2\beta/b$,   and in the last inequality we used $2\alpha>\gamma>4\beta b$, $\sigma_q\leq \sigma^{2-2/p}_{q}$ by $0\leq2-2/p\leq1$, $M_L^{1/2}(\lambda_{q+1}^{-\alpha+2\beta b}+L^{5/4}\lambda_{q}^{-\frac{\alpha}{8}\kappa_0+2\beta b^2})+M_L^{1/2}\lambda_{q+1}^{-\frac{3}2\kappa_1\alpha+2\beta b}\ll1$, and chose $a$ large enough to absorb the universal constant.

Case \uppercase\expandafter{\romannumeral3}. Let $t\in[0,\sigma_{q+1}\wedge T_L]$. By (\ref{3.8}) we know $v_l(t)=v_{q+1}^2(t)=0$ and
\begin{align}
\mathring{R}_{q+1}&=\mathring{R}_{l}+{R}_{com}+{R}_{com1}+{R}_{com2}+{R}_{com3}.
\end{align}

By the estimates of ${R}_{com1}$, ${R}_{com2}$ and ${R}_{com3}$ (\ref{est:com1})-(\ref{est:com3}) we obtain for $0\leq t\leq\sigma_{q+1}\wedge T_L\leq1$
$$\|{R}_{com1}(t)+{R}_{com2}(t)\|_{L^1}\lesssim M_L(\lambda_{q+1}^{-\alpha}+L^{5/4}\lambda_{q}^{-\frac{\alpha}{8}\kappa_0})\leq M_L,$$
$$\|{R}_{com3}(t)\|_{L^1}\lesssim M_L(\lambda_{q+1}^{-\alpha}+L^{5/4}\lambda_{q}^{-\frac{\alpha}{8}\kappa_0})(1+t^{\frac{p-2}{p}})\leq M_Lt^{1-\frac2p},$$
 where we used the conditions  on the parameters to deduce $t^{1-2/p}\geq1$ as $1-2/p<0$.

For ${R}_{com}$, the trace-free part of
${N}_{com}-{N}_{com}*_x\phi_l*_t\varphi_l,$ we do not use the mollification estimate. Instead, we bound $N_{com}
$ directly and have
\begin{align*}
\|\mathring{R}_{q+1}(t)\|_{L^1}
&\leq \|\mathring{R}_{l}(t)\|_{L^1}+ 2M_Lt^{1-
\frac{2}{p}}+\|N_{com}(t)\|_{L^1}+\|N_{com}*_t\varphi_l(t)\|_{L^1}.\end{align*}
For the last two terms, similar as (\ref{est:r0t+}) and (\ref{est:r0t-}),  we obtain for $0\leq t\leq\sigma_{q+1}\wedge T_L\leq1$
\begin{align*}
\|N_{com}(t)\|_{L^1}&\lesssim(\|v_q^1\|_{C_tB^{1-\kappa-\kappa_0}_{p,\infty}}+
\|z^{in}(t)\|_{W^{2/p-1,p}})(\|\Delta_{\leq R}z\|_{C_tL^\infty}+\|z\|_{C_tC^{-\kappa}})+\|v_q^1(t)+z^{in}(t)\|_{L^2}^2\\
&\lesssim L^{5}+L^{-3/2}M_L+N^2(1+t^{1-\frac2p})\leq M_Lt^{1-\frac2p},
\end{align*}
and for $-1\leq t_{q+1}\leq t<0$ and $1-\frac{2}{p}<0$
\begin{align*}
\|N_{com}(t)\|_{L^1}&\lesssim\|z^{in}(|t|)\|_{L^2}^2\lesssim N^2(1+|t|^{1-\frac2p})\leq M_L|t|^{1-\frac2p},
\end{align*}
which implies that for any $0\leq t_1<t_2\leq \sigma_{q+1}\wedge T_L$
 \begin{align*}
 \|N_{com}*_t&\varphi_l\|_{L^1_{[t_1,t_2]}L^1}\leq
\sup_{b\in[0,l]}\|N_{com}\|_{L^1_{[t_1-b,t_2-b]}L^1}\leq \frac{p}{p-1}M_L(t_2-t_1)^{2-
\frac{2}{p}}.
\end{align*}
Then by the choice of $A$  we obtain
\begin{align*}
\|\mathring{R}_{q+1}\|_{L^1_{[0,\sigma_{q+1}\wedge T_L]}L^1}
&\leq \|\mathring{R}_{l}\|_{L^1_{[0,\sigma_{q+1}\wedge T_L]}L^1}+\frac{3p}{p-1}M_L(\sigma_{q+1}\wedge T_L)^{2-
\frac{2}{p}}\\
&\leq \|\mathring{R}_{l}\|_{L^1_{[0,\sigma_{q+1}\wedge T_L]}L^1}+A(\sigma_{q+1}\wedge T_L)^{2-
\frac{2}{p}}.
\end{align*}

Then together with the bound in
{Case \uppercase\expandafter{\romannumeral1}} and
{Case \uppercase\expandafter{\romannumeral2}} we obtain
\begin{align*}
\|\mathring{R}_{q+1}\|_{L^1_{[0, T_L]}L^1}
&\leq  M_L\delta_{q+2}+\|\mathring{R}_{l}\|_{L^1_{[0,(2\sigma_{q+1})\wedge T_L]}L^1}+A\sigma_{q+1}^{2-\frac{2}{p}}\\
&\leq  M_L\delta_{q+2}+\sup_{b\in[-l,0]}\|\mathring{R}_{q}\|_{L^1_{[b,b+2\sigma_{q+1}]}L^1}+A\sigma_{q+1}^{2-\frac{2}{p}}\\
&\leq M_L\delta_{q+2}+ 2(q+2)A\sigma_{q+1}^{2-\frac{2}{p}},
\end{align*}
where in the last inequality we used the last line of (\ref{3.11}) and (\ref{para1}) to deduce $2\sigma_{q+1}\leq\sigma_q$. Thus we proved the second line of (\ref{3.11}) at the level $q+1$.

By the similar argument above we have for any $0\leq a\leq a+h\leq\sigma_{q+1}\wedge T_L$
\begin{align}
\|\mathring{R}_{q+1}\|_{L^1_{[a,a+h]}L^1}
&\leq\|\mathring{R}_{l}\|_{L^1_{[a,a+h]}L^1}+\frac{p}{p-1}M_Lh^{2-
\frac{2}{p}}+\|N_{com}\|_{L^1_{[a,a+h]}L^1}+\|N_{com}*_t\varphi_l\|_{L^1_{[a,a+h]}L^1}\notag\\
&\leq 2(q+1)A(\frac{h}{2})^{2-\frac{2}{p}}+A
h^{2-\frac{2}{p}}\leq 2(q+2)A(\frac{h}{2})^{2-\frac{2}{p}}.\label{rq+11}
\end{align}

Case \uppercase\expandafter{\romannumeral4}. Let $t\in[t_{q+1},0)$. We have $\mathring{R}_{q+1}=z^{in}\mathring{\otimes} z^{in}$
and by the similar argument as in (\ref{est:r0t-}) we obtain
\begin{align*}
\|\mathring{R}_{q+1}(t)\|_{L^1}&\leq M_L|t|^{1-2/p}.
\end{align*}
By the choice of $A$ we have for $t_{q+1}\leq a-h\leq a<0$
\begin{align}
\|\mathring{R}_{q+1}\|_{L^1_{[a-h,a]}L^1}\leq A((h-a)^{2-\frac{2}{p}}-(-a)^{2-\frac{2}{p}})\leq Ah^{2-\frac{2}{p}}.\label{rq+12}
\end{align}
Thus we could obtain the third line of (\ref{3.11}) at the level $q+1$. More specifically for any $t_{q+1}\leq a<a+h \leq {\sigma_{q+1}}\wedge T_L$, the case $a\geq0$ and the case $a+h\leq0$ are directly from (\ref{rq+11}) and (\ref{rq+12}). For $a<0<a+h$  also by (\ref{rq+11}) and (\ref{rq+12}) we obtain
\begin{align}
   \|\mathring{R}_{q+1}\|_{L^1_{[a,a+h]}L^1}&=\|\mathring{R}_{q+1}\|_{L^1_{[a,0)}L^1}+ \|\mathring{R}_{q+1}\|_{L^1_{[0,a+h]}L^1}\notag\\
   &\leq 2(q+1)A(\frac{a+h}{2})^{2-\frac{2}{p}} +A\left((-a)^{2-\frac{2}{p}}+(a+h)^{2-\frac{2}{p}}\right)\notag\\
   &\leq 2(q+1)A(\frac{h}{2})^{2-\frac{2}{p}}+2A(\frac{h}{2})^{2-\frac{2}{p}}
   = 2(q+2)A(\frac{h}{2})^{2-\frac{2}{p}}.\notag
\end{align}
We finish the proof of Proposition \ref{prop:3.1}.

\section{Proof of Theorem \ref{thm:unique}}\label{potuni}
This section is devoted to the proof of Theorem \ref{thm:unique}.  In the following we consider the equation on $[0, T]$ for simplicity.

\begin{proof}[Proof of Theorem \ref{thm:unique}]
By Proposition \ref{prop:Z}  it follows  $z\in C_TC^{-\kappa}$ for any $\kappa>0$. Then for any given $\zeta\in(0,1)$ we could find $\kappa>0$ such that $\zeta\in(\kappa,1-2\kappa)$.  We only need to prove uniqueness of solutions to (\ref{3.2}).  We assume $v_1,v_2\in C([0,T];L^2)\cap L^2([0,T];H^\zeta)$ be two solutions to (\ref{3.2}) with the same initial data. Now we define $w:=v_1-v_2$ which satisfies
\begin{align}w(t)=-\int_0^t\mathbb{P}_{\rm{H}}e^{(t-s)\Delta}\div((v_1+z)\otimes w+w\otimes (v_2+z))\dif s,\ \ w(0)=0,\label{mildv}\end{align}
where we recall $\mathbb{P}_{\rm{H}}$ is the Helmholtz projection.

As $v_1,v_2\in C_TL^2$, for any $\epsilon>0$ there exist $v_1^*,v_2^*\in C_TC^{2\zeta}$ such that the following holds
$$\sup_{t\in[0,T]}\|v_1-v_1^*\|_{L^2}<\epsilon,\ \sup_{t\in[0,T]}\|v_2-v_2^*\|_{L^2}<\epsilon$$ and $$\sup_{t\in[0,T]}\|v_1^*\|_{C^{2\zeta}}\lesssim_\epsilon\sup_{t\in[0,T]}\|v_1\|_{L^2},\ \sup_{t\in[0,T]}\|v_2^*\|_{C^{2\zeta}}\lesssim_\epsilon\sup_{t\in[0,T]}\|v_2\|_{L^2}.$$

We consider some $0<T^*\leq T$ to be determined later, using Lemma \ref{lem:est:leray} we bound $w$ in $L^2([0,T^*];H^\zeta)$
\begin{align}
\|w\|_{L^2_{[0,T^*]}H^\zeta}
 &\lesssim\|\int_0^te^{(t-s)\Delta}(v_1\otimes w)\dif s\|_{L^2_{[0,T^*]}H^{\zeta+1}}+\|\int_0^te^{(t-s)\Delta}(w\otimes v_2)\dif s\|_{L^2_{[0,T^*]}H^{\zeta+1}}\notag\\
 &\ \ \ \ \ \ \ +\|\int_0^te^{(t-s)\Delta}(w\otimes z+z\otimes w)\dif s\|_{L^2_{[0,T^*]}H^{\zeta+1}}=:\rm\uppercase\expandafter{\romannumeral1}_1+\uppercase\expandafter{\romannumeral1}_2+\uppercase\expandafter{\romannumeral1}_3.\notag
\end{align}
We bound each terms separately. For $\rm\uppercase\expandafter{\romannumeral1}_1$ by Lemma \ref{C.2} and Lemma \ref{lem:C.4} we obtain
\begin{align*}
\rm\uppercase\expandafter{\romannumeral1}_1
&\leq\|(v_1-v_1^*)\otimes w\|_{L^2_{[0,T^*]}H^{\zeta-1}}+\left(\int_0^{T^*}(\int_0^t(t-s)^{-\frac{1}{2}}\|v_1^*\otimes  w\|_{H^\zeta}\dif s)^2\dif t\right)^{\frac12}.
\end{align*}
We bound the first term on the right hand side by H\"older's inequality and the Sobolev embedding
\begin{align*}
\|(v_1-v_1^*)\otimes w\|_{L^2_{[0,T^*]}H^{\zeta-1}}&\lesssim \|(v_1-v_1^*)\otimes w\|_{L^2_{[0,T^*]}L^p}\lesssim(\int_0^{T^*}\|w\|_{L^{p'}}^2\|v_1-v_1^*\|_{L^2}^2\dif t)^{1/2}\\
   &\lesssim \sup_{t\in[0,T]}\|v_1(t)-v_1^*(t)\|_{L^2}\|w\|_{L^2_{[0,T^*]}H^{\zeta}}\leq C_1\epsilon\|w\|_{L^2_{[0,T^*]}H^{\zeta}},
\end{align*}
where we used the embeddings $L^p\hookrightarrow H^{\zeta-1}$ and $H^\zeta\hookrightarrow L^{p'}$ with $\frac{1}{p}-\frac{1}{2}=\frac{1-\zeta}{2}$, $\frac{1}{p}=\frac{1}{2}+\frac{1}{p'}$ and $C_1$ is used to denote the universal implicit constants. For the second term on the right hand side we apply H\"older’s inequality and derive for some $C_2>0$ depending on $\epsilon$
\begin{align*}
&\left(\int_0^{T^*}(\int_0^t(t-s)^{-\frac{1}{2}}\|v_1^*\|_{C^{2\zeta}}\|w\|_{H^\zeta}\dif s)^2\dif t\right)^{\frac12}\\
&\lesssim\sup_{s\in[0,T]}\|v_1^*(s)\|_{C^{2\zeta}}\left(\int_0^{T^*}\int_0^t(t-s)^{-\frac12}\dif s\int_0^t(t-s)^{-\frac12}\|w\|_{H^\zeta}^2\dif s\dif t\right)^{\frac12}\\
   &\lesssim\sup_{s\in[0,T]}\|v_1(s)\|_{L^2}\left(\int_0^{T^*}\int_s^{T^*}t^{\frac12}(t-s)^{-\frac12}\dif t\|w\|_{H^\zeta}^2\dif s\right)^{\frac12}\leq C_2T^{*\frac12}\|v_1\|_{C_TL^2}\|w\|_{L^2_{[0,T^*]}H^\zeta},
\end{align*}
which yields that
\begin{align*}
\rm\uppercase\expandafter{\romannumeral1}_1
&\leq (C_1\epsilon +C_2T^{*\frac12}\|v_1\|_{C_TL^2})\|w\|_{L^2_{[0,T^*]}H^\zeta}.
\end{align*}
Similarly we bound $\rm\uppercase\expandafter{\romannumeral1}_2$ by
\begin{align*}
\rm\uppercase\expandafter{\romannumeral1}_2
&\leq (C_1\epsilon +C_2T^{*\frac12}\|v_2\|_{C_TL^2})\|w\|_{L^2_{[0,T^*]}H^\zeta}.
\end{align*}
For $\rm\uppercase\expandafter{\romannumeral1}_3$ we obtain by Lemma \ref{lem:2.3} with $\zeta>\kappa$, Lemma \ref{C.2} and H\"older's inequality
\begin{align*}  \rm\uppercase\expandafter{\romannumeral1}_3& \lesssim\left(\int_0^{T^*}(\int_0^t(t-s)^{-\frac{2\kappa+1+\zeta}{2}}\|w\otimes z+z\otimes w\|_{H^{-2\kappa}}\dif s)^2\dif t\right)^{\frac12}\\
& \lesssim\left(\int_0^{T^*}(\int_0^t(t-s)^{-\frac{2\kappa+1+\zeta}{2}}\|w\|_{H^\zeta}\|z\|_{C^{-\kappa}}\dif s)^2\dif t\right)^{\frac12}\\
&\lesssim\|z\|_{C_TC^{-\kappa}}\left(\int_0^{T^*}\int_0^t(t-s)^{-\frac{2\kappa+1+\zeta}{2}}\dif s\int_0^t(t-s)^{-\frac{2\kappa+1+\zeta}{2}}\|w\|_{H^\zeta}^2\dif s\dif t\right)^{\frac12}\\
&\lesssim\|z\|_{C_TC^{-\kappa}}\left(\int_0^{T^*}\int_s^{T^*}t^{\frac{1-2\kappa-\zeta}{2}}(t-s)^{-\frac{2\kappa+\zeta+1}{2}}\dif t\|w\|_{H^\zeta}^2\dif s\right)^{\frac12}\\
&\leq C_2 T^{*\frac{1-2\kappa-\zeta}2}\|z\|_{C_TC^{-\kappa}}\|w\|_{L^2_{[0,T^*]}H^\zeta},\end{align*}
where we used $0<2\kappa+\zeta<1$.
Together with the above bounds we obtain
\begin{align*}
 \|w\|_{L^2_{[0,T^*]}H^\zeta}
 &\leq \left(2C_1\epsilon+C_2T^{*\frac{1-2\kappa-\zeta}2}(\|v_1\|_{C_TL^2}+\|v_2\|_{C_TL^2}+\|z\|_{C_TC^{-\kappa}})\right)\|w\|_{L^2_{[0,T^*]}H^\zeta}.
\end{align*}
Then if we choose $\epsilon=\frac{1}{4C_1}$ and $T^*\leq T$ small enough such that
$$T^{*\frac{1-2\kappa-\zeta}2}<\frac{1}{2C_2(\|v_1\|_{C_TL^2}+\|v_2\|_{C_TL^2}+\|z\|_{C_TC^{-\kappa}})+1},$$
we can deduce $ \|w\|_{L^2_{[0,T^*]}H^\zeta}=0$.
Additionally, since $w\in C([0,T];L^2)$, the conclusion that $w=0$ signifies that there is at most one possible solution to (\ref{3.2}) on $[0,T^*]$. By employing this method iteratively, starting from $T^*$, we can arrive at the desired result within finite steps, independent of the initial values. This implies that there is at most one possible solution in $C([0,T];L^2)\cap L^2([0,T];H^\zeta)$ $\mathbf{P}$-a.s. for any $T>0$, which implies our final result.
\end{proof}

\noindent{\bf Acknowledgement.} We are very grateful to Martina Hofmanov\'a for proposing this problem to us.

\appendix
\renewcommand{\appendixname}{Appendix~\Alph{section}}
\renewcommand{\theequation}{A.\arabic{equation}}

  \section{Estimates of $\rho$ and $a_{(\xi)}$}
  \label{s:appA}
  For completeness, we include here the detailed proof of the estimates (\ref{4.6}), (\ref{4.8}) employed in Section \ref{covq+12}.
We first give the estimates of $\rho$. By the definition of $\rho$ we have
$$\|\rho\|_{L^p}\leq 2l+ 2\|\mathring{R}_l\|_{L^p}+\gamma_{q+1}.$$
Furthermore, by mollification estimates, the embedding $W^{4,1}\subset L^\infty$ (in fact the time-space dimension is three), $l\leq \sigma_{q-1}$ and (\ref{3.11}) we obtain for $N\geq 0$
$$\|\mathring{R}_l\|_{C_{[ (2\sigma_{q-1})\wedge T_L,T_L],x}^N}\lesssim l^{-4-N}\|\mathring{R}_q\|_{L^1_{(\sigma_{q-1}\wedge T_L,T_L]}L^1}\lesssim l^{-4-N}\delta_{q+1}M_L.$$
Thus by $l<\delta_{q+1}M_L$ we have
\begin{align}
   \|\rho\|_{C_{{[(2\sigma_{q-1})\wedge T_L,T_L]},x}^0}\lesssim l+l^{-4}\|\mathring{R}_q\|_{L^1_{(\sigma_{q-1}\wedge T_L,T_L]}L^1}+\gamma_{q+1}\lesssim l^{-4}\delta_{q+1}M_L+\gamma_{q+1}. \notag
\end{align}
Next we estimate the $C_{t,x}^N$-norm for $N\in\mathbb{N}$. We apply the chain rule in \cite[Proposition C.1]{BDLPS15} to $f(z)=\sqrt{l^2+z^2},|D^mf(z)|\lesssim l^{-m+1}$ to obtain
\begin{align*}
\|\sqrt{l^2+|\mathring{R}_l|^2}\|_{C_{{[(2\sigma_{q-1})\wedge T_L,T_L]},x}^N}&\lesssim\|\sqrt{l^2+|\mathring{R}_l|^2}\|_{C_{{[(2\sigma_{q-1})\wedge T_L,T_L]},x}^0}+\|Df\|_{C^0}\|\mathring{R}_l\|_{C_{{[(2\sigma_{q-1})\wedge T_L,T_L]},x}^N}\\
&+\|Df\|_{C^{N-1}}\|\mathring{R}_l\|_{C_{{[(2\sigma_{q-1})\wedge T_L,T_L]},x}^1}^N\\
&\lesssim l^{-N-4}\delta_{q+1}M_L+l^{-N+1}l^{-5N}\delta_{q+1}^NM_L^N.
\end{align*}
Then by $M_L\leq l^{-1}$ we deduce for $N\geq1$
\begin{align}\label{A.1}
\|\rho\|_{C_{{[(2\sigma_{q-1})\wedge T_L,T_L]},x}^N}&\lesssim\|\sqrt{l^2+|\mathring{R}_l|^2}\|_{C_{{[(2\sigma_{q-1})\wedge T_L,T_L]},x}^N}+\gamma_{q+1}\notag\\
&\lesssim(l^{-N-4}+l^{-2N+2}l^{-5N})\delta_{q+1}M_L+\gamma_{q+1}\notag\\
&\lesssim l^{-7N+2}\delta_{q+1}M_L+\gamma_{q+1}.
\end{align}
Likewise, using the inequalities $l\leq \delta_{q+1}^{1/2}\leq -t_q$ and (\ref{3.11}), we deduce
$$\|\mathring{R}_l\|_{C_{[ 0,T_L],x}^N}\lesssim l^{-4-N}\|\mathring{R}_q\|_{L^1_{[t_q,T_L]}L^1}\lesssim l^{-4-N}(M_L+qA).$$
Then following the above calculation we obtain
\begin{align}
\|\rho\|_{C_{[0,T_L],x}^0}
\lesssim  l^{-4}(M_L+qA)+\gamma_{q+1},\label{rho0}
\end{align}
and for $N\geq1$
\begin{align}
\|\rho\|_{C_{[0,T_L],x}^N}
\lesssim  l^{-7N+2}(M_L+qA)+\gamma_{q+1},\notag
\end{align}
which implies (\ref{4.6}).

Next we estimate $a_{(\xi)}$
in $C_{t,x}^N$-norm. By Leibniz rule we get
$$\|a_{(\xi)}\|_{C_{{[(2\sigma_{q-1})\wedge T_L,T_L]},x}^N}\lesssim\sum_{m=0}^N\|\rho^{1/2}\|_{C_{{[(2\sigma_{q-1})\wedge T_L,T_L]},x}^m}\|\gamma_\xi(\Id-\frac{\mathring{R}_l}{\rho})\|_{C_{{[(2\sigma_{q-1})\wedge T_L,T_L]},x}^{N-m}}.$$
Apply \cite[Proposition C.1]{BDLPS15} to $f(z)=z^{1/2}$, $|D^mf(z)|\lesssim|z|^{1/2-m}$, for $m=1,...,N$, and using (\ref{A.1}) we obtain for $m\geq1$
\begin{align*}
\|\rho^{1/2}&\|_{C_{{[(2\sigma_{q-1})\wedge T_L,T_L]},x}^m}\notag\\
&\lesssim \|\rho^{1/2}\|_{C_{{[(2\sigma_{q-1})\wedge T_L,T_L]},x}^0}+l^{-1/2}\|\rho\|_{C_{{[(2\sigma_{q-1})\wedge T_L,T_L]},x}^m}+l^{1/2-m}\|\rho\|_{C_{{[(2\sigma_{q-1})\wedge T_L,T_L]},x}^1}^m\\
&\lesssim l^{-2}\delta_{q+1}^{1/2}M_L^{1/2}+\gamma_{q+1}^{1/2}+l^{-1/2}(l^{-7m+2}\delta_{q+1}M_L+\gamma_{q+1})+l^{1/2-m}(l^{-5m}\delta_{q+1}^mM_L^{m}+\gamma_{q+1}^m)\\
&\lesssim l^{-7m+1}(\delta_{q+1}^{1/2}M_L^{1/2}+\gamma_{q+1}^{1/2}),
\end{align*}
where we used $M_L+K\leq l^{-1}$.

Next we estimate $\gamma_\xi(\Id-\frac{\mathring{R}_l}{\rho})$, by \cite[Proposition C.1]{BDLPS15} we need to estimate
$$\|\frac{\mathring{R}_l}{\rho}\|_{C_{{[(2\sigma_{q-1})\wedge T_L,T_L]},x}^{N-m}}+\|\frac{\nabla_{t,x}\mathring{R}_l}{\rho}\|_{C_{{[(2\sigma_{q-1})\wedge T_L,T_L]},x}^0}^{N-m}+\|\frac{\mathring{R}_l}{\rho^2}\|_{C_{{[(2\sigma_{q-1})\wedge T_L,T_L]},x}^0}^{N-m}\|\rho\|_{C_{{[(2\sigma_{q-1})\wedge T_L,T_L]},x}^1}^{N-m}.$$
We use $\rho\geq l$ to have
$$\|\frac{\nabla_{t,x}\mathring{R}_l}{\rho}\|_{C_{{[(2\sigma_{q-1})\wedge T_L,T_L]},x}^0}^{N-m}\lesssim l^{-N+m}l^{-5(N-m)}\delta_{q+1}^{N-m}M_L^{N-m}\lesssim l^{-7(N-m)},$$
and in view of $|\frac{\mathring{R}_l}{\rho}|\leq 1$
$$\|\frac{\mathring{R}_l}{\rho^2}\|_{C_{{[(2\sigma_{q-1})\wedge T_L,T_L]},x}^0}^{N-m}\lesssim\|\frac{1}{\rho}\|_{C_{{[(2\sigma_{q-1})\wedge T_L,T_L]},x}^0}^{N-m}\lesssim l^{-N+m},$$
and by (\ref{A.1}) and $M_L+K\leq l^{-1}$
$$\|\rho\|_{C_{{[(2\sigma_{q-1})\wedge T_L,T_L]},x}^1}^{N-m}\lesssim l^{-5(N-m)}\delta_{q+1}^{N-m}M_L^{N-m}+\gamma_{q+1}^{N-m}\lesssim l^{-6(N-m)}.$$
Moreover, we write
\begin{align*}
\|\frac{\mathring{R}_l}{\rho}\|_{C_{{[(2\sigma_{q-1})\wedge T_L,T_L]},x}^{N-m}}&\lesssim  \sum_{k=0}^{N-m}\|\mathring{R}_l\|_{C_{{[(2\sigma_{q-1})\wedge T_L,T_L]},x}^k}\|\frac{1}{\rho}\|_{C_{{[(2\sigma_{q-1})\wedge T_L,T_L]},x}^{N-m-k}},
\end{align*}
using (\ref{A.1}) and $M_L+K\leq l^{-1}$
\begin{align*}
\|\frac{1}{\rho}\|&_{C_{{[(2\sigma_{q-1})\wedge T_L,T_L]},x}^{N-m-k}}\\
&\lesssim \|\frac{1}{\rho}\|_{C_{{[(2\sigma_{q-1})\wedge T_L,T_L]},x}^{0}}+l^{-2}\|{\rho}\|_{C_{{[(2\sigma_{q-1})\wedge T_L,T_L]},x}^{N-m-k}}+l^{-(N-m-k)-1}\|{\rho}\|_{C_{{[(2\sigma_{q-1})\wedge T_L,T_L]},x}^1}^{N-m-k}\\
&\lesssim l^{-2}(l^{1-7(N-m-k)}+\gamma_{q+1})+l^{-(N-m-k)-1}(l^{-5(N-m-k)}M_L^{N-m-k}+\gamma_{q+1}^{N-m-k})\\
&\lesssim l^{-1-7(N-m-k)}.
\end{align*}
Thus we obtain
\begin{align*}
\|\frac{\mathring{R}_l}{\rho}\|_{C_{{[(2\sigma_{q-1})\wedge T_L,T_L]},x}^{N-m}}&\lesssim  \sum_{k=0}^{N-m-1}l^{-5-k}l^{-1-7(N-m-k)}+l^{-5-(N-m)}l^{-1}\lesssim l^{-6-7(N-m)}.
\end{align*}
Finally the above bounds leads to
$$\|\gamma_\xi(\Id-\frac{\mathring{R}_l}{\rho})\|_{C_{{[(2\sigma_{q-1})\wedge T_L,T_L]},x}^{N-m}}\lesssim l^{-6-7(N-m)}.$$
Combining this with the bounds for $\rho^{1/2}$ above yields for $N\in\mathbb{N}$
\begin{align*}
\|a_{(\xi)}\|_{C_{{[(2\sigma_{q-1})\wedge T_L,T_L]},x}^N}&\lesssim (l^{-2}l^{-6-7N}+\sum_{m=1}^{N-1}l^{1-7m}l^{-6-7(N-m)}+l^{1-7N})(\delta_{q+1}^{1/2}M_L^{1/2}+\gamma_{q+1}^{1/2})\\
&\lesssim l^{-8-7N}(\delta_{q+1}^{1/2}M_L^{1/2}+\gamma_{q+1}^{1/2}),
\end{align*}
where the final bound is also valid for $N=0$.
Similarly we obtain
\begin{align*}
\|a_{(\xi)}\|_{C_{[0,T_L],x}^N}
\lesssim l^{-8-7N}((M_L+qA)^{1/2}+\gamma_{q+1}^{1/2}),
\end{align*}
which implies (\ref{4.8}).

 \renewcommand{\appendixname}{Appendix~\Alph{section}}
  \renewcommand{\theequation}{B.\arabic{equation}}

  \section{Introduction to accelerating jet flows}
  \label{s:appB}

In this part we recall the construction of accelerating jets from \cite[Section 3.1-3.3]{CL20}. We point
out that the construction is entirely deterministic, that is, none of the functions below depends on
$\omega$. Let us begin with the following geometric lemma. Recall that $\mathcal{S}^{2\times 2}$ is the set of symmetric $2\times 2$ matrics.

\bl$($\cite[Lemma 3.1]{CL20}$)$\label{lem:B.1}
Denote by $\bar{B}_{1/2}(\mathrm{Id})$ the closed ball of radius $1/2$ around the identity matrix Id, in the
space of $\mathcal{S}^{2\times 2}$. There exists $\Lambda\in \mathbb{S}^{1}\cap \mathbb{Q}^2$ such that for each $\xi\in\Lambda$ there exists
a $C^\infty$-function $\gamma_\xi$: $\bar{B}_{1/2}(\mathrm{Id})\to \mathbb{R}$ such that
$$R=\sum_{\xi\in\Lambda}\gamma_\xi^2(R)(\xi\otimes\xi)$$
for every symmetric matrix satisfying $|R-\mathrm{Id}|\leq 1/2$.
\el

We choose a collection of distinct points $p_{(\xi)} \in \mathbb{T}^2$
for $\xi\in\Lambda$ and a number $\mu_0 >2$ such
that
$$\cup_{\xi\in\Lambda}B_{\frac{2}{\mu_0}}(p_{(\xi)})\subset[0,1]^2,$$
and for $\xi\neq\xi'\in\Lambda$
$$\dif(p_{(\xi)},p_{(\xi')})>\frac{2}{\mu_0},
$$
the points $p_{(\xi)}$ will be the centers of $W_{(\xi)}$ and $W_{(\xi)}^{(c)}$.\\
For $\xi\in\Lambda$, let us introduce unit vectors $\xi^{\perp}:=(\xi_2,-\xi_1),$
and their associated coordinates
$$x_{(\xi)}=(x-p_{(\xi)})\cdot\xi,\ \ y_{(\xi)}=(x-p_{(\xi)})\cdot\xi^{\perp}\ \ for\ \ x\in \mathbb{T}^2 .$$
We first construct a family of stationary jets. To this end, we introduce the parameters:
$\nu,\mu\in\mathbb{N}$ with $\mu_0<\nu\leq \mu.$ Now we choose compactly supported nontrivial $\varphi, \psi\in C_c^\infty((-\frac{1}{\mu_0} , \frac{1}{\mu_0}
))$ and define non-periodic potentials $\tilde{\Psi}_{(\xi)}\in C_c^\infty(\mathbb{R}^2)$ and vector fields $\tilde{W}_{(\xi)},\tilde{W}_{(\xi)}^{(c)}\in C_c^{\infty}(\mathbb{R}^2)$
$$\tilde{\Psi}_{(\xi)}(x)=c_{(\xi)}\mu^{-1}(\nu\mu)^{1/2}\varphi(\nu x_{(\xi)})\psi(\mu y_{(\xi)}),$$
$$\tilde{W}_{(\xi)}(x)=-c_{(\xi)}(\nu\mu)^{1/2}\varphi(\nu x_{(\xi)})\psi'(\mu y_{(\xi)})\xi,$$
$$\tilde{W}^{(c)}_{(\xi)}(x)=c_{(\xi)}\nu\mu^{-1}(\nu\mu)^{1/2}\varphi'(\nu x_{(\xi)})\psi(\mu y_{(\xi)})\xi^{\perp},$$
where $c_{(\xi)}>0$ are normalizing constants such that
\begin{align}
\int_{\mathbb{T}^2}\!\!\!\!\!\!\!\!\!\!\!\!\; {}-{}  \tilde{W}_{(\xi)}\otimes\tilde{W}_{(\xi)}\dif x=\xi\otimes\xi.\notag
\end{align}
We periodize so that $\tilde{\Psi}_{(\xi)},\tilde{W}_{(\xi)}$ and $\tilde{W}^{(c)}_{(\xi)}$ are viewed as periodic functions on $\mathbb{T}^2$.
Finally note that each $\tilde{W}_{(\xi)}$ has disjoint support
$$\supp \tilde{W}_{(\xi)} \cap \supp \tilde{W}_{(\xi')} = \emptyset\ \ if\ \ \xi\neq\xi',$$
and for any $\xi\in\Lambda$, the following identities hold
\begin{align}
\div(\tilde{W}_{(\xi)}\otimes \tilde{W}_{(\xi)})=\xi\cdot\nabla|\tilde{W}_{(\xi)}|^2\xi,\notag
\end{align}
$$\nabla^{\perp}\tilde{\Psi}_{(\xi)}=-\partial_{y_{(\xi)}}\tilde{\Psi}_{(\xi)}\xi+\partial_{x_{(\xi)}}\tilde{\Psi}_{(\xi)}\xi^{\perp}=\tilde{W}_{(\xi)}+\tilde{W}^{(c)}_{(\xi)}.$$

Next, we introduce a simple method to avoid the collision of the
support sets of different $\tilde{W}_{(\xi)}$. Let us first choose temporal functions $g_{(\xi)}$ and $h_{(\xi)}$ to oscillate the
building blocks $\tilde{W}_{(\xi)}$ intermittently in time. Let $G \in C_c^\infty(0, 1)$ be such that
$$\int_0^1G^2(t)\dif t=1,\ \ \int_0^1G(t)\dif t=0.$$
For any $\eta \geq 1$, we define $\tilde{g}_{(\xi)}: \mathbb{T}\to \mathbb{R}$ as the 1-periodic extension of $\eta^{1/2}G(\eta(t -t_\xi))$, where $t_\xi$ are chosen so that $\tilde{g}_{(\xi)}$ have disjoint supports for different $\xi$. In other words,
$$\tilde{g}_{(\xi)}(t)=\sum_{n\in\mathbb{Z}}\eta^{1/2}G(\eta(n+t -t_\xi)).$$
We will also oscillate the velocity perturbation at a large
frequency $\sigma\in\mathbb{N}$. So we define
$$g_{(\xi)}(t)=\tilde{g}_{(\xi)}(\sigma t).$$
Then we have
\begin{align}\label{gwnp}
    \|g_{(\xi)}\|_{W^{n,p}([0,1])}\lesssim(\sigma\eta)^n\eta^{1/2-1/p}.
\end{align}
For the corrector term we define $h_{(\xi)}:\mathbb{T}\to\mathbb{R}$ by
\begin{align}
h_{(\xi)}(t)=\int_0^{\sigma t}(\tilde{g}^2_{(\xi)}(s)-1)\dif s.\notag
\end{align}

In view of the zero-mean condition for ${\tilde{g}}^2_{(\xi)}(t)-1$, these $h_{(\xi)}$ are $\mathbb{T}/\sigma$-periodic  and we have
\begin{align}
    \|h_{(\xi)}\|_{L^\infty}\leq1.\label{hlinfty}
\end{align}
Moreover, we have the identity
\begin{align}
\partial_t(\sigma^{-1}h_{(\xi)})=g^2_{(\xi)}(t)-1.    \label{parth}
\end{align}

\br\label{rem:gxi}
We emphasize that  $g_{(\xi)}$ is $\mathbb{T}/\sigma$-periodic and thus for any $n,n_1<n_2\in\mathbb{N}_0,p\in[1,\infty)$ we have
$$\|g_{(\xi)}\|_{ W^{n,p}[n_1/\sigma,n_2/\sigma]}^p=\frac{n_2-n_1}\sigma\|g_{(\xi)}\|_{ W^{n,p}[0,1]}^p.$$
When dealing with arbitrary values of $0<a_0<b_0$, it is possible to find  $n_1,n_2\in\mathbb{N}_0$ satisfying $\frac{n_1}{\sigma}<a_0\leq\frac{n_1+1}{\sigma},\frac{n_2}{\sigma}\leq b_0<\frac{n_2+1}{\sigma}$. Then we have
\begin{align}
\|g_{(\xi)}\|_{W^{n,p}[a_0,b_0]}^p\leq \|g_{(\xi)}\|_{ W^{n,p}[\frac{n_1}\sigma,\frac{n_2+1}{\sigma}]}^p=\frac{n_2-n_1+1}\sigma\|g_{(\xi)}\|_{ W^{n,p}[0,1]}^p\leq (b_0-a_0+\frac{2}{\sigma})\|g_{(\xi)}\|_{ W^{n,p}[0,1]}^p.\notag
\end{align}
\er
Now we will let the stationary flows $\tilde{W}_{(\xi)}$ travel along $\mathbb{T}^2$ in time, relating the velocity of the
moving support sets to the intermittent oscillator $g_{(\xi)}$. More precisely, we define
\begin{align}\label{Wxi}
W_{(\xi)}(x,t):&=\tilde{W}_{(\xi)}(\sigma x+\phi_{(\xi)}(t)\xi),\notag\\
W_{(\xi)}^{(c)}(x,t):&=\tilde{W}_{(\xi)}^{(c)}(\sigma x+\phi_{(\xi)}(t)\xi),\\
\Psi_{(\xi)}(x,t):&=\tilde{\Psi}_{(\xi)}(\sigma x+\phi_{(\xi)}(t)\xi),\notag
\end{align}
where $\phi_{(\xi)}(t)$ is defined by
$$\phi'_{(\xi)}=\theta g_{(\xi)}(t)$$
for some $\theta\in\mathbb{N}$. Hence we have
\begin{align}
\int_{\mathbb{T}^2}\!\!\!\!\!\!\!\!\!\!\!\!\; {}-{}  {W}_{(\xi)}\otimes{W}_{(\xi)}\dif x=\xi\otimes\xi,\label{wpingjun}
\end{align}
and each $g_{(\xi)}{W}_{(\xi)}$ has disjoint support
\begin{align}
    \supp (g_{(\xi)}{W}_{(\xi)}) \cap \supp (g_{(\xi')}{W}_{(\xi')}) = \emptyset\ \ if\ \ \xi\neq\xi'.\label{wbujiao}
\end{align}

Then, we claim that for $N\geq0,p\in [1, \infty]$ the
following holds
\begin{align}
\|\nabla^NW_{(\xi)}\|_{C_tL^p}+\nu^{-1}\mu\|\nabla^NW^{(c)}_{(\xi)}\|_{C_tL^p}+\mu\|\nabla^N\Psi_{(\xi)}\|_{C_tL^p}\lesssim(\sigma\mu)^N (\nu\mu)^{1/2-1/p},\label{2.4}
\end{align}
\begin{align}
\|(\xi\cdot\nabla)\Psi_{(\xi)}\|_{C_tL^p}\lesssim\mu^{-1}\sigma\nu (\nu\mu)^{1/2-1/p},\label{2.4*}
\end{align}
where the implicit constants may depend on $N,p$, but are independent of $\nu,\sigma,\theta,\eta $ and $\mu$. These estimates can be easily deduced from the definitions.

Moreover we have
\begin{align}
\sigma^{-1}\nabla^{\perp}\Psi_{(\xi)}&=W_{(\xi)}+W^{(c)}_{(\xi)},\label{2.6}\\
\partial_t|W_{(\xi)}|^2\xi&=\sigma^{-1}\theta g_{(\xi)}\div(W_{(\xi)}\otimes W_{(\xi)}),\label{B.9}\\
\partial_t\Psi_{(\xi)}&=\sigma^{-1}\theta g_{(\xi)}(\xi\cdot\nabla)\Psi_{(\xi)}.\label{B.10}
\end{align}

 \renewcommand{\appendixname}{Appendix~\Alph{section}}
  \renewcommand{\theequation}{C.\arabic{equation}}
  \section{Some technical tools}
  \label{s:appC}
\subsection{Improved H\"older’s inequality on $\mathbb{T}^d$}\label{ihiot2}This theorem allows us to quantify the decorrelation in the usual Hölder’s
inequality when we increase the oscillation of one function. The proof follows from similar argument as \cite[Lemma 2.1]{MS18}.
\bt\label{ihiot}
Let $p \in [1, \infty],m<n\in\mathbb{N}_0$ and $a\in C^1(\mathbb{R}^d; \mathbb{R}), f\in L^p(\mathbb{T}^d ; \mathbb{R})$. Then for any $\sigma\in\mathbb{N}$,
$$| \|af(\sigma\cdot)\|_{L^p([\frac{m}{\sigma},\frac{n}{\sigma}]^d)}- \|a\|_{L^p([\frac{m}{\sigma},\frac{n}{\sigma}]^d)}\|f\|_{L^p(\mathbb{T}^d)}|\lesssim \sigma^{-1/p}(\frac{n-m}{\sigma})^{d/p}\|a\|_{C^{0,1}([\frac{m}{\sigma},\frac{n}{\sigma}]^d)}\|f\|_{L^p(\mathbb{T}^d)}.$$
Here $C^{0,1}$ is the Lipschitz norm.
In particular, if $d=1$ and $f$ is mean-zero we have
$$|\int_{\frac{m}{\sigma}}^{\frac{n}{\sigma}}a(t)f(\sigma t)\dif t|\lesssim\frac{n-m}{\sigma^2}\|a\|_{C^{0,1}([\frac{m}{\sigma},\frac{n}{\sigma}])}\|f\|_{L^1
(\mathbb{T})}.$$
\et
\begin{proof}
To begin, we partition the domain $[\frac{m}{\sigma},\frac{n}{\sigma}]^d$ into $(n-m)^d$ smaller cubes, denoted by $Q_{i},i=1,...,(n-m)^d$. Each of these cubes has an edge length of $\frac{1}{\sigma}$.
\begin{align*}
 \int_{[\frac{m}{\sigma},\frac{n}{\sigma}]^d}|a(x)f(\sigma x)|^p\dif x&=\sum_{k=1}^{(n-m)^d} \int_{Q_k}  \left(|a(x)|^p-\sigma^d\int_{Q_k} |a(y)|^p\dif y \right)|f(\sigma x)|^p\dif x\\
 &\ \ \ \ +\sum_{k=1}^{(n-m)^d}\sigma^d \int_{Q_k}|a(y)|^p\dif y \int_{Q_k}|f(\sigma x)|^p\dif x\\
 &=\sum_{k=1}^{(n-m)^d}  \sigma^d\int_{Q_k} \left(\int_{Q_k} (|a(x)|^p-|a(y)|^p)\dif y \right)|f(\sigma x)|^p\dif x\\
 &\ \ \ \ +\sum_{k=1}^{(n-m)^d} \int_{Q_k}|a(y)|^p\dif y \int_{\mathbb{T}^d}|f(x)|^p\dif x.
\end{align*}
Thus we have
\begin{align*}
    |\|af(\sigma\cdot)\|_{L^p([\frac{m}{\sigma},\frac{n}{\sigma}]^d)}^p&- \|a\|^p_{L^p([\frac{m}{\sigma},\frac{n}{\sigma}]^d)}\|f\|^p_{L^p(\mathbb{T}^d)}|\\
    &\leq \left|\sum_{k=1}^{(n-m)^d}  \sigma^d\int_{Q_k} \left(\int_{Q_k} (|a(x)|^p-|a(y)|^p)\dif y \right)|f(\sigma x)|^p\dif x\right|\\
    &\lesssim\sum_{k=1}^{(n-m)^d}  \int_{Q_k} \frac{1}{\sigma}\|a\|_{C^{0,1}([\frac{m}{\sigma},\frac{n}{\sigma}]^d)}^{p} )|f(\sigma x)|^p\dif x\\
&\lesssim  \frac{1}{\sigma}\|a\|_{C^{0,1}([\frac{m}{\sigma},\frac{n}{\sigma}]^d)}^{p}(\frac{n-m}{\sigma})^d\|f\|_{L^p(\mathbb{T}^d)}^p,
\end{align*}
where in the second inequality we used the fact that for any $x,y\in Q_k$
\begin{align*}
   \left| |a(x)|^p-|a(y)|^p\right|\lesssim \frac{1}{\sigma}\|a\|_{C^0([\frac{m}{\sigma},\frac{n}{\sigma}]^d)}^{p-1}\| a\|_{C^{0,1}([\frac{m}{\sigma},\frac{n}{\sigma}]^d)}\lesssim \frac{1}{\sigma}\|a\|_{C^{0,1}([\frac{m}{\sigma},\frac{n}{\sigma}]^d)}^{p}.
\end{align*}
Then using the fact that for $A,B>0$ $|A-B|^p\leq\big{|} |A|^p-|B|^p\big{|}$ we obtain the first result.

For the second result, we have the same decomposition as above:
\begin{align*}
 \int_{\frac{m}{\sigma}}^{\frac{n}{\sigma}}a(t)f(\sigma t)\dif t
 &=\sigma\sum_{k=m+1}^n \int_{\frac{k-1}{\sigma}}^{\frac{k}{\sigma}} \int_{\frac{k-1}{\sigma}}^{\frac{k}{\sigma}} [a(t)- a(s)]f(\sigma t)\dif s \dif t+\sum_{k=m+1}^n  \sigma\int_{\frac{k-1}{\sigma}}^{\frac{k}{\sigma}} a(s)\dif s \int_{\frac{k-1}{\sigma}}^{\frac{k}{\sigma}} f(\sigma t)\dif t\\
  &=\sigma\sum_{k=m+1}^n \int_{\frac{k-1}{\sigma}}^{\frac{k}{\sigma}} \int_{\frac{k-1}{\sigma}}^{\frac{k}{\sigma}} [a(t)- a(s)]f(\sigma t)\dif s \dif t
\end{align*}
where we used the fact that $f$ is mean-zero. The following estimate is same as above.
\end{proof}

\subsection{Estimates for the heat operator}\label{eftho}
In this section we give the basic estimates for the heat semigroup $P_t:=e^{t\Delta}$. Let $T\geq1$.

\bl$($\cite[Lemma 2.8]{ZZZ20}$)$\label{C.2}
For any $\theta > 0$, $\alpha \in \mathbb{R},\ p, q \in [1, \infty]$ and all $t \in (0, T ]$,
$$ \|P_tf\|_{ B^{\theta+\alpha}_{p,q}} \lesssim T^{\theta/2}t^{-\theta/2}\| f\|_{ B^\alpha_{p,q}},\ \ \|P_tf\|_{ L^p} \lesssim T^{\theta/2}t^{-\theta/2}\| f\|_{ B^{-\theta}_{p,\infty}}.$$
\el

Let $\mathcal{I}f = \int_0^\cdot  P_{t-s}f\dif s.$ Then we have
\bl$($\cite[Lemma 2.9]{ZZZ20}$)$\label{C.3}
Let $\alpha \in(0, 2), p\in[1,\infty]$.
Then
$$\| \mathcal{I}f\|_{C_TB^\alpha_{p,\infty}}+\|\mathcal{I}f\|_{C_T^{\frac\alpha2}L^p} \lesssim T\| f\|_{C_TB^{\alpha-2}_{p,\infty}}.$$
\el

\bl\label{lem:C.4}$($\cite[Lemma A.5]{SZZ22}$)$
Let $\beta\in\mathbb{R}$ we have
\begin{align*}
    \|\mathcal{I}f\|_{L^2_{[0,T]}H^\beta}\lesssim\|f\|_{L^2_{[0,T]}H^{\beta-2}}.
\end{align*}
\el

\end{document}